# A catalog of matchstick graphs

Raffaele Salvia

**Abstract**

Classification of planar unit-distance graphs with up to 9 edges, by homeomorphism and isomorphism classes. With exactly nine edges, there are 633 nonisomorphic connected matchstick graphs, of which 19 are topologically distinct from each other. Increasing edges' number, their quantities rise more than exponentially, in a still unclear way.

## 1. Introduction

Two remarkable classes of graphs are planar graphs and unit distance graphs. Their intersection, that is the class of the planar graphs drawable using edges of unitary length, has been defined as the "matchstick graphs" class[1].

It is an NP-hard problem to test if a given graph can be represented as a matchstick graph (Eades, Wormland 1990 [2]; Cabello et al. 2007 [3]), and no formulas are known to determine how many distinct matchstick graphs exist for a given number of edges.

The problem of finding the number of topologically distinct matchstick graph with $n$ edges is known as the "match problem" (Gardner 1969 [4]). R. Read (1968 [5]) solved the problem for $n=8$, reporting also the total number of matchstick graphs –up to isomorphism- with 7 and 8 edges.

Here are listed all the distinct matchstick graphs with less of 10 edges, grouped by topological equivalence.

## 2. Overview

Table 1 summarizes the number of non-homeomorphic and non-isomorphic connected planar unit distance graphs. For a quicker referring, henceforth will be assumed the following definitions:

**Definition 2.1.** Given a nonnegative integer $n$, $\varrho(n)$ indicates the number of distinct connected matchstick graph, up to topological homeomorphism, with $n$ edges.

**Definition 2.2.** Given a nonnegative integer $n$, $\MYRUNE(n)$ indicates the number of distinct connected matchstick graph, up to isomorphism, with $n$ edges.

| $n$ | $\varrho(n)$ | $\MYRUNE(n)$ |
|---|---|---|
| 1 | 1 | 1 |
| 2 | 1 | 1 |
| 3 | 3 | 3 |
| 4 | 5 | 5 |
| 5 | 10 | 12 |
| 6 | 19 | 28 |
| 7 | 39 | 74 |
| 8 | 84 | 207 |
| 9 | 196 | 633 |

Table 1.

So far, the first 8 terms of both sequences have been known [6]. Here are first showed the values of $\varrho(9)$ and $\MYRUNE(9)$. The evolution of both sequences shall be discussed in section 4.

In table 2, topologically different matchstick graphs are reported by number of edges and number of faces (including the outer region).

In the next section, every one of the 964 distinct matchstick graph with $|E|<10$ will be showed. They are listed and progressively numbered, in order of priority, by:

1. Number of edges $|E|$;
2. Number of faces $F$;
3. Maximum degree $\Delta$.

|  | F | | | | | |
|---|---|---|---|---|---|---|
| $\|E\|$ | 1 | 2 | 3 | 4 | 5 | 6 |
| 1 | 1 | 0 | 0 | 0 | 0 | 0 |
| 2 | 1 | 0 | 0 | 0 | 0 | 0 |
| 3 | 2 | 1 | 0 | 0 | 0 | 0 |
| 4 | 3 | 2 | 0 | 0 | 0 | 0 |
| 5 | 5 | 4 | 1 | 0 | 0 | 0 |
| 6 | 7 | 8 | 4 | 0 | 0 | 0 |
| 7 | 11 | 15 | 12 | 1 | 0 | 0 |
| 8 | 16 | 29 | 31 | 8 | 0 | 0 |
| 9 | 26 | 56 | 75 | 37 | 3 | 0 |

Table 2.

# 3. Classification

Drawings in the left cells of tables are aimed to illustrate the topological homeomorphism between graphs in the same box. The number on the bottom right of right cells indicates the number of graphs contained therein. Graphs 8-81-1, 9-84-1, 9-85-1, 9-179-1, 9-179-2 and 9-190-1 are represented magnified by a factor of 1.5 respect to others, and graphs 9-159-4 and 9-164-1 by a factor of 2.5, to avoid apparent overlaps.

## 3.1. Matchstick graphs with $|E|=1$
### 3.1.1. Matchstick graphs with $|E|=1$, $F=1$, $|V|=2$
#### 3.1.1.1. Matchstick graphs with $|E|=1$, $F=1$, $|V|=2$, $\Delta=1$

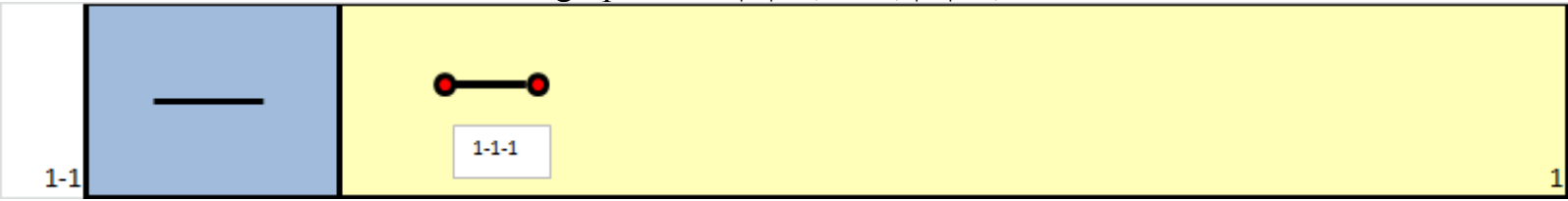

## 3.2. Matchstick graphs with $|E|=2$
### 3.2.1. Matchstick graphs with $|E|=2$, $F=1$, $|V|=3$
#### 3.2.1.1. Matchstick graphs with $|E|=2$, $F=1$, $|V|=3$, $\Delta=2$

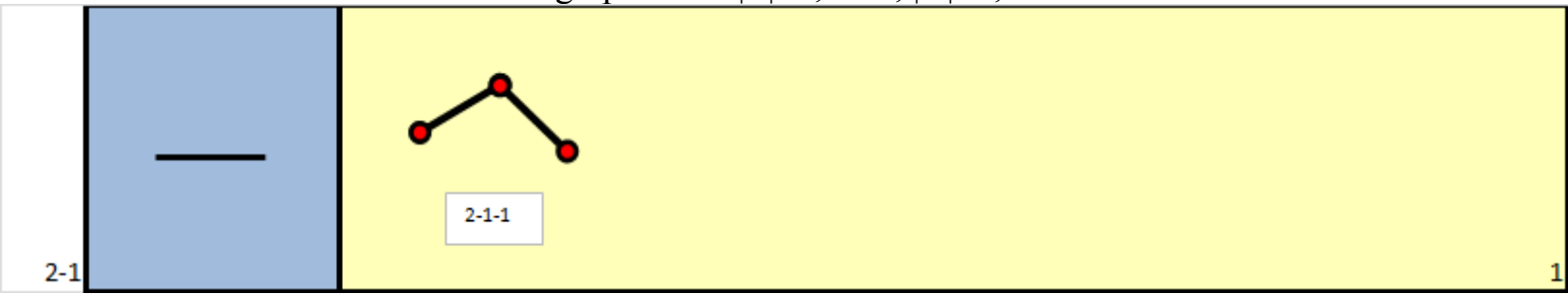

## 3.3. Matchstick graphs with $|E|=3$
### 3.3.1. Matchstick graphs with $|E|=3$, $F=1$, $|V|=4$
#### 3.3.1.1. Matchstick graphs with $|E|=3$, $F=1$, $|V|=4$, $\Delta=2$

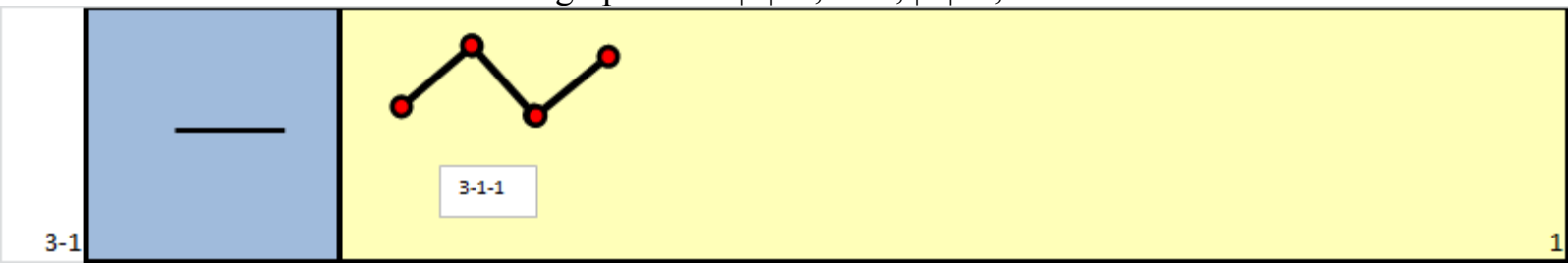

#### 3.3.1.2. Matchstick graphs with $|E|=3$, $F=1$, $|V|=4$, $\Delta=3$

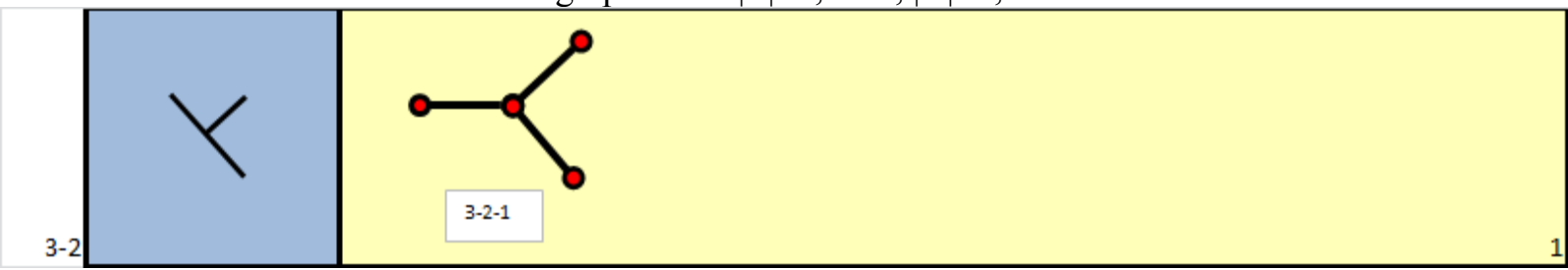

### 3.3.2. Matchstick graphs with $|E|=3$, $F=2$, $|V|=3$
#### 3.3.2.1. Matchstick graphs with $|E|=3$, $F=2$, $|V|=3$, $\Delta=2$

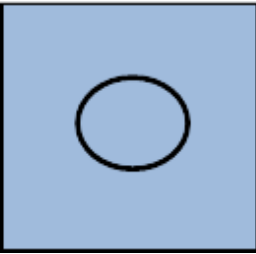 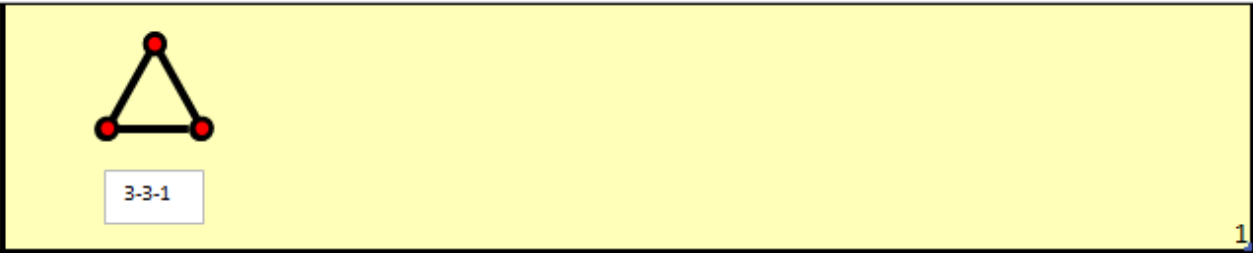

## 3.4. Matchstick graphs with |E|=4
### 3.4.1. Matchstick graphs with |E|=4, F=1, |V|=5
#### 3.4.1.1. Matchstick graphs with |E|=4, F=1, |V|=5, Δ=2

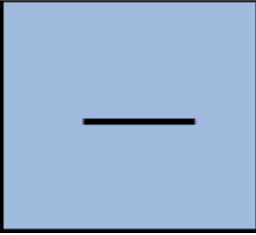 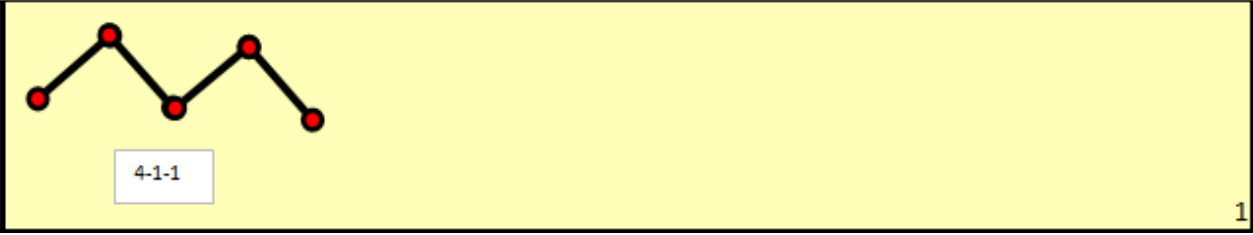

#### 3.4.1.2. Matchstick graphs with |E|=4, F=1, |V|=5, Δ=3

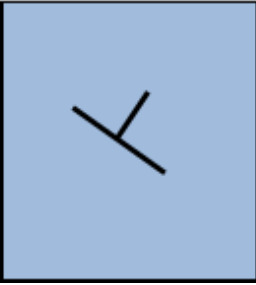 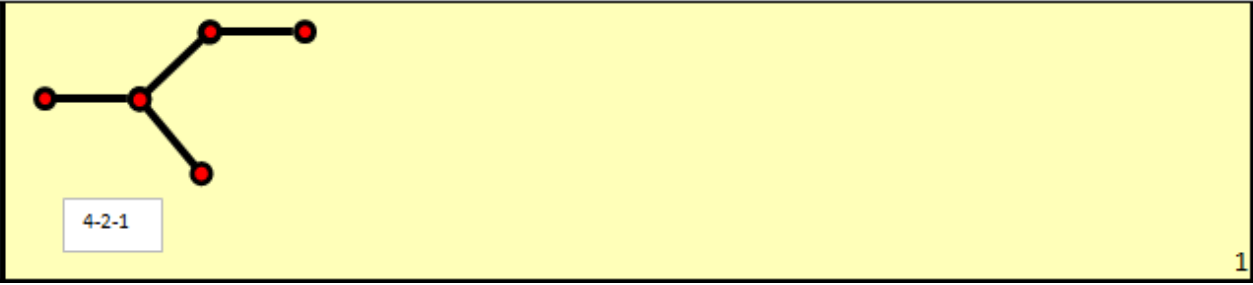

#### 3.4.1.3. Matchstick graphs with |E|=4, F=1, |V|=5, Δ=4

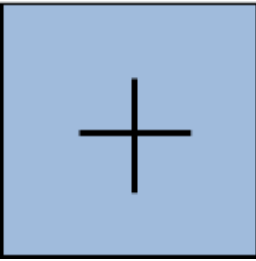 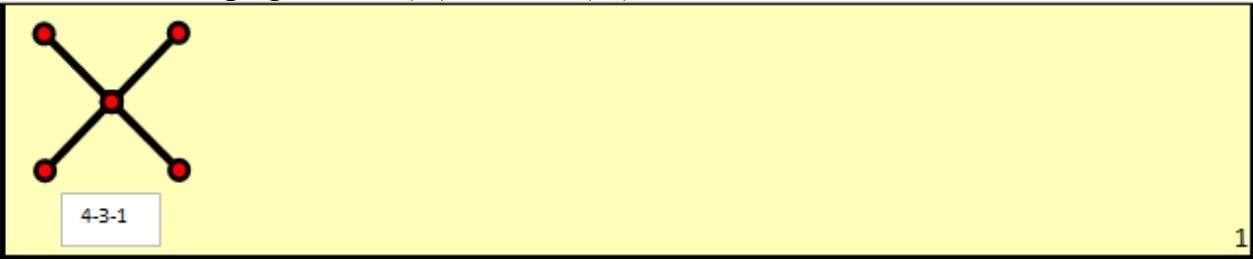

### 3.4.2. Matchstick graphs with |E|=4, F=2, |V|=4
#### 3.4.2.1. Matchstick graphs with |E|=4, F=2, |V|=4, Δ=2

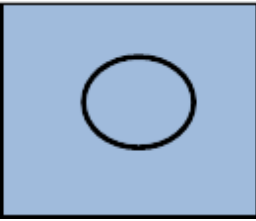 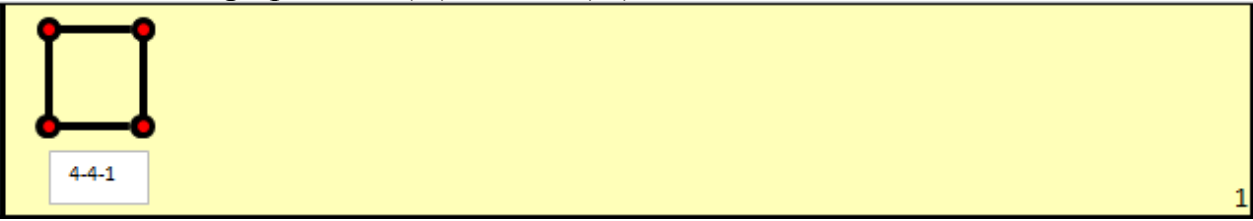

#### 3.4.2.2. Matchstick graphs with |E|=4, F=2, |V|=4, Δ=3

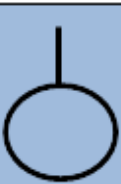
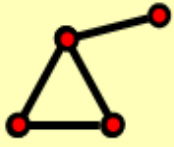

4-5-1

4-5

## 3.5. Matchstick graphs with |E|=5
### 3.5.1. Matchstick graphs with |E|=5, $F$=1, |V|=6
#### 3.5.1.1. Matchstick graphs with |E|=5, $F$=1, |V|=6, $\Delta$=2

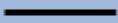
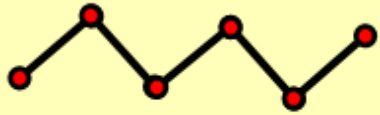

5-1-1

5-1

#### 3.5.1.2. Matchstick graphs with |E|=5, $F$=1, |V|=6, $\Delta$=3

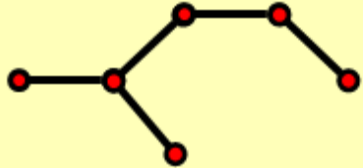
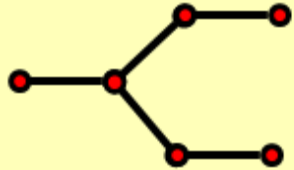

5-2-1    5-2-2

5-2

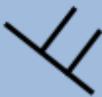
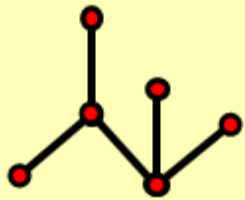

5-3-1

5-3

#### 3.5.1.3. Matchstick graphs with |E|=5, $F$=1, |V|=6, $\Delta$=4

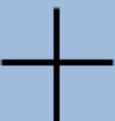
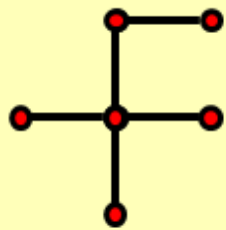

5-4-1

5-4

#### 3.5.1.3. Matchstick graphs with |E|=5, $F$=1, |V|=6, $\Delta$=5

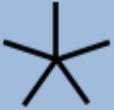
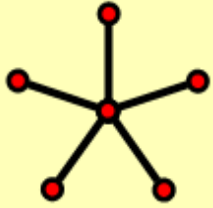

5-5-1

5-5

### 3.5.2. Matchstick graphs with |E|=5, F=2, |V|=5
#### 3.5.2.1. Matchstick graphs with |E|=5, F=2, |V|=5, Δ=2

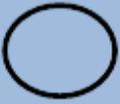
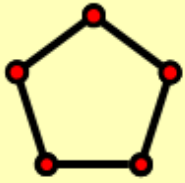

5-6-1

5-6

#### 3.5.2.2. Matchstick graphs with |E|=5, F=2, |V|=5, Δ=3

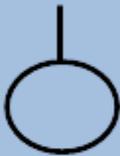
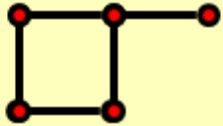
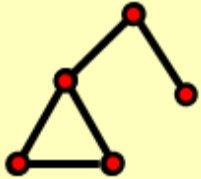

5-7-1    5-7-2

5-7

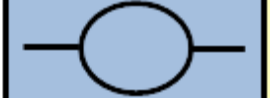
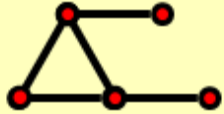

5-8-1

5-8

#### 3.5.2.3. Matchstick graphs with |E|=5, F=2, |V|=5, Δ=4

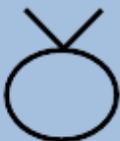
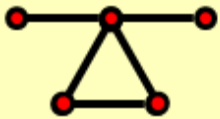

5-9-1

5-9

### 3.5.3. Matchstick graphs with |E|=5, F=3, |V|=4
#### 3.5.2.1. Matchstick graphs with |E|=5, F=3, |V|=4, Δ=3

| 5-10 | 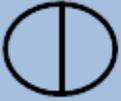 | 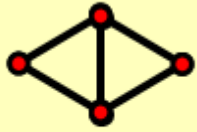 |
| --- | --- | --- |
| | | 5-10-1 |

## 3.6. Matchstick graphs with |E|=6
### 3.6.1. Matchstick graphs with |E|=6, $F$=1, |V|=7
#### 3.6.1.1. Matchstick graphs with |E|=6, $F$=1, |V|=7, $\Delta$=2

| 6-1 | 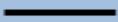 | 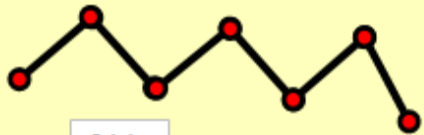 |
| --- | --- | --- |
| | | 6-1-1 |

#### 3.6.1.2. Matchstick graphs with |E|=6, $F$=1, |V|=7, $\Delta$=3

| 6-2 | 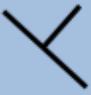 | 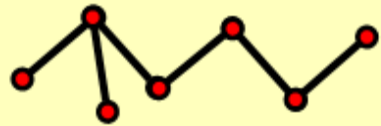 | 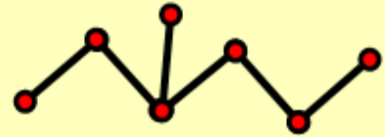 | 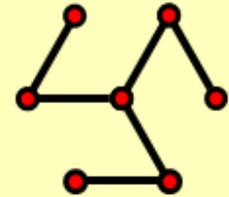 |
| --- | --- | --- | --- | --- |
| | | 6-2-1 | 6-2-2 | 6-2-3 |

| 6-3 | 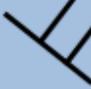 | 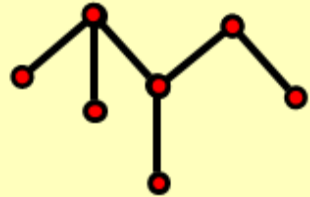 | 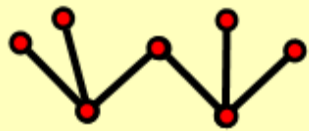 |
| --- | --- | --- | --- |
| | | 6-3-1 | 6-3-2 |

#### 3.6.1.3. Matchstick graphs with |E|=6, $F$=1, |V|=7, $\Delta$=4

| 6-4 | 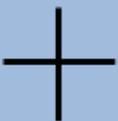 | 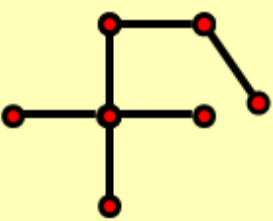 | 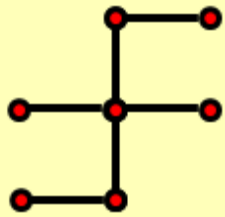 |
| --- | --- | --- | --- |
| | | 6-4-1 | 6-4-2 |

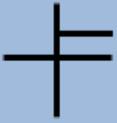
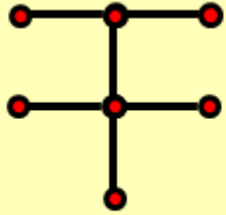

6-5-1

6-5

### 3.6.1.4. Matchstick graphs with $|E|=6$, $F=1$, $|V|=7$, $\Delta=5$

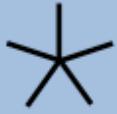
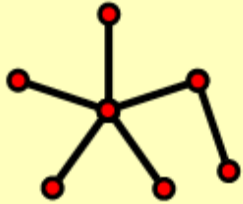

6-6-1

6-6

### 3.6.1.5. Matchstick graphs with $|E|=6$, $F=1$, $|V|=7$, $\Delta=6$

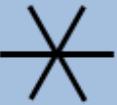
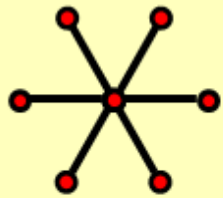

6-7-1

6-7

## 3.6.2. Matchstick graphs with $|E|=6$, $F=2$, $|V|=6$
### 3.6.2.1. Matchstick graphs with $|E|=6$, $F=2$, $|V|=6$, $\Delta=2$

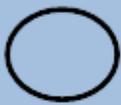
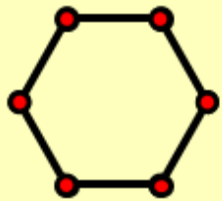

6-8-1

6-8

### 3.6.2.2. Matchstick graphs with $|E|=6$, $F=2$, $|V|=6$, $\Delta=3$

| 6-9 | | 6-9-1 | 6-9-2 | 6-9-3 | 3 |
| 6-10 | | 6-10-1 | 6-10-2 | 6-10-3 | 3 |
| 6-11 | | 6-11-1 | | | 1 |
| 6-12 | | 6-12-1 | | | 1 |

3.6.2.3. Matchstick graphs with $|E|=6$, $F=2$, $|V|=6$, $\Delta=4$

| 6-13 | | 6-13-1 | 6-13-2 | 2 |
| 6-14 | | 6-14-1 | | 1 |

3.6.2.4. Matchstick graphs with $|E|=6$, $F=2$, $|V|=6$, $\Delta=5$

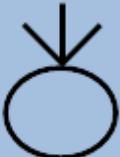

6-15

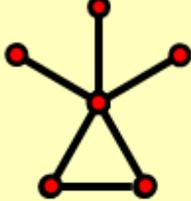

6-15-1

### 3.6.3. Matchstick graphs with |E|=6, $F$=3, |V|=5

#### 3.6.3.1. Matchstick graphs with |E|=6, $F$=3, |V|=5, $\Delta$=3

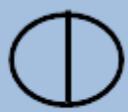

6-16

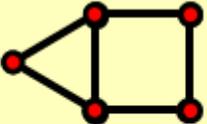

6-16-1

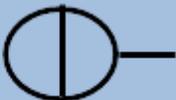

6-17

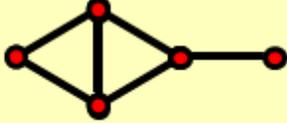

6-17-1

#### 3.6.3.2. Matchstick graphs with |E|=6, $F$=3, |V|=5, $\Delta$=4

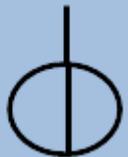

6-18

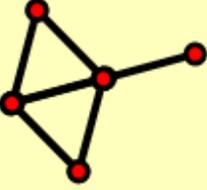

6-18-1

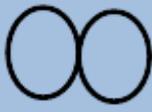

6-19

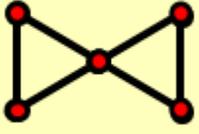

6-19-1

### 3.7. Matchstick graphs with |E|=7

#### 3.7.1. Matchstick graphs with |E|=7, $F$=1, |V|=8

##### 3.7.1.1. Matchstick graphs with |E|=7, $F$=1, |V|=8, $\Delta$=2

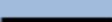

7-1

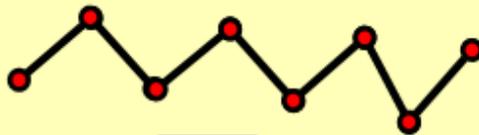

7-1-1

### 3.7.1.2. Matchstick graphs with |E|=7, F=1, |V|=8, Δ=3

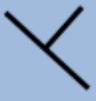

### 3.7.1.3. Matchstick graphs with |E|=7, F=1, |V|=8, Δ=4

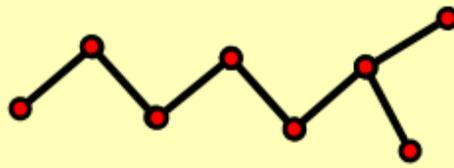

| 7-6 | 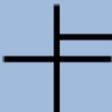 | 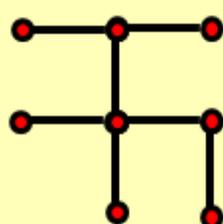 7-6-1     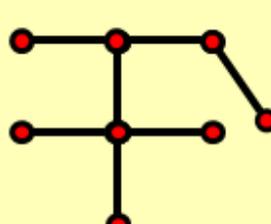 7-6-2     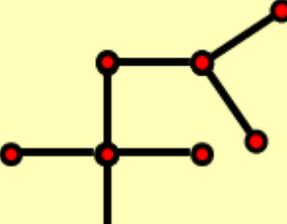 7-6-3 | 3 |

| 7-7 | 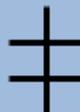 | 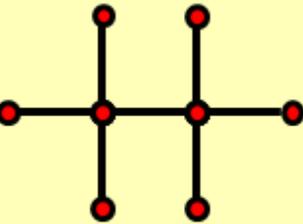 7-7-1 | 1 |

### 3.7.1.4. Matchstick graphs with $|E|=7$, $F=1$, $|V|=8$, $\Delta=5$

| 7-8 | 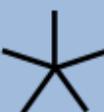 | 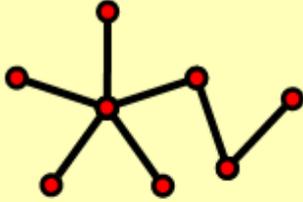 7-8-1     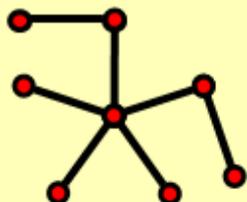 7-8-2 | 2 |

| 7-9 | 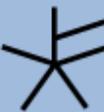 | 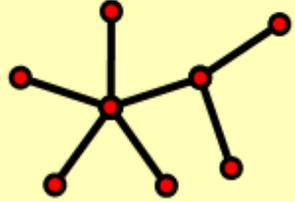 7-9-1 | 1 |

### 3.7.1.5. Matchstick graphs with $|E|=7$, $F=1$, $|V|=8$, $\Delta=6$

| 7-10 | 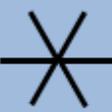 | 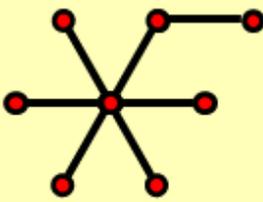 7-10-1 | 1 |

### 3.7.1.6. Matchstick graphs with $|E|=7$, $F=1$, $|V|=8$, $\Delta=7$

| 7-11 | 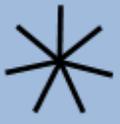 | 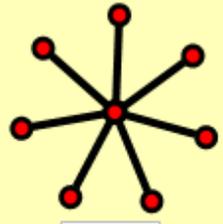 7-11-1 | 1 |

## 3.7.2. Matchstick graphs with |E|=7, F=2, |V|=7
### 3.7.2.1. Matchstick graphs with |E|=7, F=2, |V|=7, Δ=2

| 7-12 | 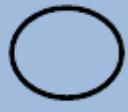 | 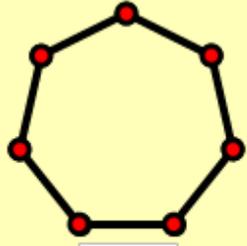 7-12-1 | 1 |

### 3.7.2.2. Matchstick graphs with |E|=7, F=2, |V|=7, Δ=3

| 7-13 | 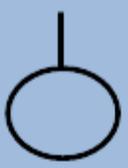 | 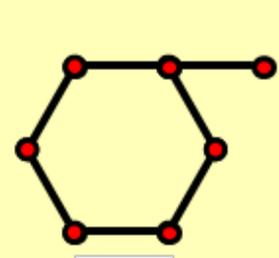 7-13-1 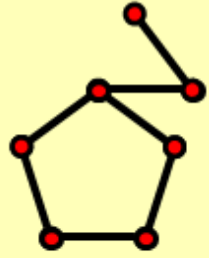 7-13-2 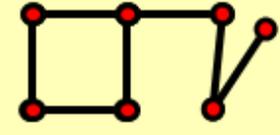 7-13-3 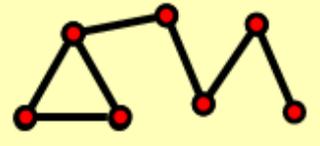 7-13-4 | 4 |
| 7-14 | 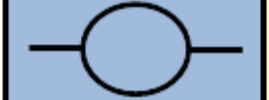 | 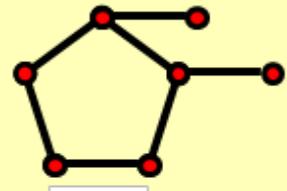 7-14-1 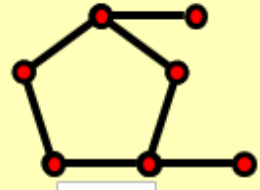 7-14-2 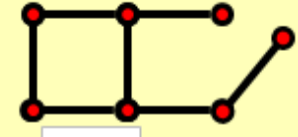 7-14-3 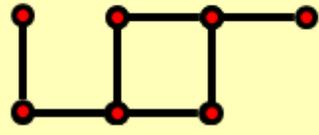 7-14-4 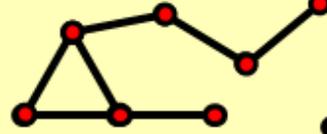 7-14-5 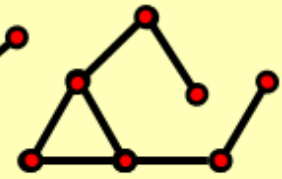 7-14-6 | 6 |
| 7-15 | 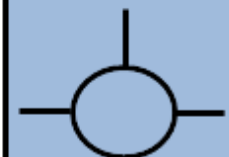 | 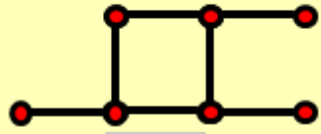 7-15-1 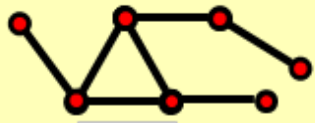 7-15-2 | 2 |

| 7-16 | 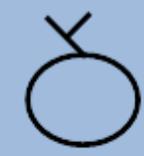 | 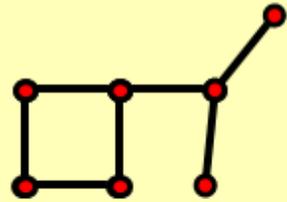 7-16-1    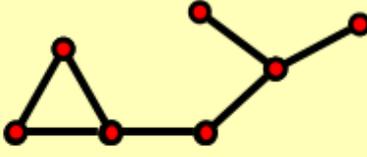 7-16-2    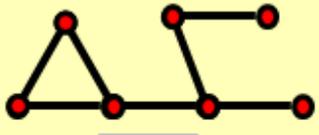 7-16-3 | 3 |
|---|---|---|---|
| 7-17 | 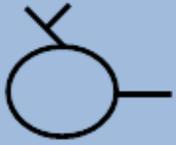 | 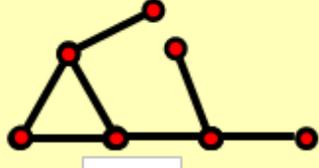 7-17-1 | 1 |

### 3.7.2.3. Matchstick graphs with $|E|=7$, $F=2$, $|V|=7$, $\Delta=4$

| 7-18 | 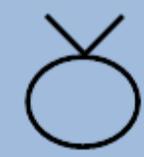 | 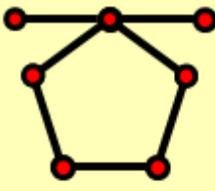 7-18-1   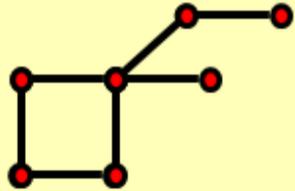 7-18-2   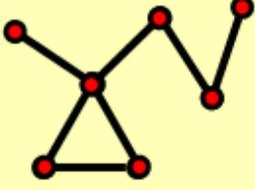 7-18-3   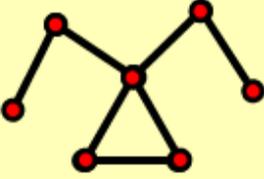 7-18-4 | 4 |
|---|---|---|---|
| 7-19 | 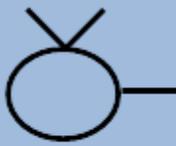 | 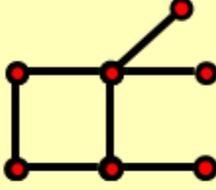 7-19-1   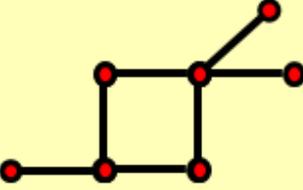 7-19-2   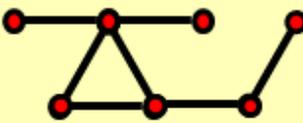 7-19-3 <br> 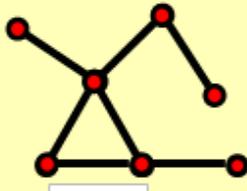 7-19-4 | 4 |
| 7-20 | 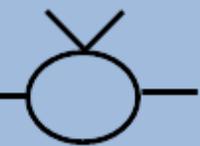 | 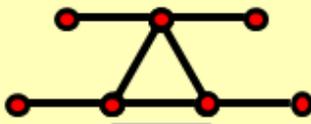 7-20-1 | 1 |
| 7-21 | 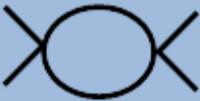 | 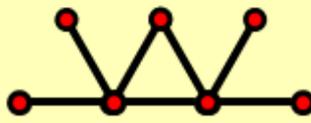 7-21-1 | 1 |

| 7-22 | 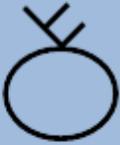 | 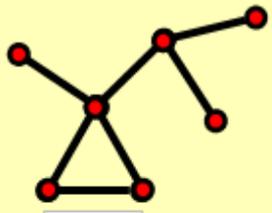 7-22-1 | 1 |
| 7-23 | 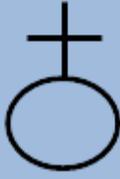 | 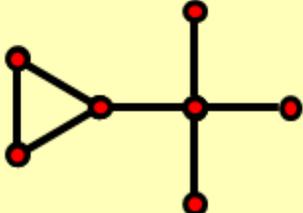 7-23-1 | 1 |

#### 3.7.2.4. Matchstick graphs with |E|=7, F=2, |V|=7, Δ=5

| 7-24 | 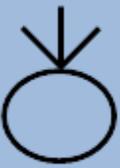 | 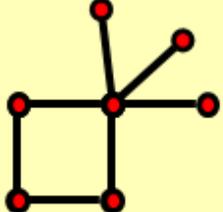 7-24-1    7-24-2 | 2 |
| 7-25 | 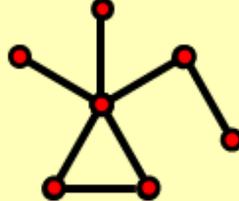 | 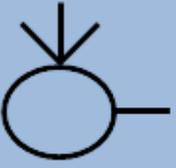 7-25-1 | 1 |

#### 3.7.2.5. Matchstick graphs with |E|=7, F=2, |V|=7, Δ=6

| 7-26 | 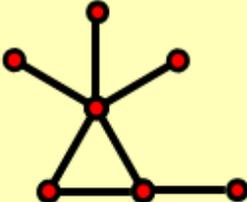 | 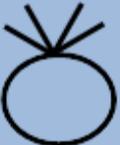 7-26-1 | 1 |

### 3.7.3. Matchstick graphs with |E|=7, F=3, |V|=6
#### 3.7.3.1. Matchstick graphs with |E|=7, F=3, |V|=6, Δ=3

| 7-27 | 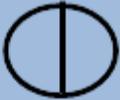 | 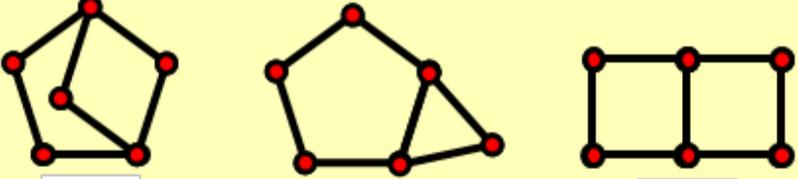 7-27-1    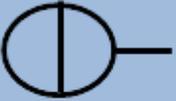 7-27-2    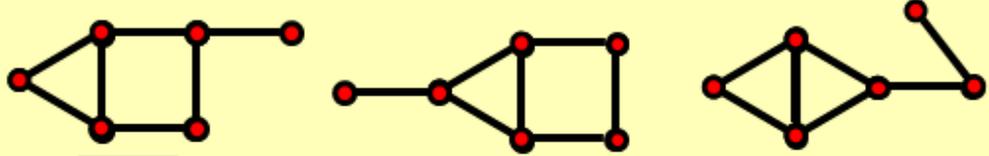 7-27-3 | 3 |

| 7-28 | 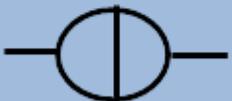 | 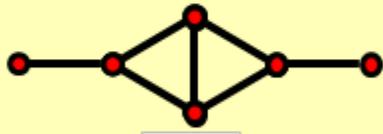 7-28-1    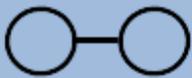 7-28-2    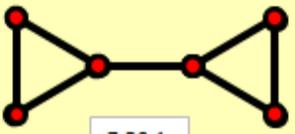 7-28-3 | 3 |

| 7-29 | 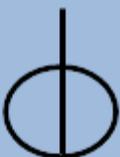 | 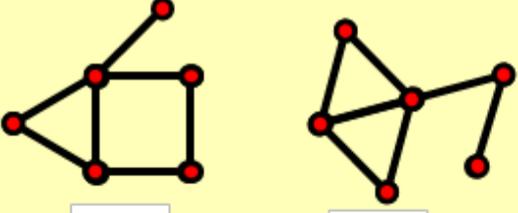 7-29-1 | 1 |

| 7-30 | 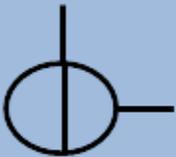 | 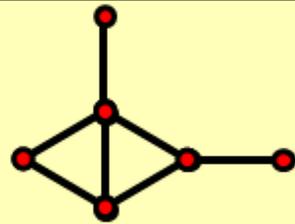 7-30-1 | 1 |

### 3.7.3.2. Matchstick graphs with |E|=7, $F$=3, |V|=6, $\Delta$=4

| 7-31 | 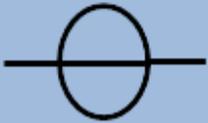 | 7-31-1    7-31-2 | 2 |

| 7-32 | | 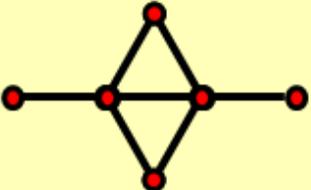 7-32-1 | 1 |

| 7-33 | | 7-33-1 | 1 |

| 7-34 | 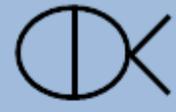 | 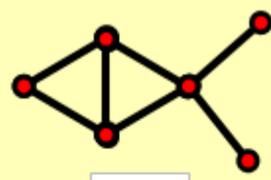 7-34-1 | 1 |
| 7-35 | 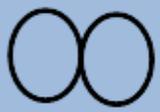 | 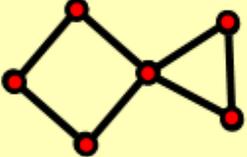 7-35-1 | 1 |
| 7-36 | 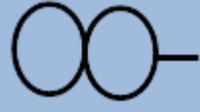 | 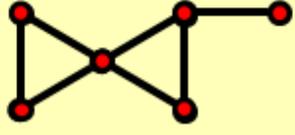 7-36-1 | 1 |

#### 3.7.3.3. Matchstick graphs with |E|=7, F=3, |V|=6, Δ=5

| 7-37 | 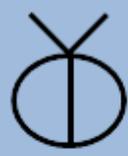 | 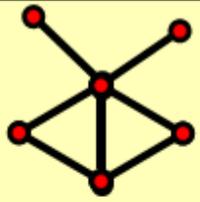 7-37-1 | 1 |
| 7-38 | 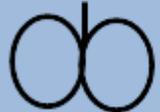 | 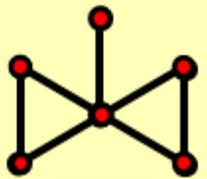 7-38-1 | 1 |

### 3.7.4. Matchstick graphs with |E|=7, F=4, |V|=5
#### 3.7.4.1. Matchstick graphs with |E|=7, F=4, |V|=5, Δ=4

| 7-39 | 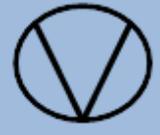 | 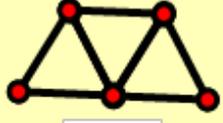 7-39-1 | 1 |

## 3.8. Matchstick graphs with |E|=8
### 3.8.1. Matchstick graphs with |E|=8, F=1, |V|=9
#### 3.8.1.1. Matchstick graphs with |E|=8, F=1, |V|=9, Δ=2

| | | |
|---|---|---|
| 8-1 | 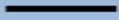 | 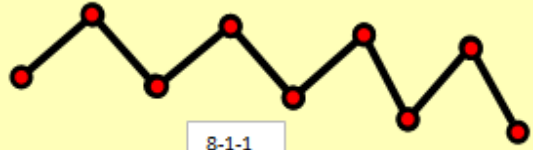  8-1-1 |

### 3.8.1.2. Matchstick graphs with |E|=8, $F$=1, |V|=9, $\Delta$=3

| | | |
|---|---|---|
| 8-2 | 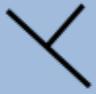 | 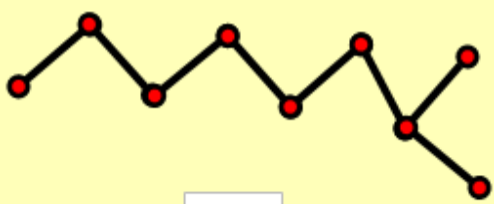 8-2-1  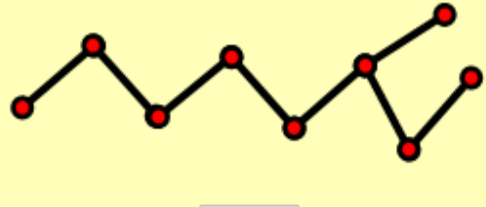 8-2-2 <br> 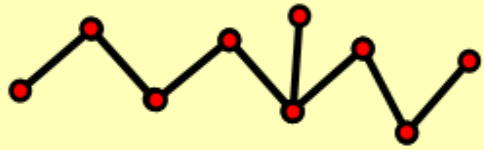 8-2-3  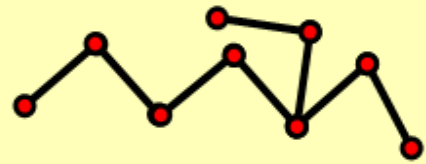 8-2-4 <br> 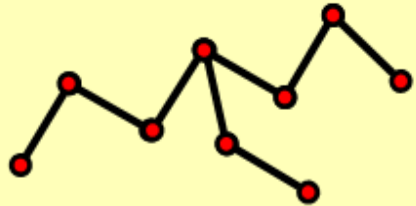 8-2-5 |
| 8-3 | | 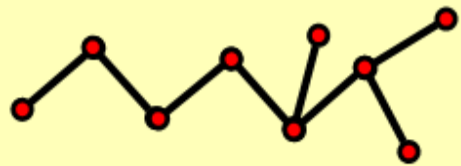 8-3-1  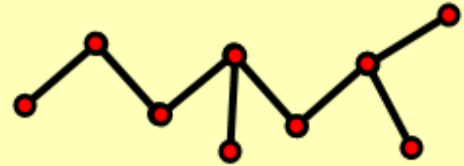 8-3-2 <br> 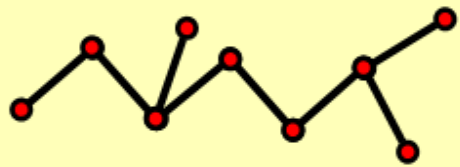 8-3-3  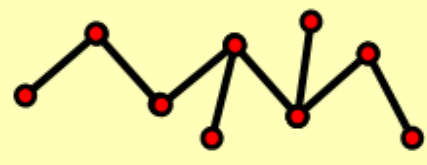 8-3-4 <br> 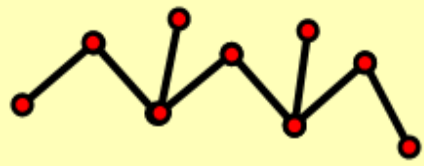 8-3-5  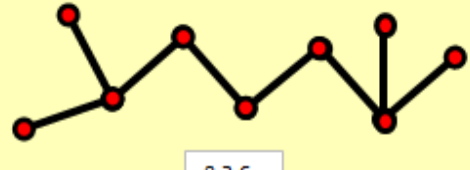 8-3-6 |

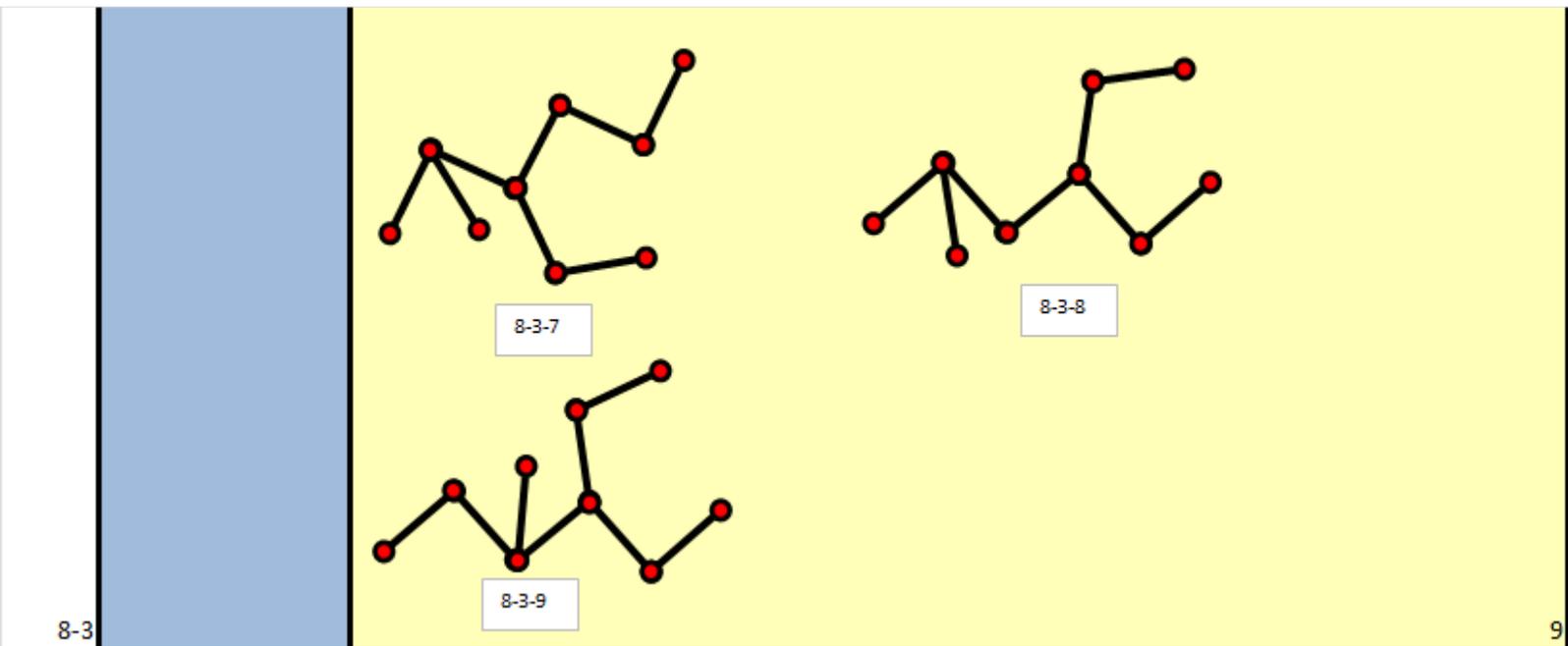

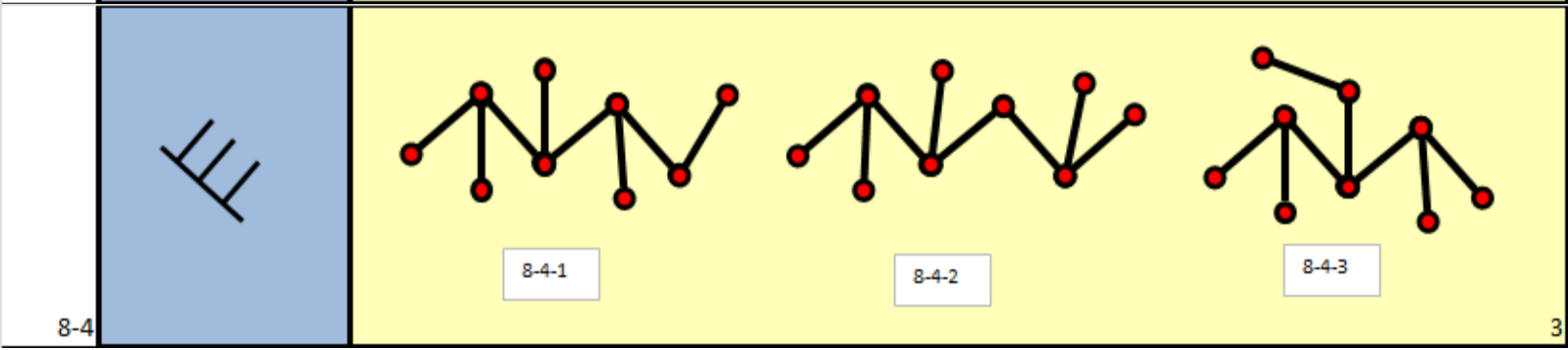

### 3.8.1.3. Matchstick graphs with |E|=8, $F$=1, |V|=9, $\Delta$=4

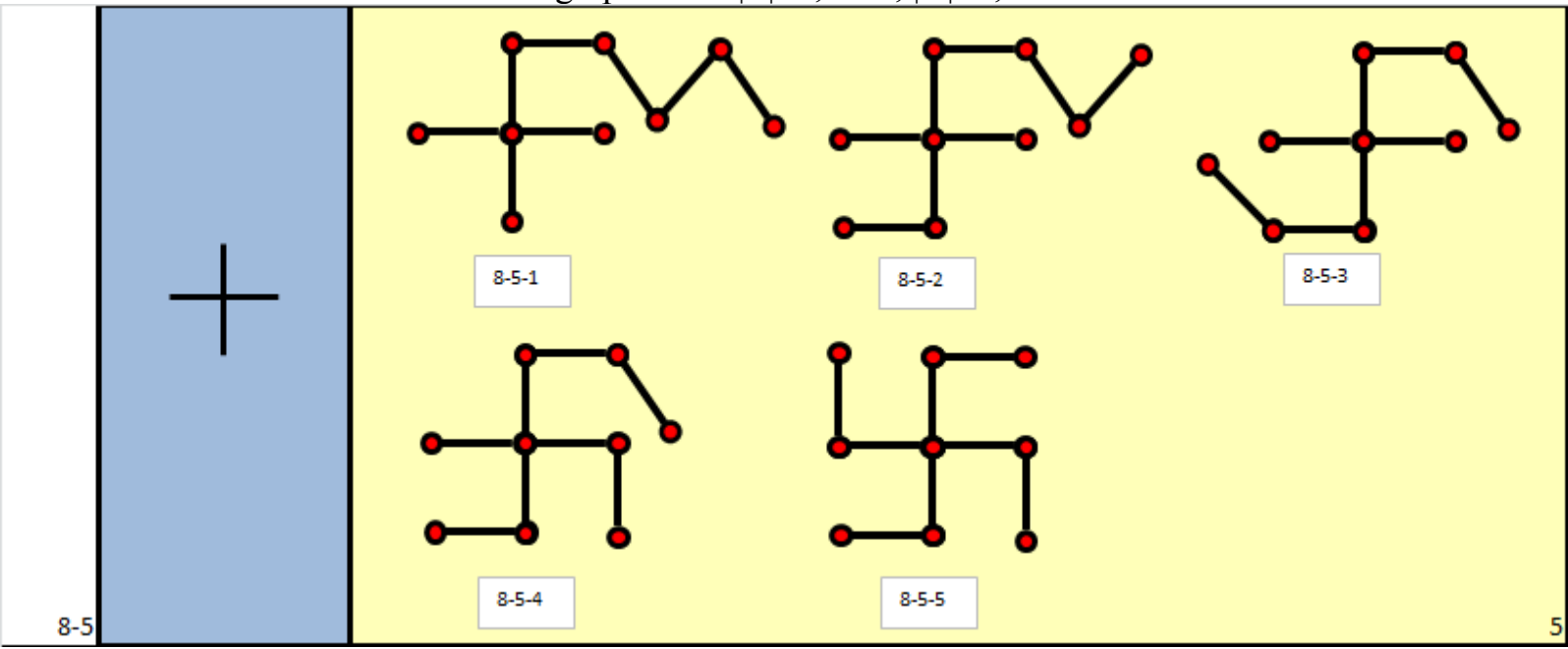

| | | |
|---|---|---|
| 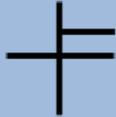 | 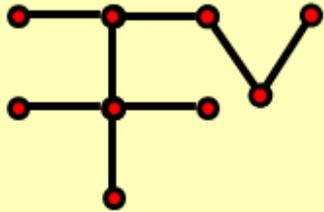 | 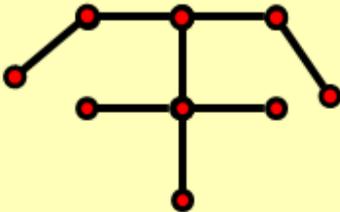 |
| 8-6-1 | 8-6-2 | 8-6-3 |
| 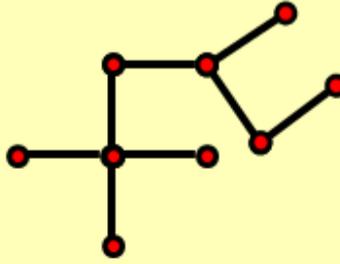 | 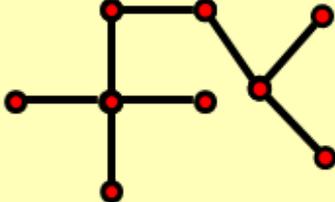 | 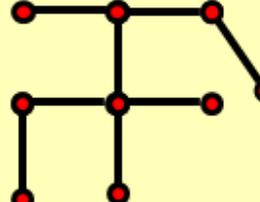 |
| 8-6-4 | 8-6-5 | 8-6-6 |
| 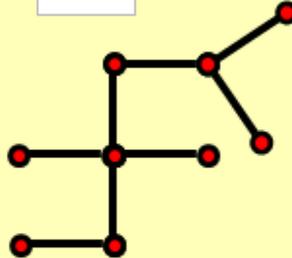 | 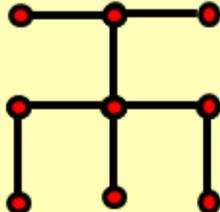 | |
| 8-6-7 | 8-6-8 | |

8-6

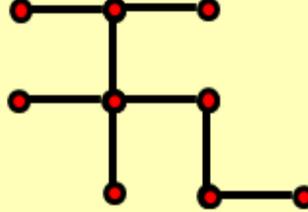
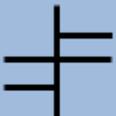

8-7-1

8-7

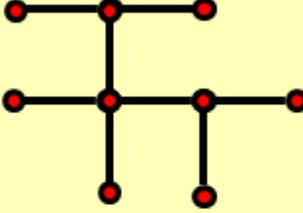
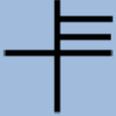

8-8-1

8-8

| 8-9 | 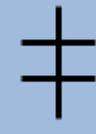 | 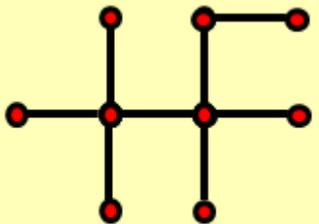 | 2 |

8-9-1, 8-9-2

### 3.8.1.4. Matchstick graphs with |E|=8, F=1, |V|=9, Δ=5

| 8-10 | 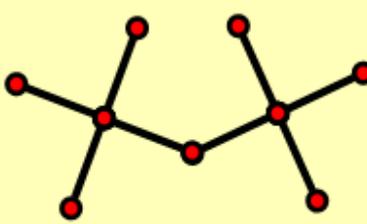 | 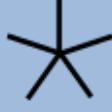 | 3 |

8-10-1, 8-10-2, 8-10-3

| 8-11 | 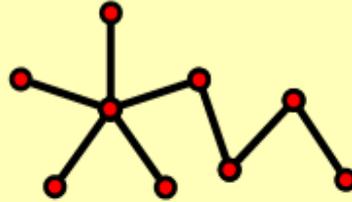 | 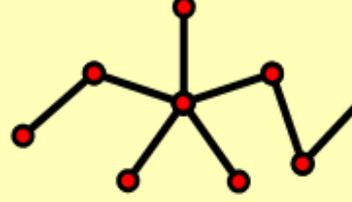 | 3 |

8-11-1, 8-11-2, 8-11-3

| 8-12 | 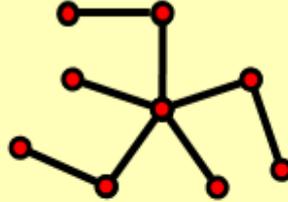 | 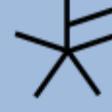 | 1 |

8-12-1

### 3.8.1.5. Matchstick graphs with |E|=8, F=1, |V|=9, Δ=6

| 8-13 | 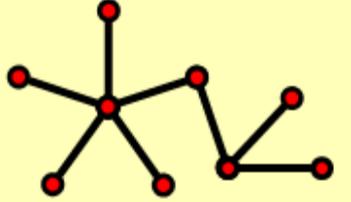 | 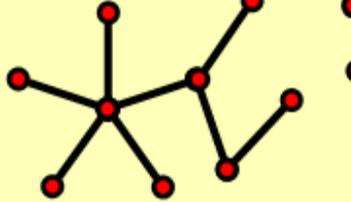 | 2 |

8-13-1, 8-13-2

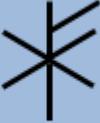

8-14

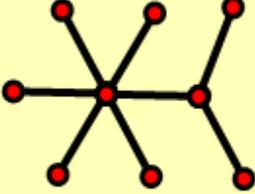

8-14-1

### 3.8.1.6. Matchstick graphs with |E|=8, F=1, |V|=9, Δ=7

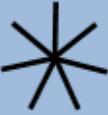

8-15

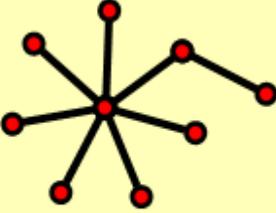

8-15-1

### 3.8.1.7. Matchstick graphs with |E|=8, F=1, |V|=9, Δ=8

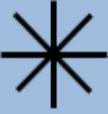

8-16

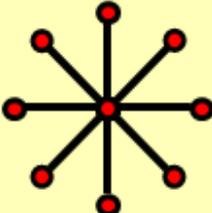

8-16-1

## 3.8.2. Matchstick graphs with |E|=8, F=2, |V|=8
### 3.8.2.1. Matchstick graphs with |E|=8, F=2, |V|=8, Δ=2

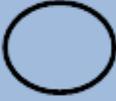

8-17

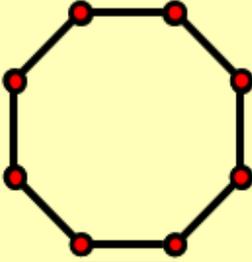

8-17-1

### 3.8.2.2. Matchstick graphs with |E|=8, F=2, |V|=8, Δ=3

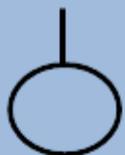

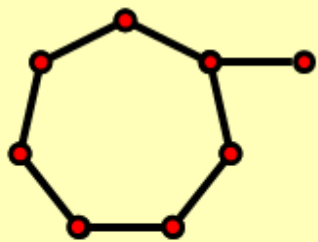
8-18-1

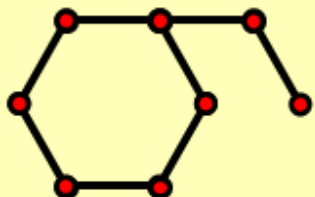
8-18-2

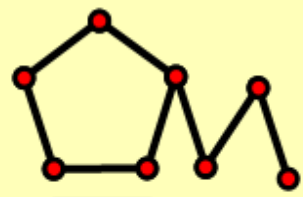
8-18-3

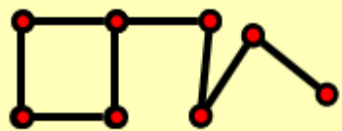
8-18-4

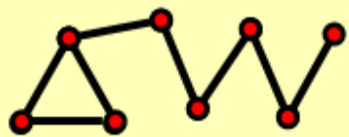
8-18-5

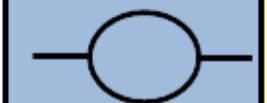

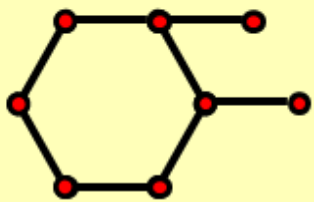
8-19-1

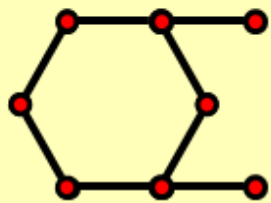
8-19-2

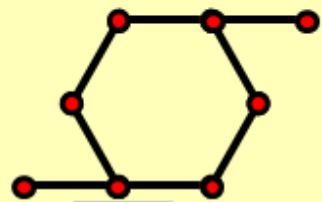
8-19-3

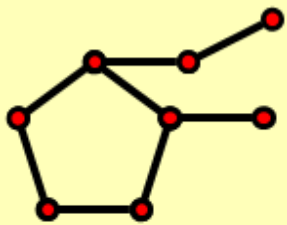
8-19-4

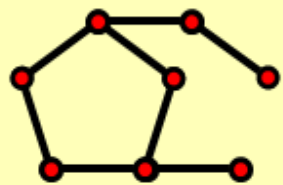
8-19-5

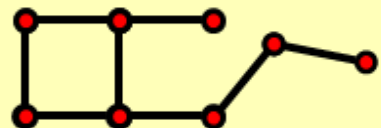
8-19-6

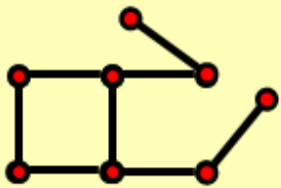
8-19-7

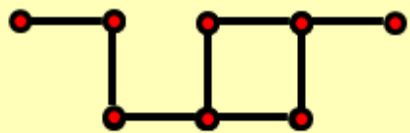
8-19-8

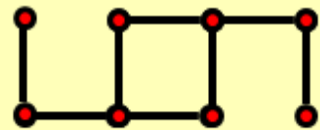
8-19-9

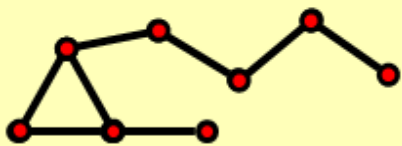
8-19-10

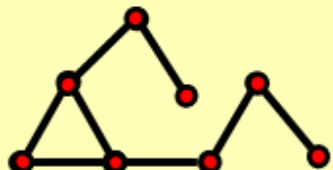
8-19-11

8-18

8-19

| | | |
|---|---|---|
| 8-20 | 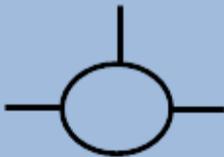 | 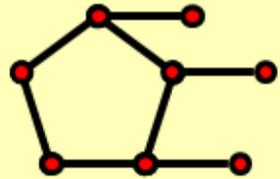 8-20-1    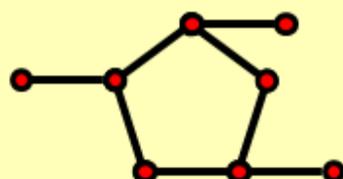 8-20-2    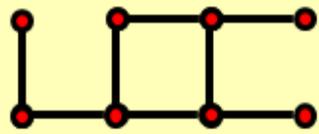 8-20-3 <br> 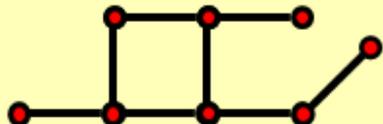 8-20-4    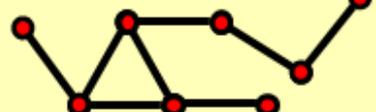 8-20-5    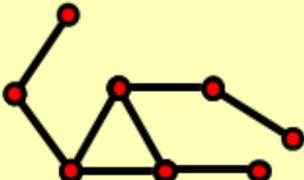 8-20-6 |
| 8-21 | 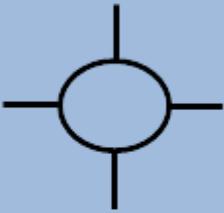 | 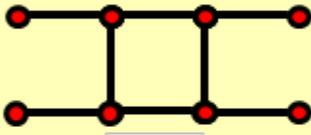 8-21-1 |
| 8-22 | 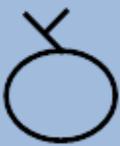 | 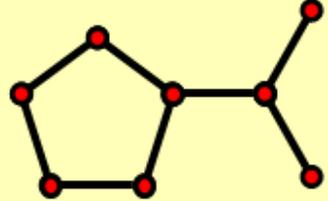 8-22-1    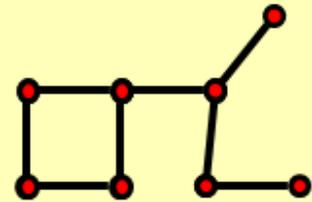 8-22-2    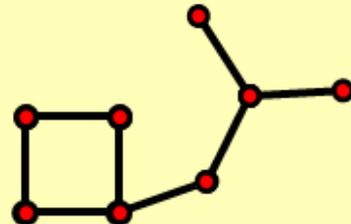 8-22-3 <br> 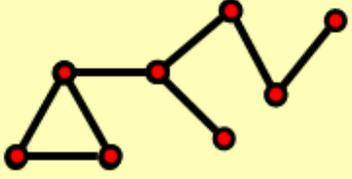 8-22-4    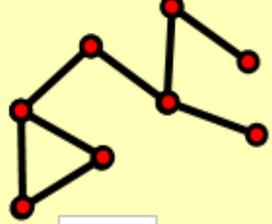 8-22-5    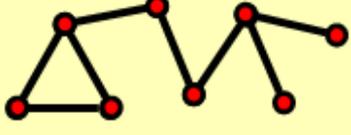 8-22-6 <br> 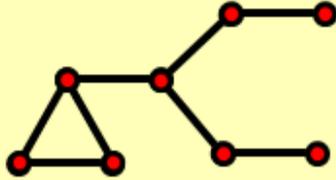 8-22-7 |

| 8-23 | | 8-23-1 | 8-23-2 | 8-23-3 | | |
| --- | --- | --- | --- | --- | --- | --- |
| | | 8-23-4 | 8-23-5 | | | 5 |
| 8-24 | | 8-24-1 | | | | 1 |
| 8-25 | | 8-25-1 | | | | 1 |

3.8.2.3. Matchstick graphs with $|E|=8$, $F=2$, $|V|=8$, $\Delta=4$

| 8-26 | | 8-26-1 | 8-26-2 | 8-26-3 | | |
| --- | --- | --- | --- | --- | --- | --- |
| | | 8-26-4 | 8-26-5 | 8-26-6 | | 6 |

| | | |
|---|---|---|
| 8-27 | 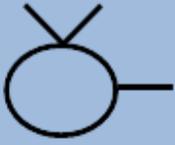 | 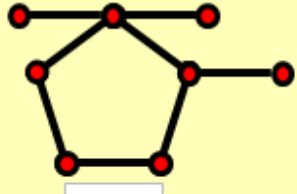 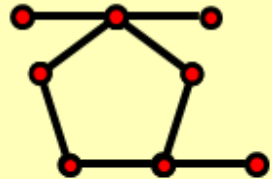 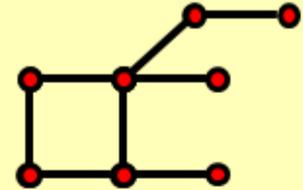  8-27-1  8-27-2  8-27-3  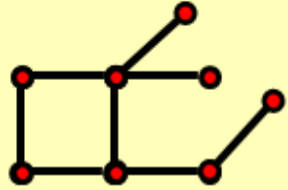 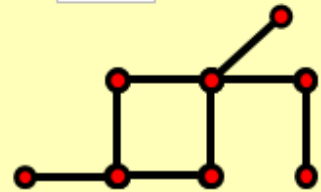 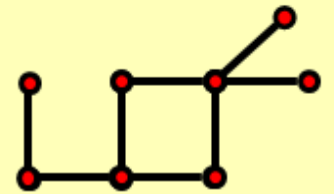  8-27-4  8-27-5  8-27-6  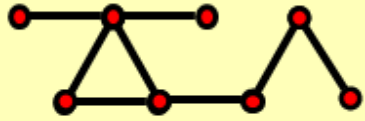 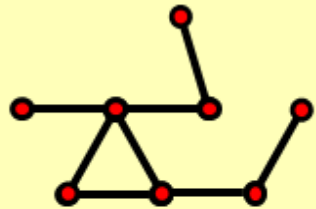 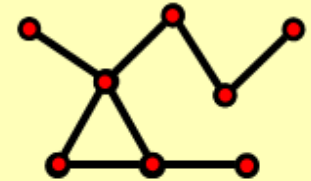  8-27-7  8-27-8  8-27-9  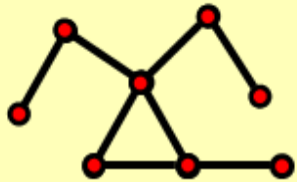  8-27-10 |
| 8-28 | 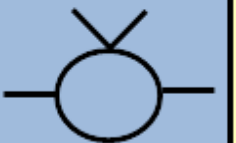 | 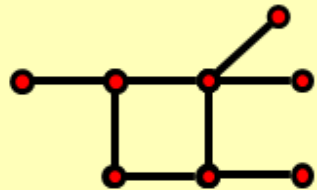 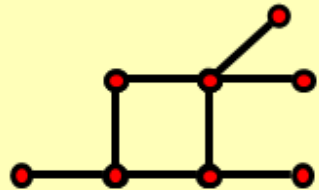 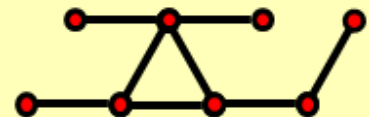  8-28-1  8-28-2  8-28-3  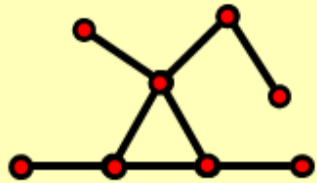  8-28-4 |
| 8-29 | 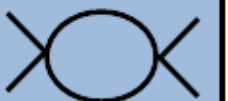 | 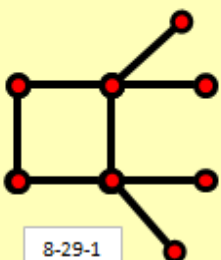 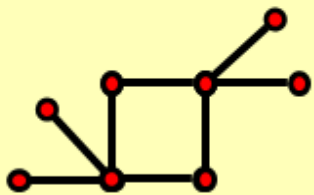 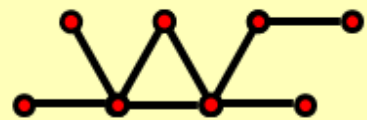  8-29-1  8-29-2  8-29-3 |
| 8-30 | 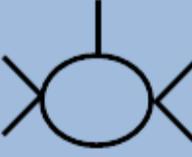 | 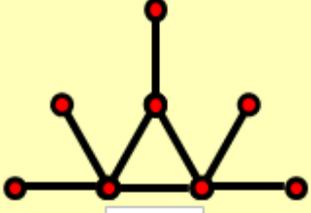 8-30-1 |
| 8-31 | 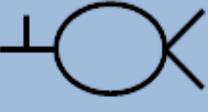 | 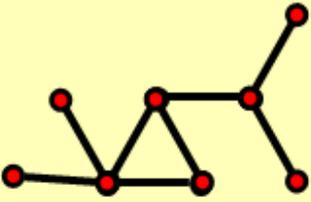 8-31-1 |
| 8-32 | 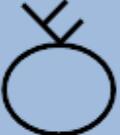 | 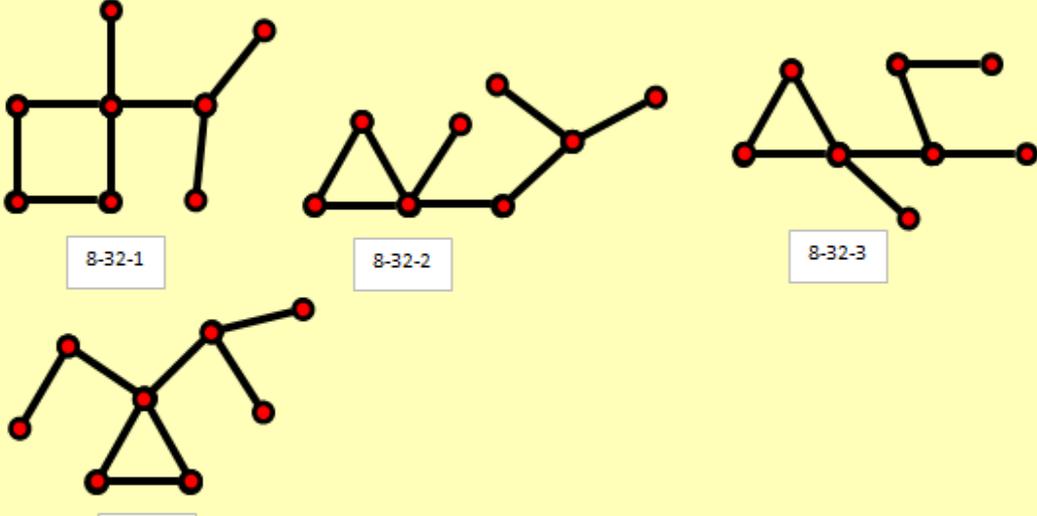 8-32-1, 8-32-2, 8-32-3, 8-32-4 |
| 8-33 | 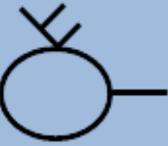 | 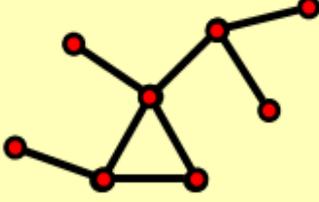 8-33-1 |
| 8-34 | 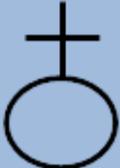 | 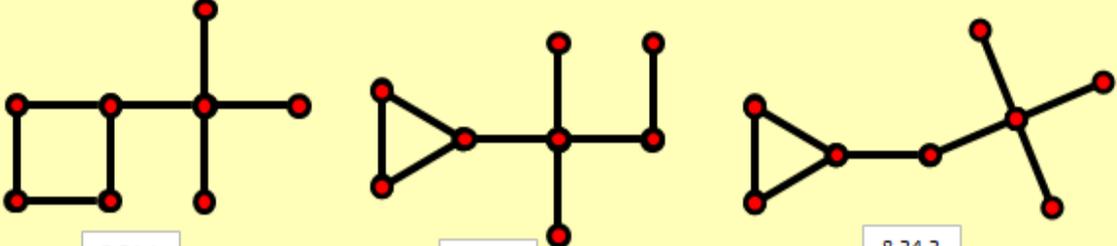 8-34-1, 8-34-2, 8-34-3 |

| 8-35 | 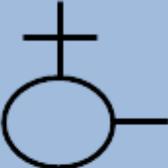 | 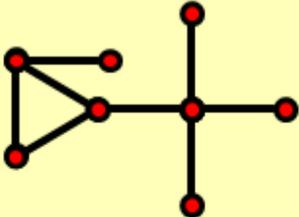  8-35-1 | | | 1 |

| 8-36 | 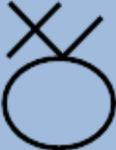 | 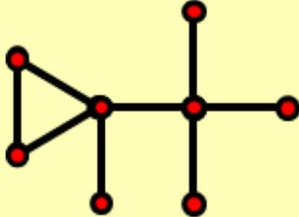  8-36-1 | | | 1 |

### 3.8.2.4. Matchstick graphs with |E|=8, F=2, |V|=8, Δ=5

| 8-37 | 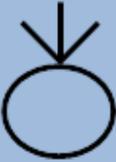 | 8-37-1, 8-37-2, 8-37-3, 8-37-4 | 4 |

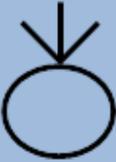

| 8-38 | 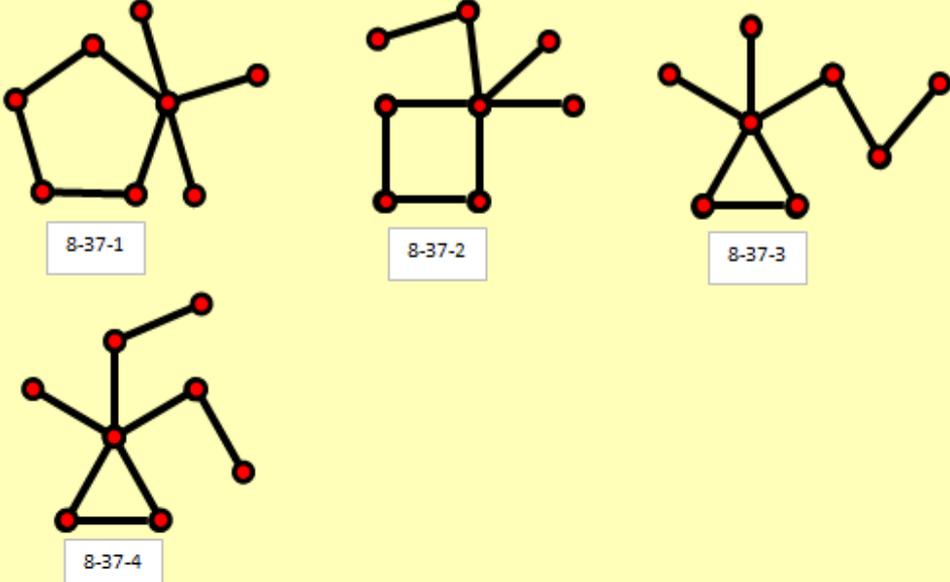 | 8-38-1, 8-38-2, 8-38-3, 8-38-4 | 4 |

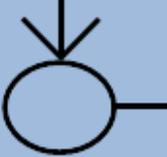
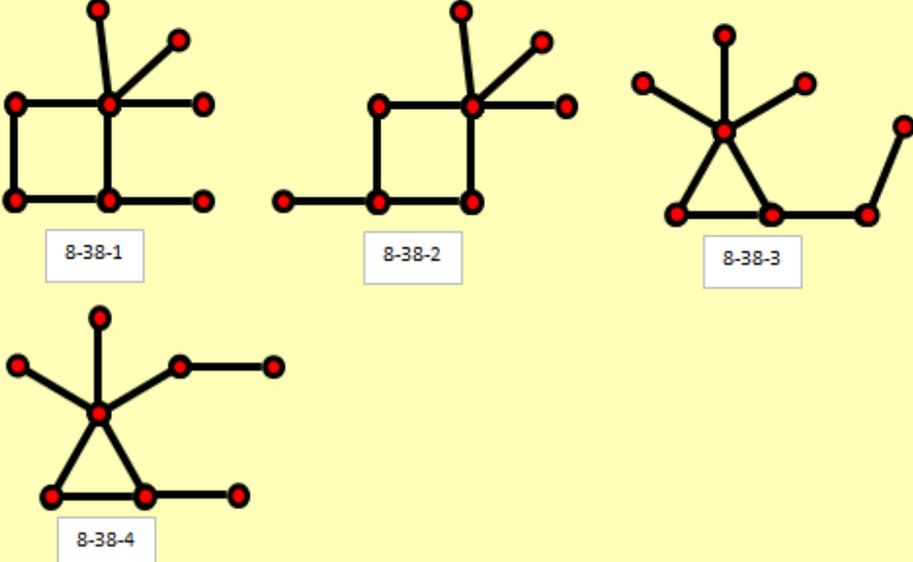
27

| | | |
|---|---|---|
| 8-39 | 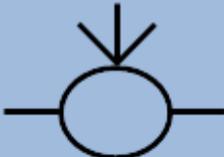 | 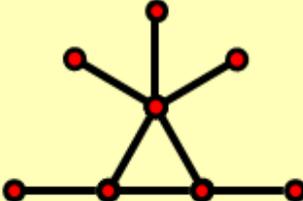<br>8-39-1 |
| 8-40 | 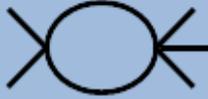 | 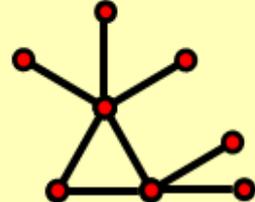<br>8-40-1 |
| 8-41 | 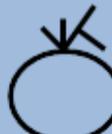 | 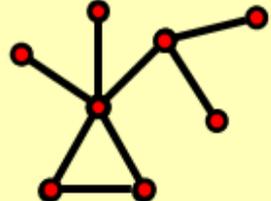<br>8-41-1 |
| 8-42 | 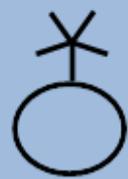 | 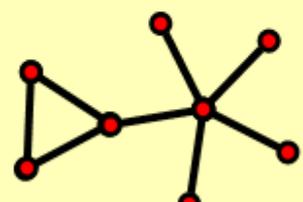<br>8-42-1 |

3.8.2.5. Matchstick graphs with $|E|=8$, $F=2$, $|V|=8$, $\Delta=6$

| | | |
|---|---|---|
| 8-43 | 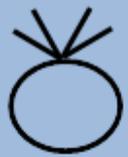 | 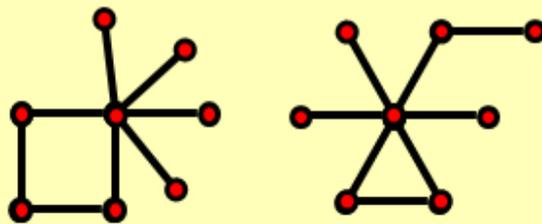<br>8-43-1    8-43-2 |
| 8-44 | 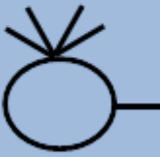 | 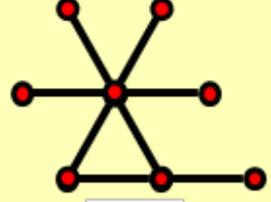<br>8-44-1 |

### 3.8.2.6. Matchstick graphs with |E|=8, F=2, |V|=8, Δ=7

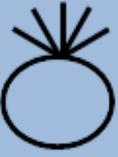

8-45-1

8-45 | 1

### 3.8.3. Matchstick graphs with |E|=8, F=3, |V|=7
#### 3.8.3.1. Matchstick graphs with |E|=8, F=3, |V|=7, Δ=3

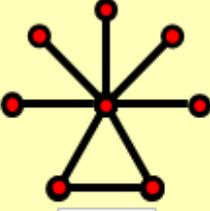

8-46-1  8-46-2  8-46-3  8-46-4

8-46 | 4

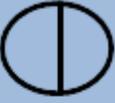

8-47-1  8-47-2  8-47-3

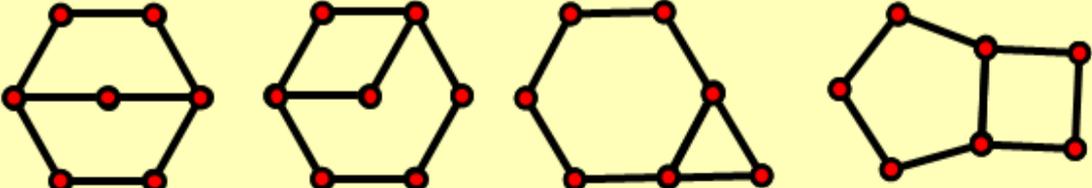

8-47-4  8-47-5  8-47-6

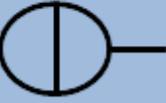

8-47-7  8-47-8  8-47-9

8-47 | 9

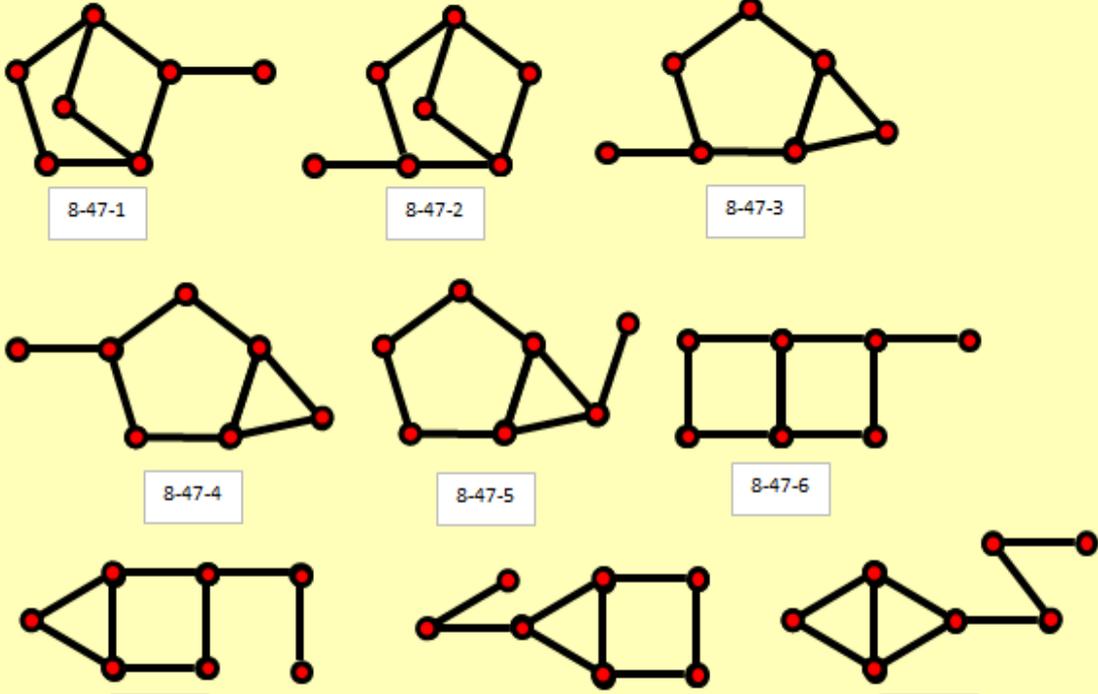

8-48-1

8-48 | 1

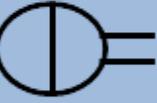

8-49-1  8-49-2

8-49 | 2

| | | |
|---|---|---|
| 8-50 | 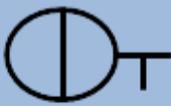 | 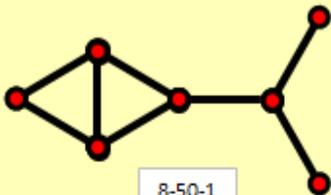 8-50-1 |
| 8-51 | 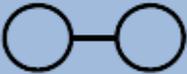 | 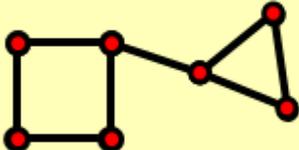 8-51-1    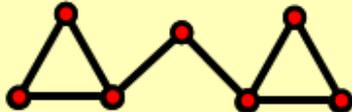 8-51-2 |
| 8-52 | 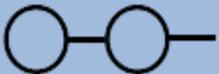 | 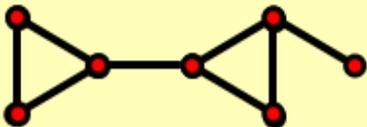 8-52-1 |

### 3.8.3.2. Matchstick graphs with |E|=8, $F$=3, |V|=7, $\Delta$=4

| | | |
|---|---|---|
| 8-53 | 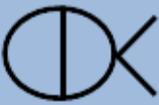 | 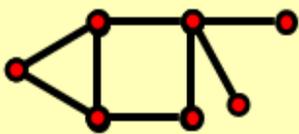 8-53-1   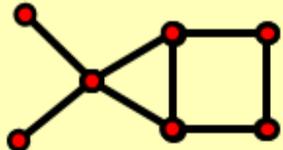 8-53-2   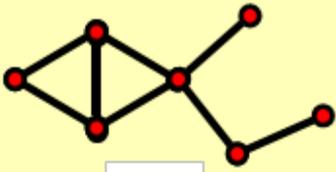 8-53-3 |
| 8-54 | 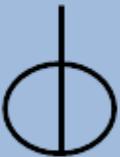 | 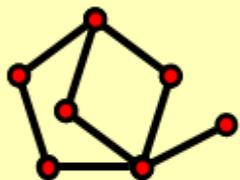 8-54-1   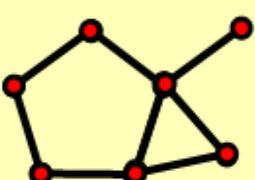 8-54-2   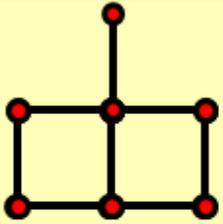 8-54-3   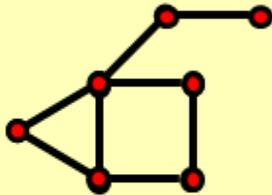 8-54-4 <br> 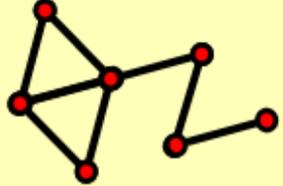 8-54-5 |
| 8-55 | 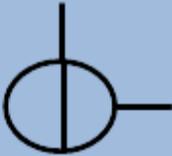 | 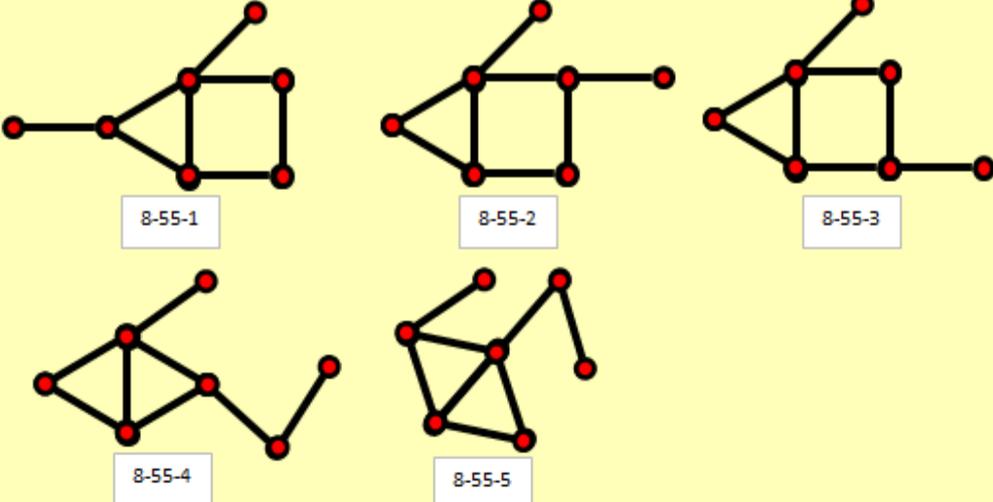 8-55-1   8-55-2   8-55-3 <br> 8-55-4   8-55-5 | 5
| 8-56 | 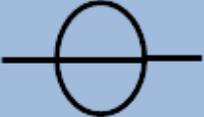 | 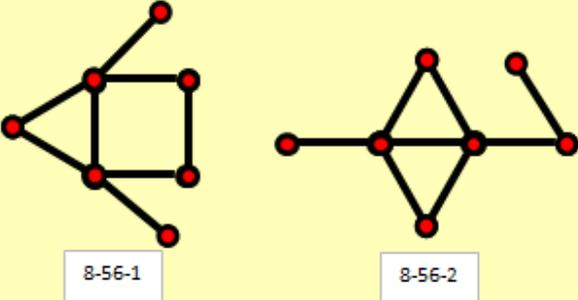 8-56-1   8-56-2 | 2
| 8-57 | 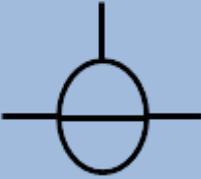 | 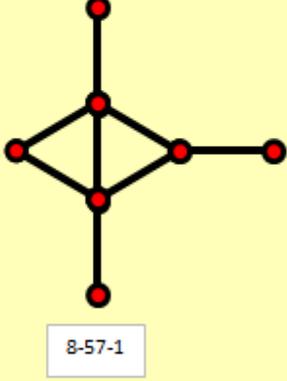 8-57-1 | 1
| 8-58 | 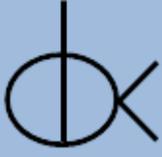 | 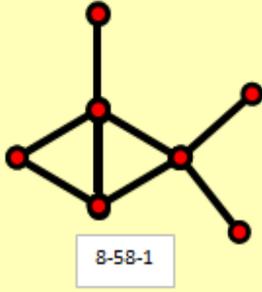 8-58-1 | 1
| 8-59 | 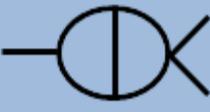 | 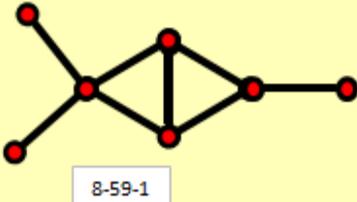 8-59-1 | 1

| | | |
|---|---|---|
| 8-60 | 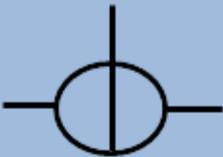 | 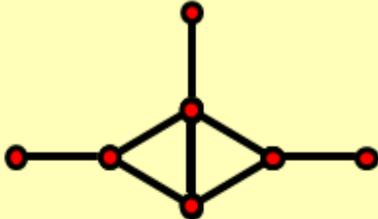 8-60-1 |
| 8-61 | 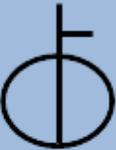 | 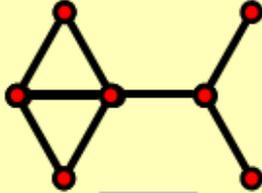 8-61-1 |
| 8-62 | 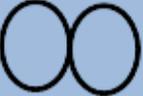 | 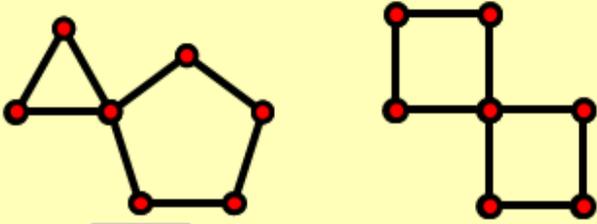 8-62-1   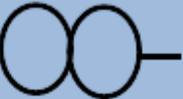 8-62-2 |
| 8-63 | 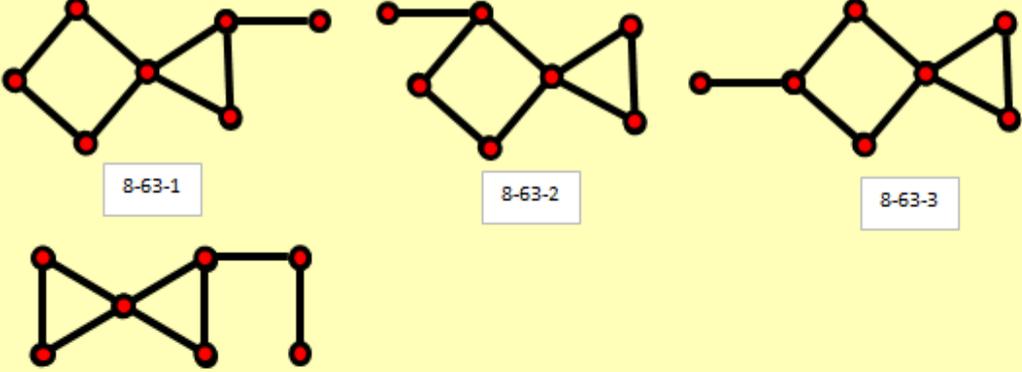 | 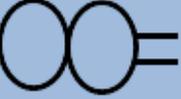 8-63-1   8-63-2   8-63-3   <br> 8-63-4 |
| 8-64 | 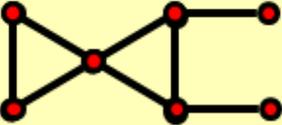 | 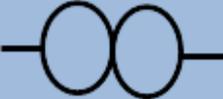 8-64-1 |
| 8-65 | 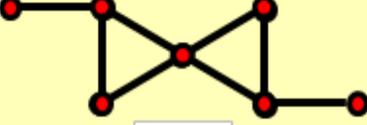 | 8-65-1 |

| 8-66 | 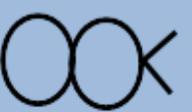 | 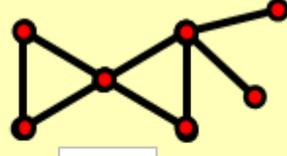<br>8-66-1 | 1 |
| 8-67 | 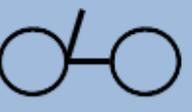 | 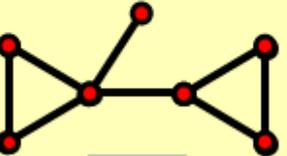<br>8-67-1 | 1 |

### 3.8.3.3. Matchstick graphs with |E|=8, *F*=3, |V|=7, Δ=5

| 8-68 | 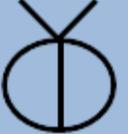 | 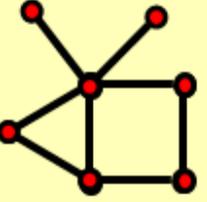<br>8-68-1     8-68-2 | 2 |
| 8-69 | 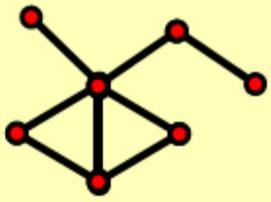 | 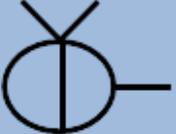<br>8-69-1 | 1 |
| 8-70 | 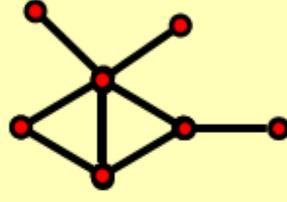 | 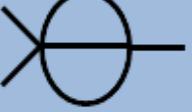<br>8-70-1 | 1 |
| 8-71 | 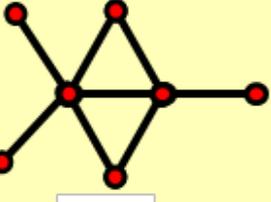 | 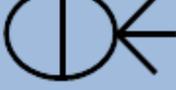<br>8-71-1 | 1 |
| 8-72 | 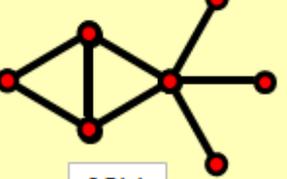 | 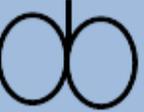 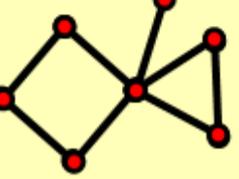<br>8-72-1     8-72-2 | 2 |

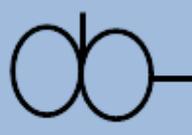

8-73

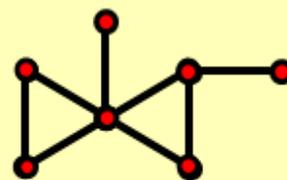

8-73-1

### 3.8.3.4. Matchstick graphs with |E|=8, F=3, |V|=7, Δ=6

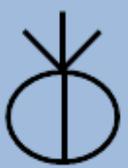

8-74

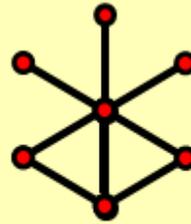

8-74-1

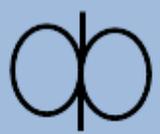

8-75

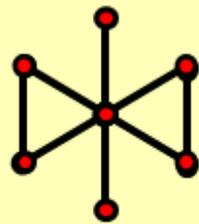

8-75-1

### 3.8.4. Matchstick graphs with |E|=8, F=4, |V|=6
### 3.8.4.1. Matchstick graphs with |E|=8, F=4, |V|=6, Δ=3

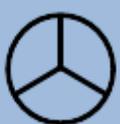

8-76

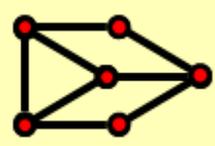

8-76-1

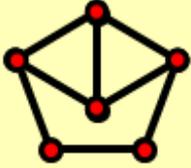

8-76-2

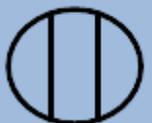

8-77

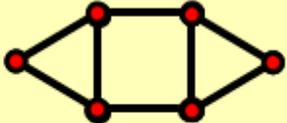

8-77-1

### 3.8.4.2. Matchstick graphs with |E|=8, F=4, |V|=6, Δ=4

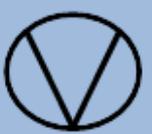

8-78

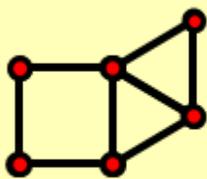

8-78-1

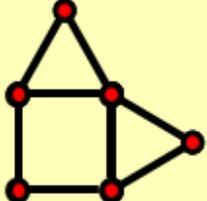

8-78-2

| 8-79 | 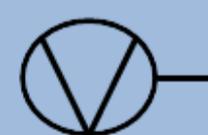 | 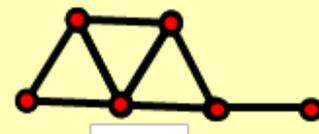 8-79-1 | 1 |

| 8-80 | 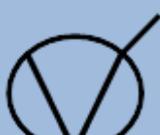 | 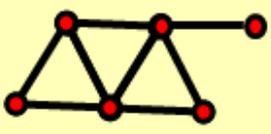 8-80-1 | 1 |

| 8-81 | 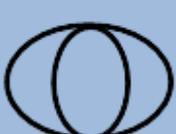 | 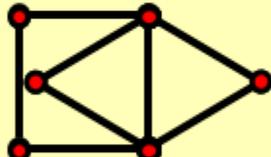 8-81-1 | 1 |

| 8-82 | 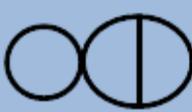 | 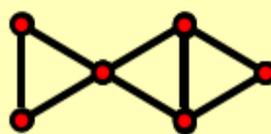 8-82-1 | 1 |

### 3.8.4.3. Matchstick graphs with |E|=8, F=4, |V|=6, Δ=5

| 8-83 | 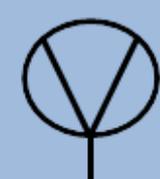 | 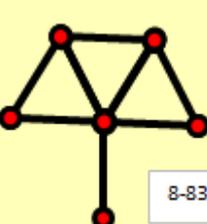 8-83-1 | 1 |

| 8-84 | 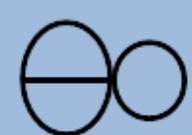 | 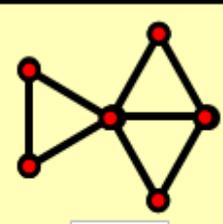 8-84-1 | 1 |

## 3.9. Matchstick graphs with |E|=9
### 3.9.1. Matchstick graphs with |E|=9, F=1, |V|=10
#### 3.9.1.1. Matchstick graphs with |E|=9, F=1, |V|=10, Δ=2

| 9-1 | 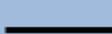 | 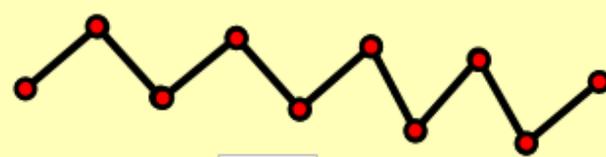 9-1-1 | 1 |

## 3.9.1.2. Matchstick graphs with |E|=9, F=1, |V|=10, Δ=3

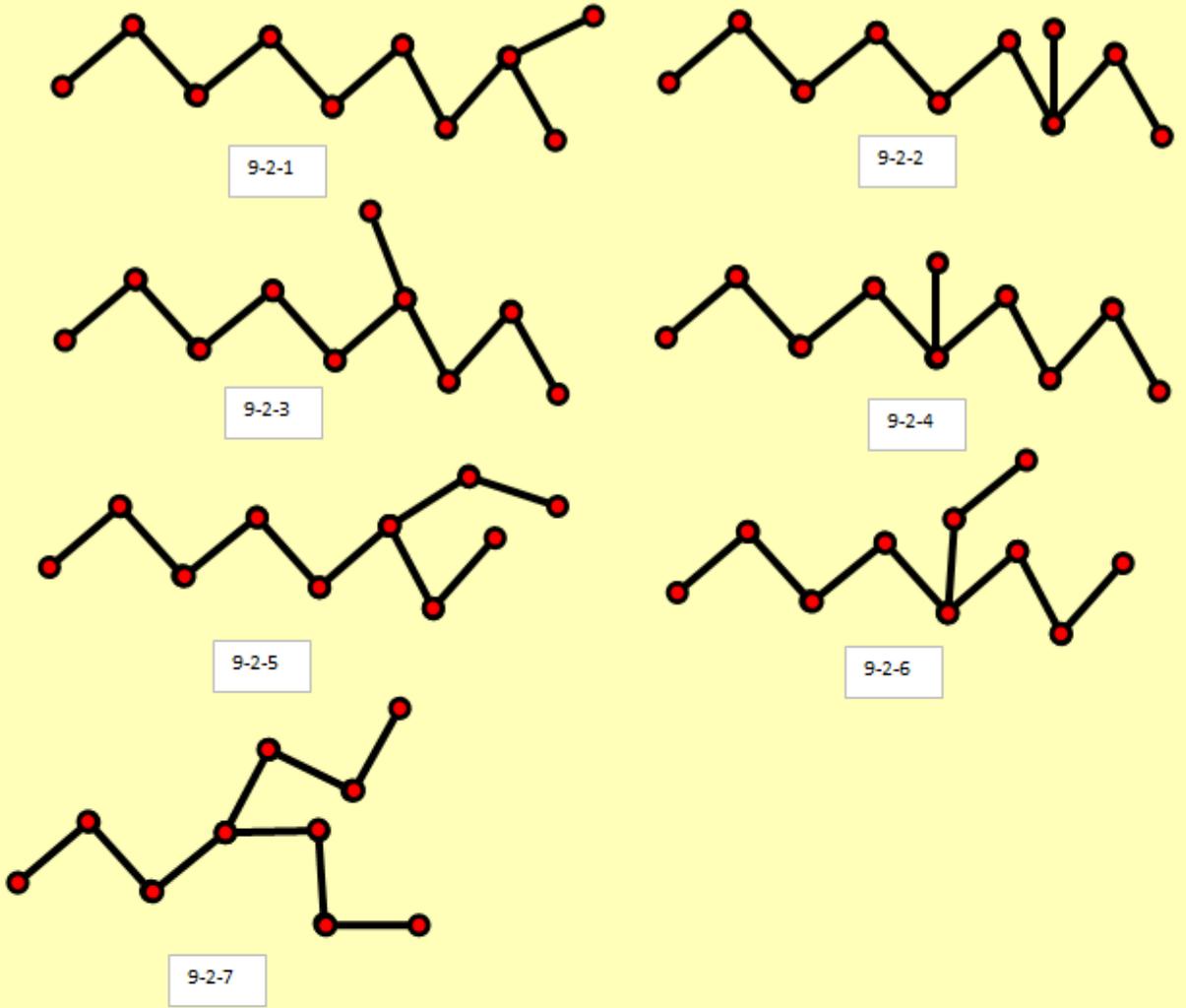

9-2-1   9-2-2   9-2-3   9-2-4   9-2-5   9-2-6   9-2-7
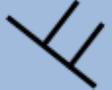

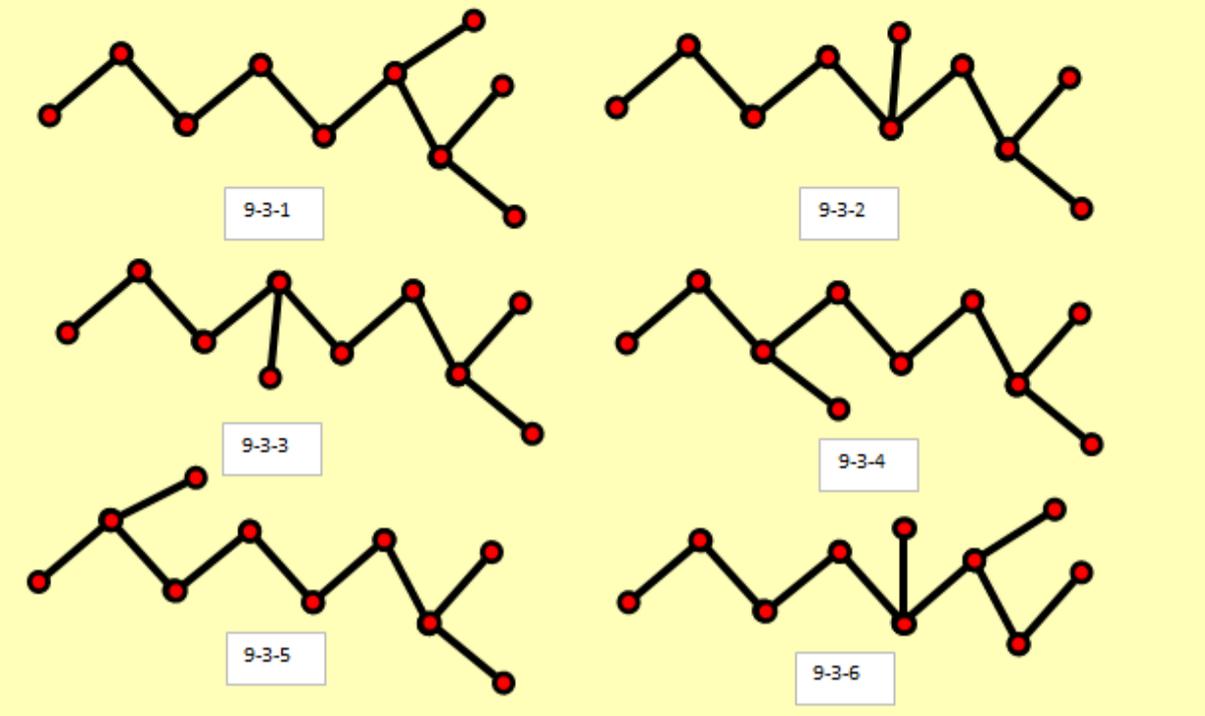

9-3-1   9-3-2   9-3-3   9-3-4   9-3-5   9-3-6

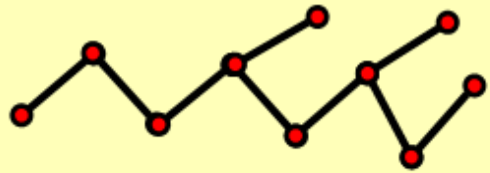

9-3-7

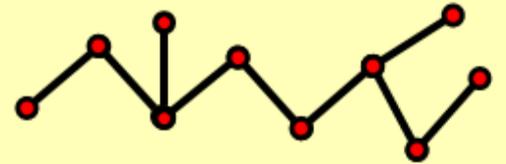

9-3-8

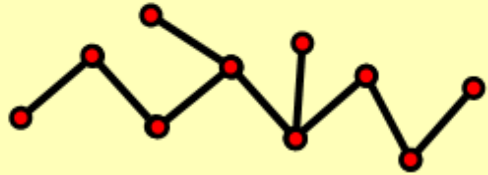

9-3-9

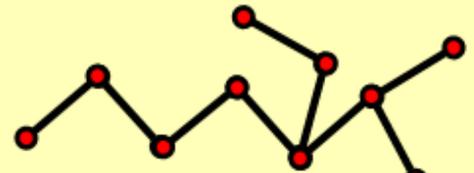

9-3-10

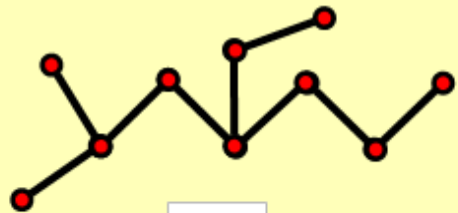

9-3-11

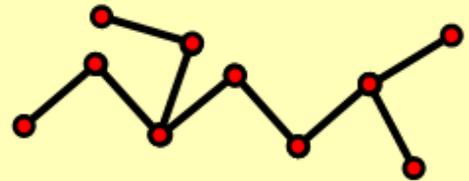

9-3-12

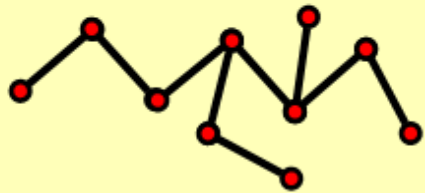

9-3-13

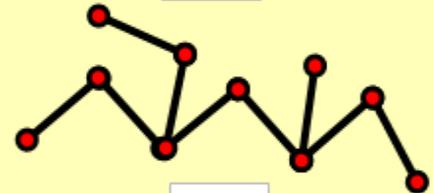

9-3-14

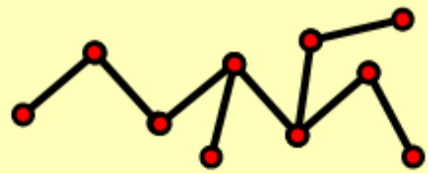

9-3-15

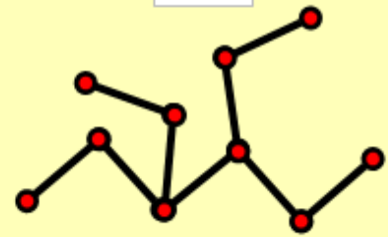

9-3-16

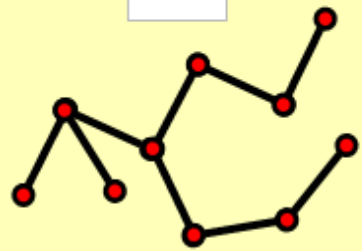

9-3-17

9-3    1737

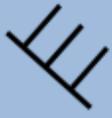

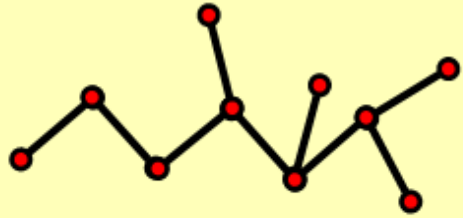
9-4-1

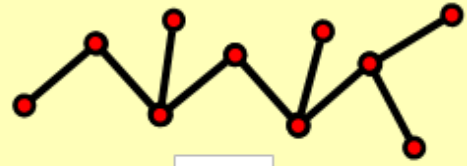
9-4-2

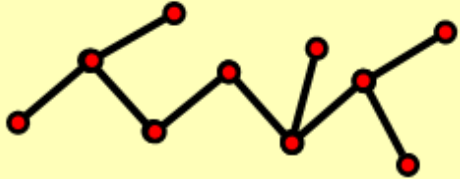
9-4-3

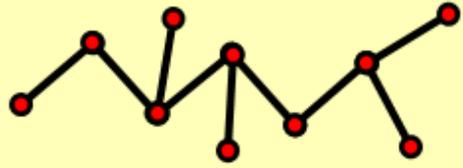
9-4-4

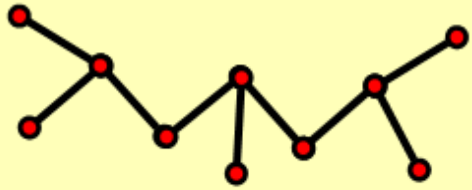
9-4-5

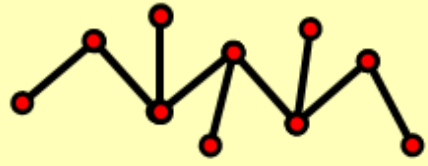
9-4-6

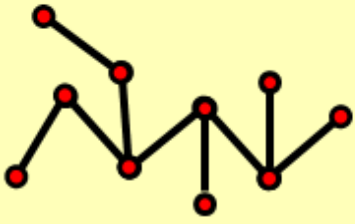
9-4-7

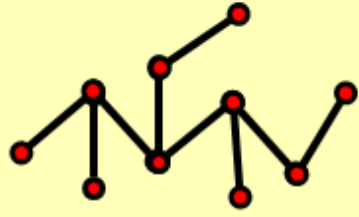
9-4-8

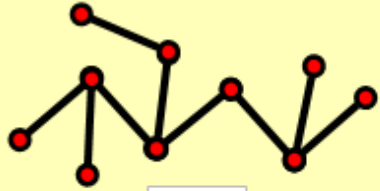
9-4-9

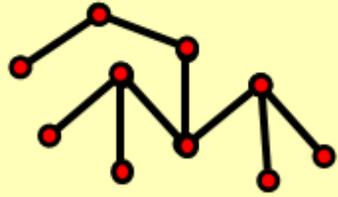
9-4-10

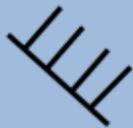

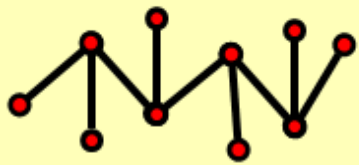
9-5-1

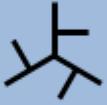
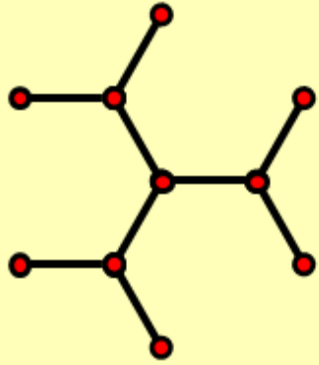

9-6-1
9-6

### 3.9.1.3. Matchstick graphs with |E|=9, $F$=1, |V|=10, $\Delta$=4

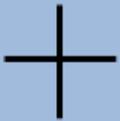

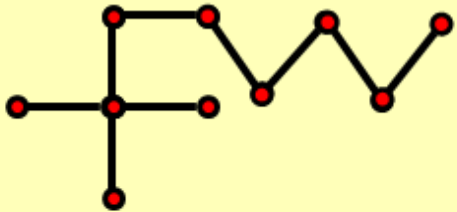
9-7-1

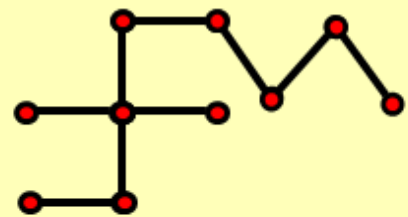
9-7-2

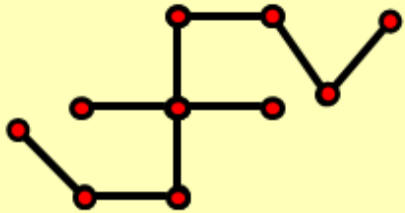
9-7-3

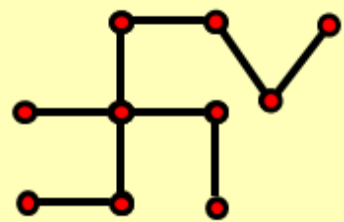
9-7-4

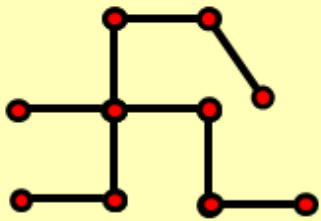
9-7-5

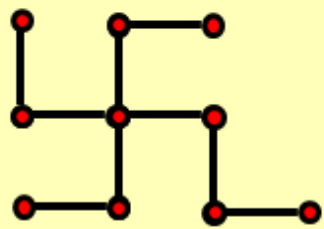
9-7-6

9-7

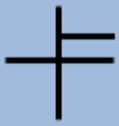

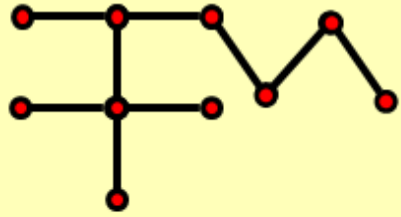

9-8-1

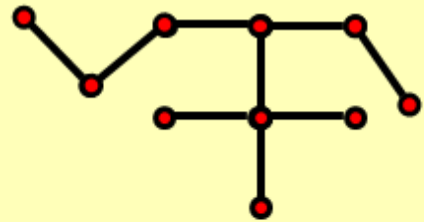

9-8-2

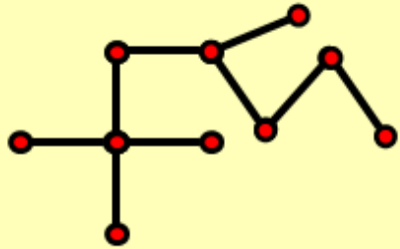

9-8-3

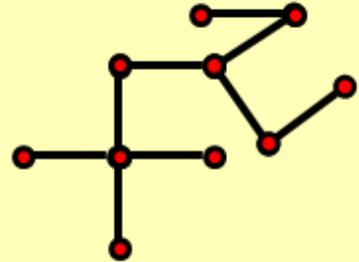

9-8-4

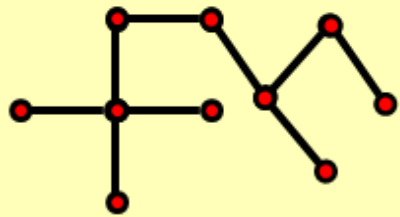

9-8-5

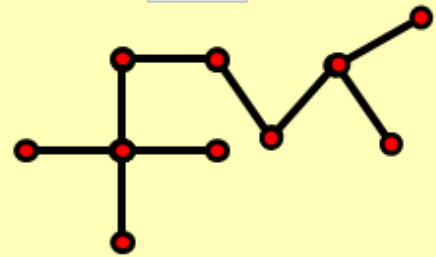

9-8-6

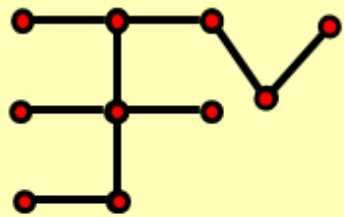

9-8-7

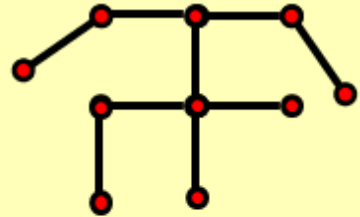

9-8-8

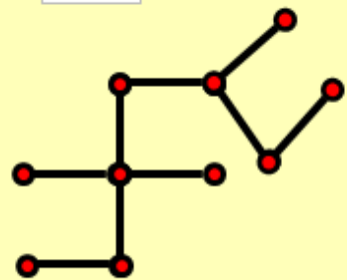

9-8-9

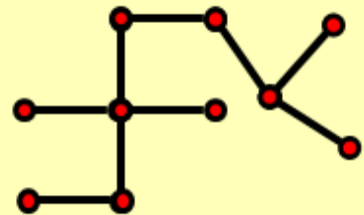

9-8-10

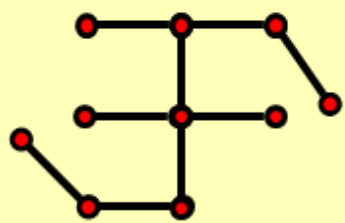

9-8-11

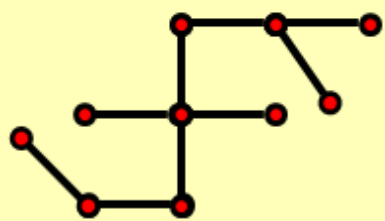

9-8-12

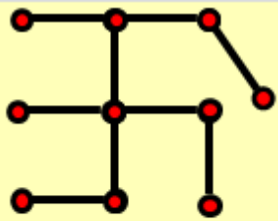
9-8-13

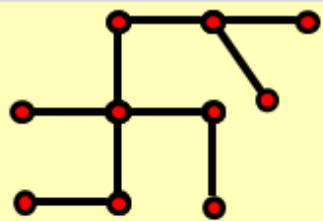
9-8-14

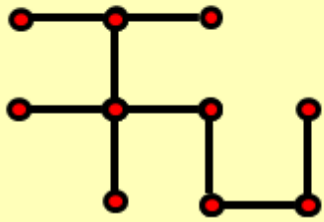
9-8-15

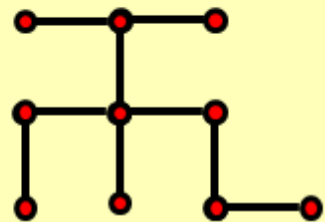
9-8-16

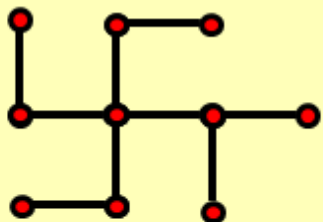
9-8-17

9-8    17

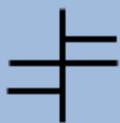

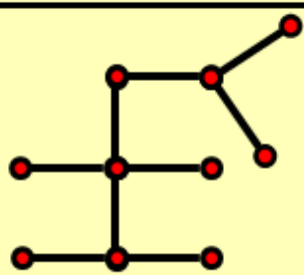
9-9-1

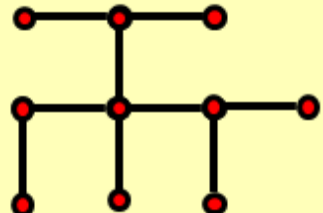
9-9-2

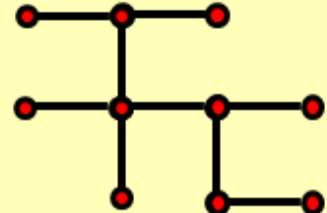
9-9-3

9-9    3

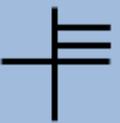

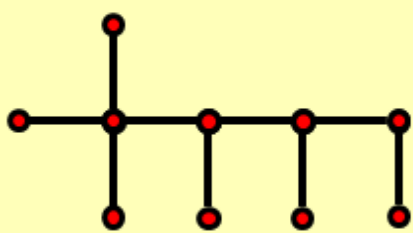
9-10-1

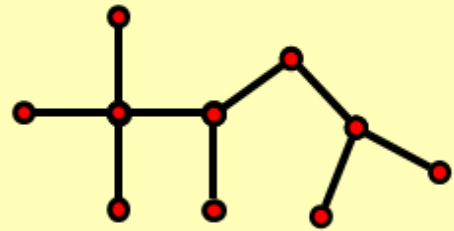
9-10-2

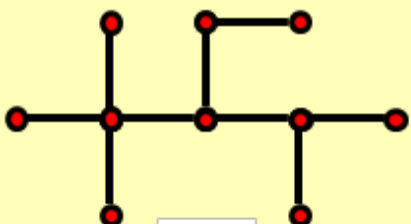
9-10-3

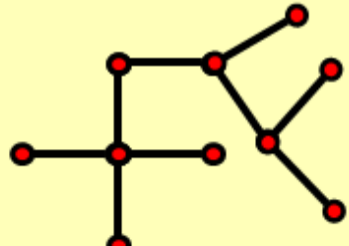
9-10-4

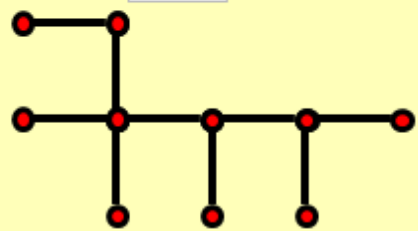
9-10-5

9-10    5

| | | |
|---|---|---|
| 9-11 | 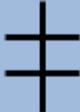 | 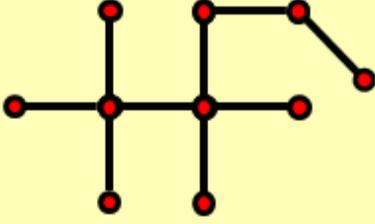 9-11-1  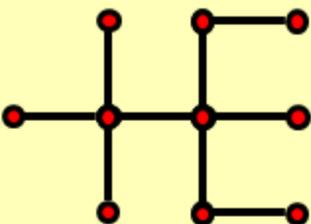 9-11-2  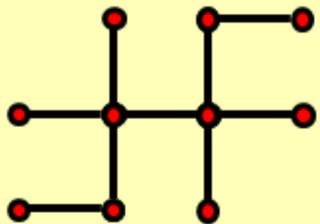 9-11-3 <br> 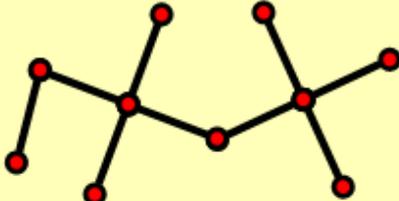 9-11-4  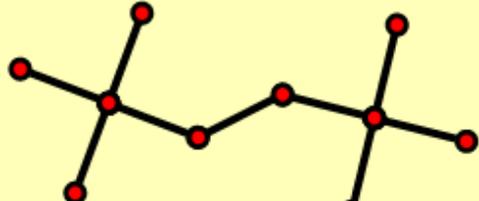 9-11-5 |
| 9-12 | 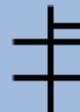 | 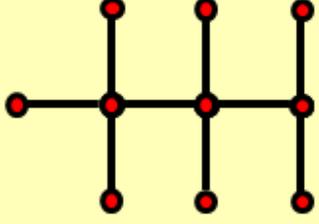 9-12-1 |
| 9-13 | 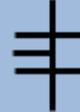 | 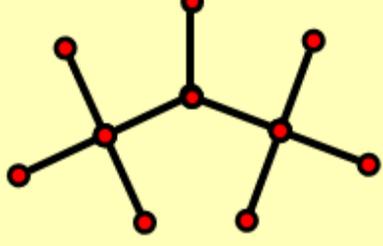 9-13-1 |

3.9.1.4. Matchstick graphs with |E|=9, $F$=1, |V|=10, $\Delta$=5

| | |
|---|---|
| 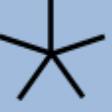 | 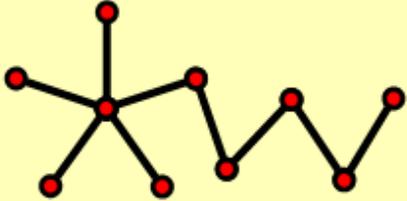 9-14-1  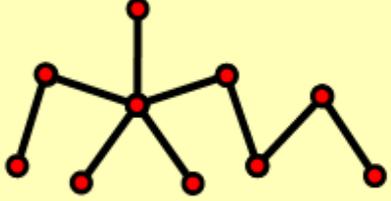 9-14-2 |

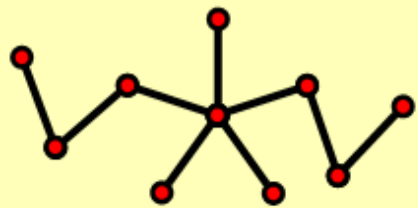

9-14-3

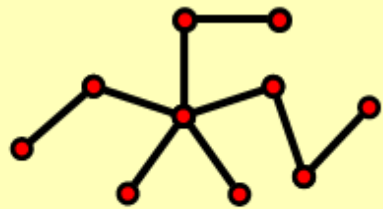

9-14-4

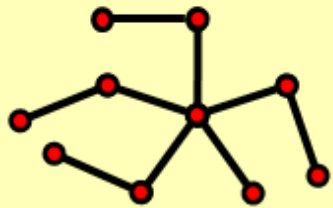

9-14-5

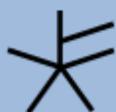

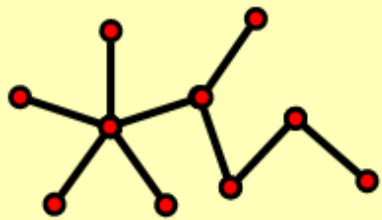

9-15-1

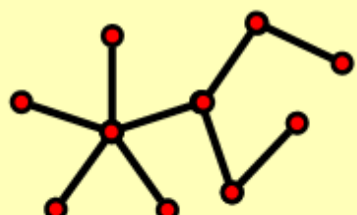

9-15-2

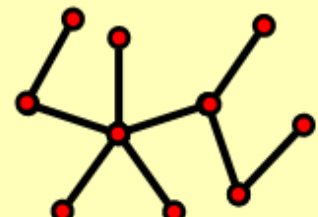

9-15-5

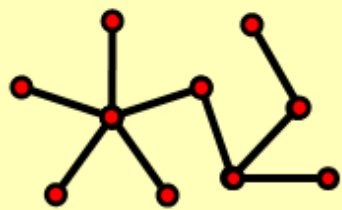

9-15-3

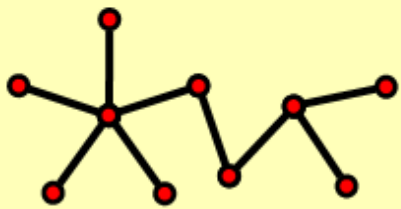

9-15-4

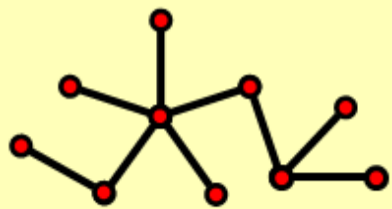

9-15-6

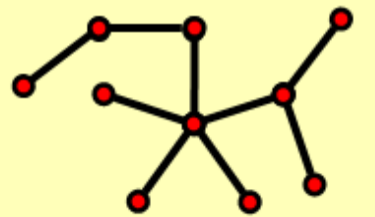

9-15-7

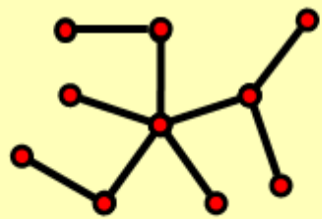

9-15-8

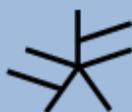

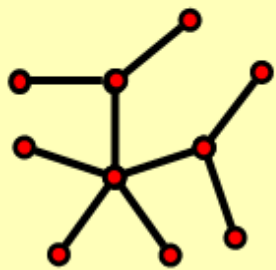

9-16-1

9-14

9-15

9-16

| | | |
|---|---|---|
| 9-17 | 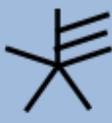 | 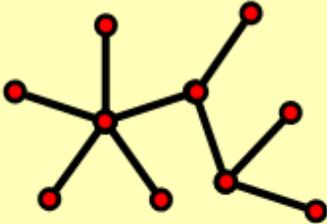 9-17-1 |

| | | | | |
|---|---|---|---|---|
| 9-18 | 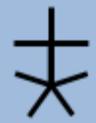 | 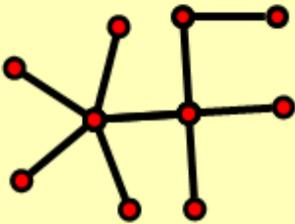 9-18-1 | 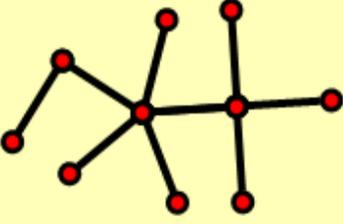 9-18-2 | 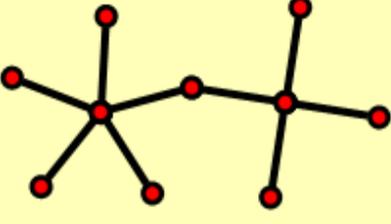 9-18-3 |

| | | |
|---|---|---|
| 9-19 | 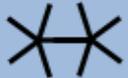 | 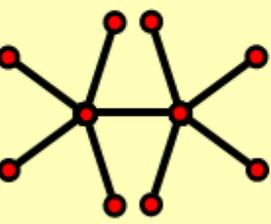 9-19-1 |

### 3.9.1.5. Matchstick graphs with $|E|=9$, $F=1$, $|V|=10$, $\Delta=6$

| | | | | |
|---|---|---|---|---|
| 9-20 | 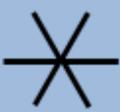 | 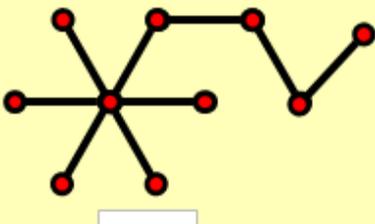 9-20-1 | 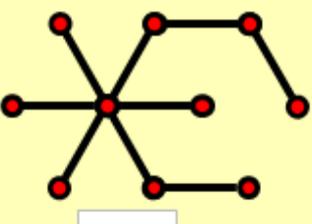 9-20-2 | 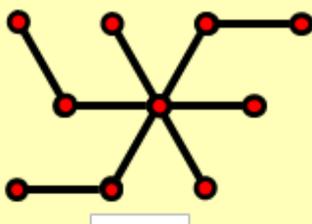 9-20-3 |
| 9-21 | 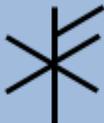 | 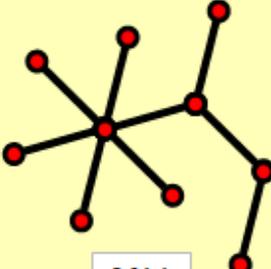 9-21-1 | 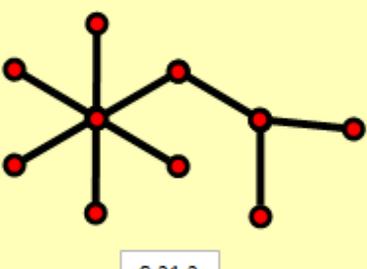 9-21-2 | 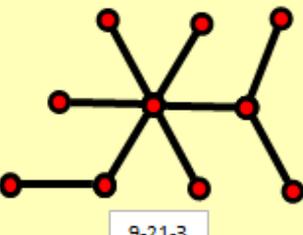 9-21-3 |

| | | |
|---|---|---|
| 9-22 | 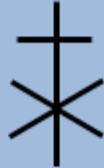 | 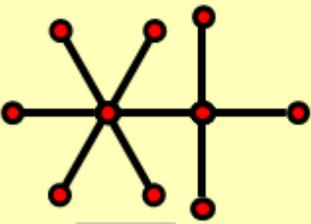 9-22-1 |

### 3.9.1.6. Matchstick graphs with |E|=9, F=1, |V|=10, Δ=7

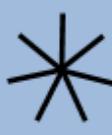

9-23

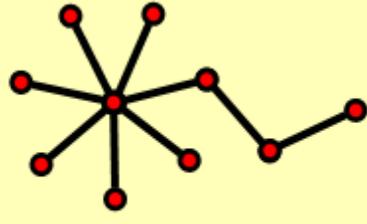

9-23-1

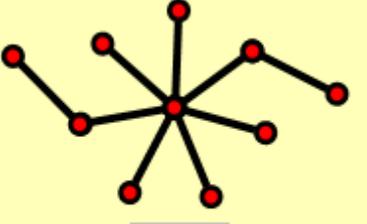

9-23-2

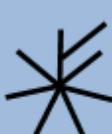

9-24

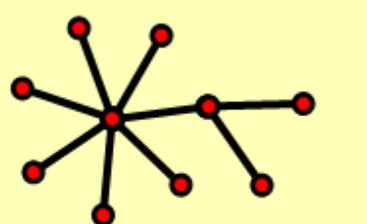

9-24-1

### 3.9.1.7. Matchstick graphs with |E|=9, F=1, |V|=10, Δ=8

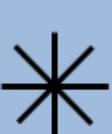

9-25

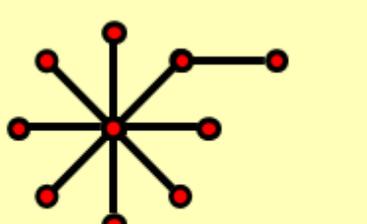

9-25-1

### 3.9.1.8. Matchstick graphs with |E|=9, F=1, |V|=10, Δ=9

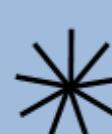

9-26

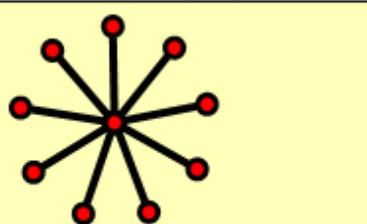

9-26-1

## 3.9.2. Matchstick graphs with |E|=9, F=2, |V|=9
### 3.9.2.1. Matchstick graphs with |E|=9, F=2, |V|=9, Δ=2

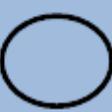

9-27

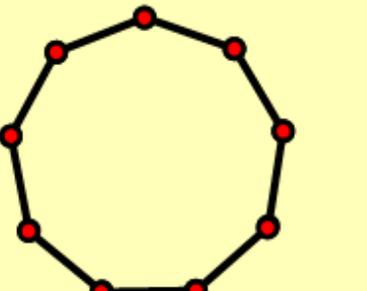

9-27-1

## 3.9.2.2. Matchstick graphs with |E|=9, F=2, |V|=9, Δ=3

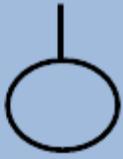

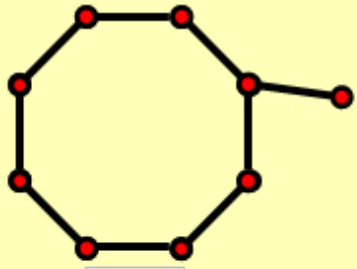

9-28-1

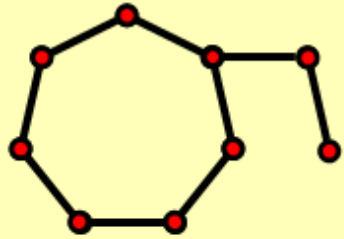

9-28-2

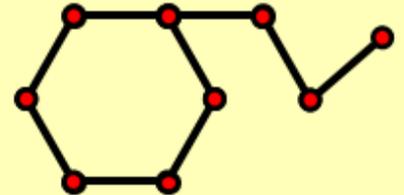

9-28-3

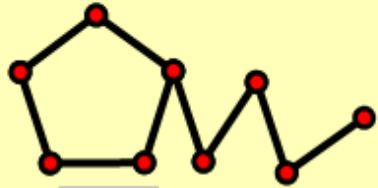

9-28-4

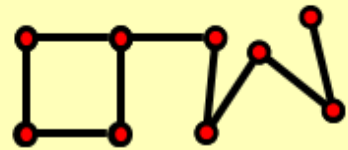

9-28-5

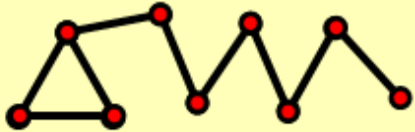

9-28-6

9-28

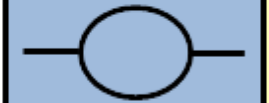

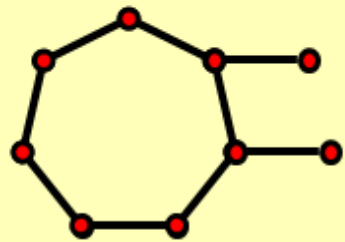

9-29-1

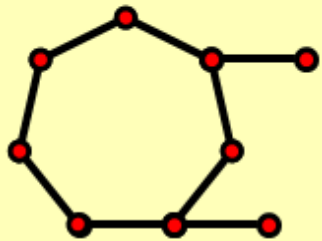

9-29-2

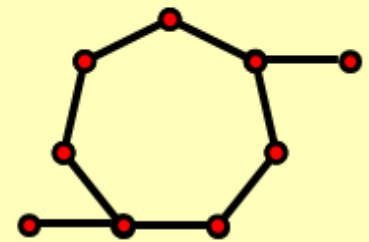

9-29-3

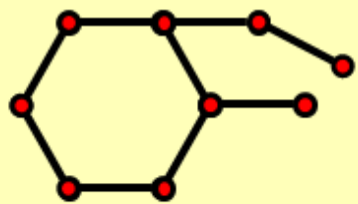

9-29-4

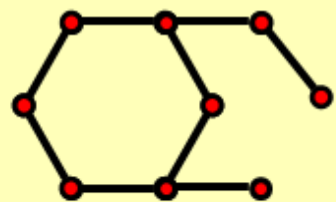

9-29-5

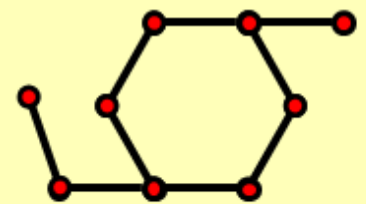

9-29-6

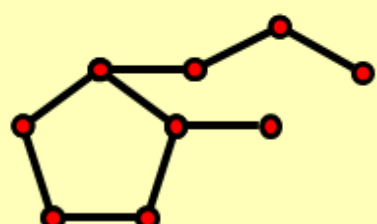

9-29-7

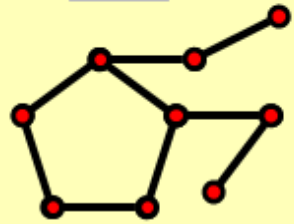

9-29-8

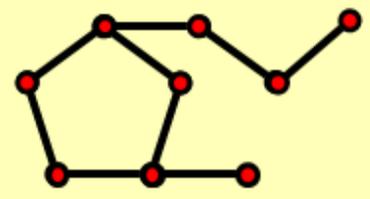

9-29-9

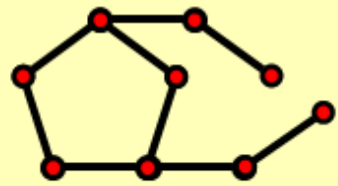
9-29-10

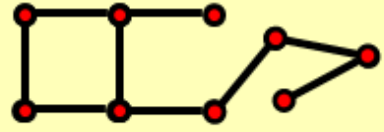
9-29-11

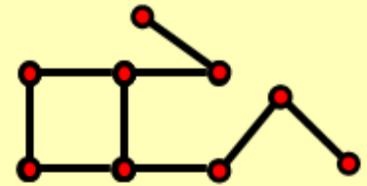
9-29-12

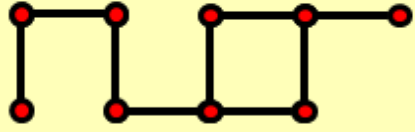
9-29-13

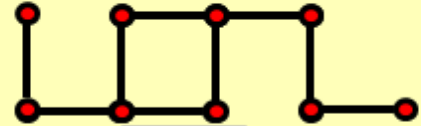
9-29-14

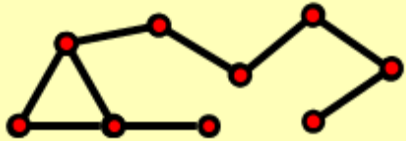
9-29-15

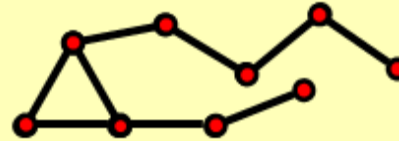
9-29-16

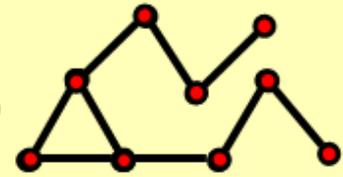
9-29-17

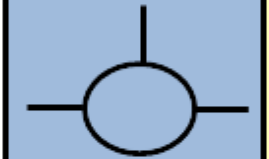

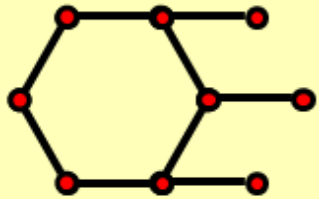
9-30-1

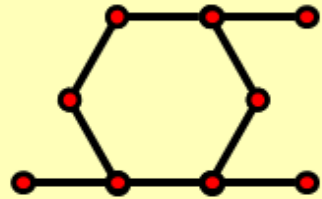
9-30-2

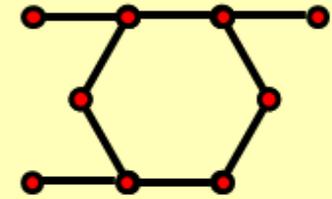
9-30-3

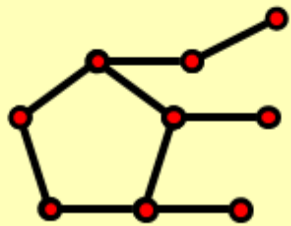
9-30-4

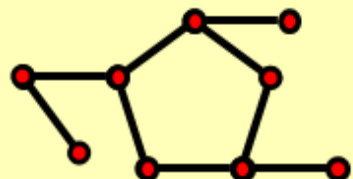
9-30-5

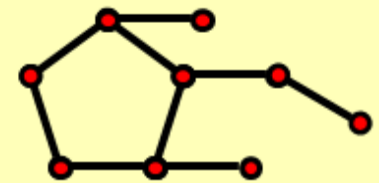
9-30-6

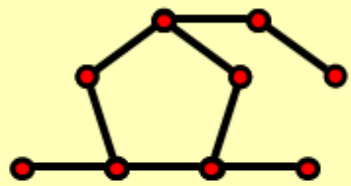
9-30-7

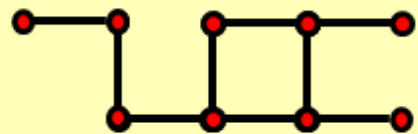
9-30-8

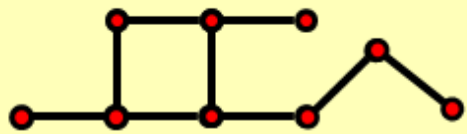
9-30-9

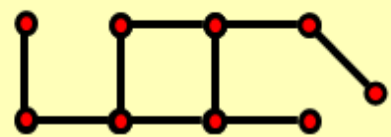
9-30-10

| 9-30 | | 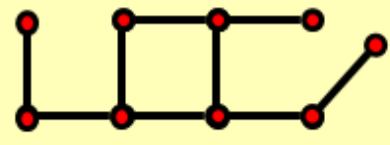 9-30-11 | 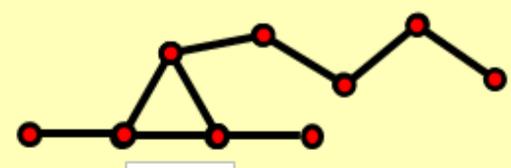 9-30-12 | 14 |
|---|---|---|---|---|
| | | 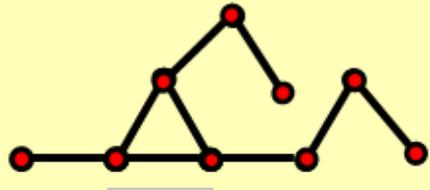 9-30-13 | 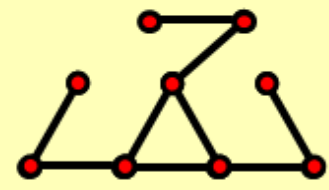 9-30-14 | |
| 9-31 | 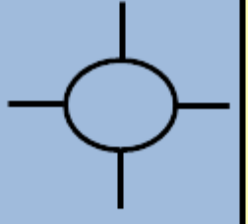 | 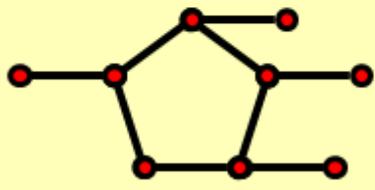 9-31-1 | 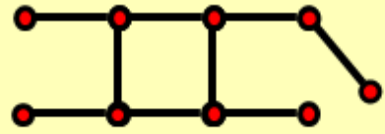 9-31-2 | 2 |
| | 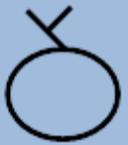 | 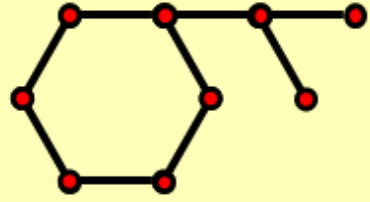 9-32-1 | 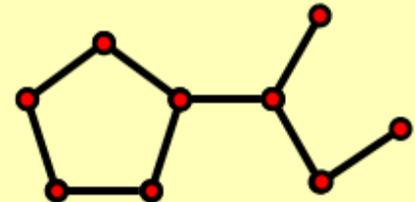 9-32-2 | |
| | | 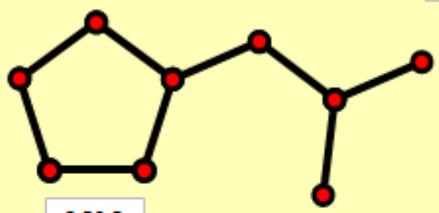 9-32-3 | 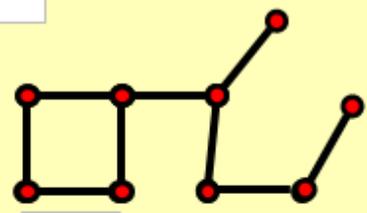 9-32-4 | |
| | | 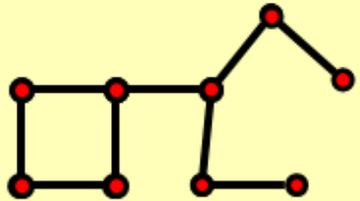 9-32-5 | 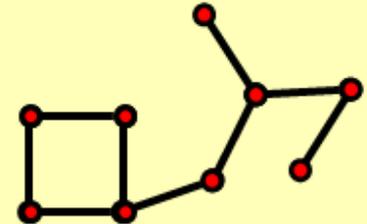 9-32-6 | |
| | | 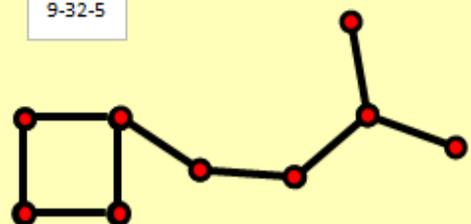 9-32-7 | 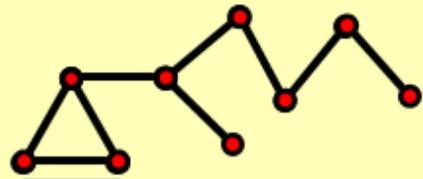 9-32-8 | |

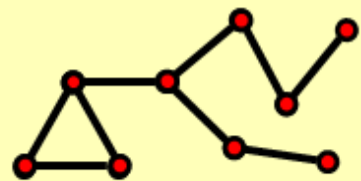

9-32-9

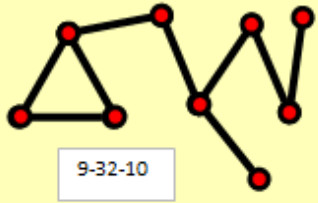

9-32-10

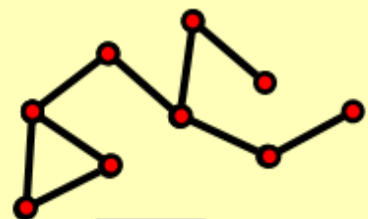

9-32-11

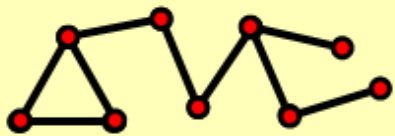

9-32-12

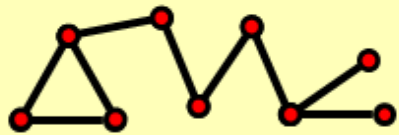

9-32-14

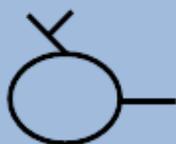

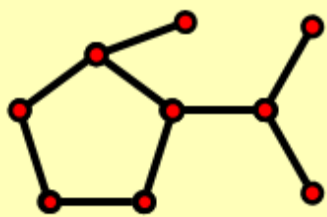

9-33-1

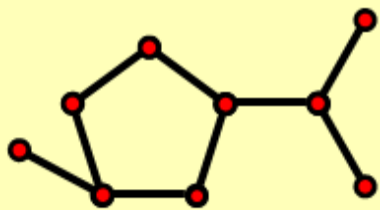

9-33-2

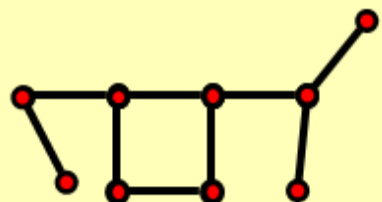

9-33-3

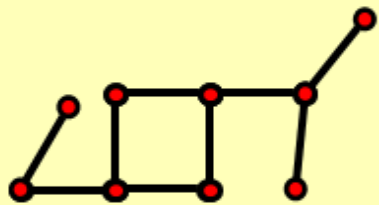

9-33-4

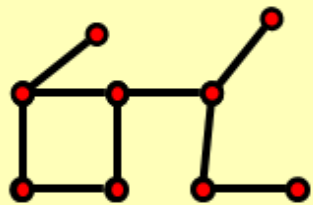

9-33-5

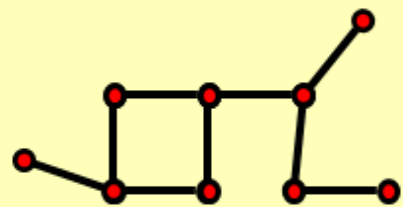

9-33-6

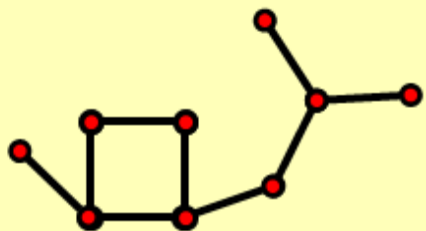

9-33-7

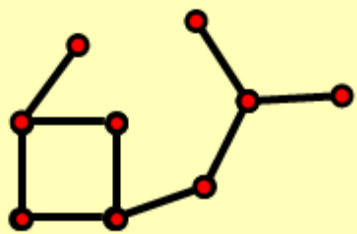

9-33-8

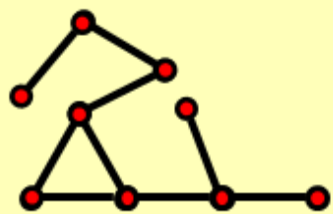

9-33-9

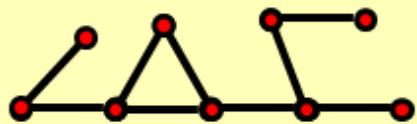

9-33-10

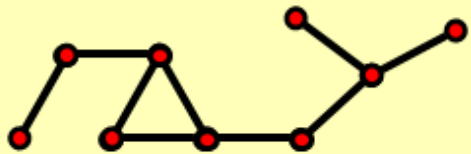

9-33-11

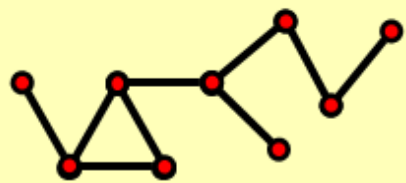

9-33-12

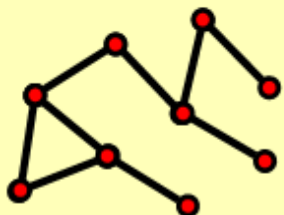

9-33-13

9-32

| | | |
|---|---|---|
| 9-33 | 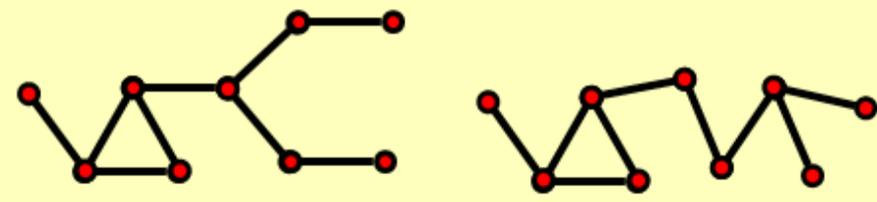 | 9-33-14, 9-33-15 — 15 |
| 9-34 | 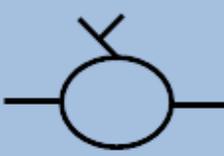 | 9-34-1, 9-34-2, 9-34-3, 9-34-4, 9-34-5 — 5 |
| 9-35 | 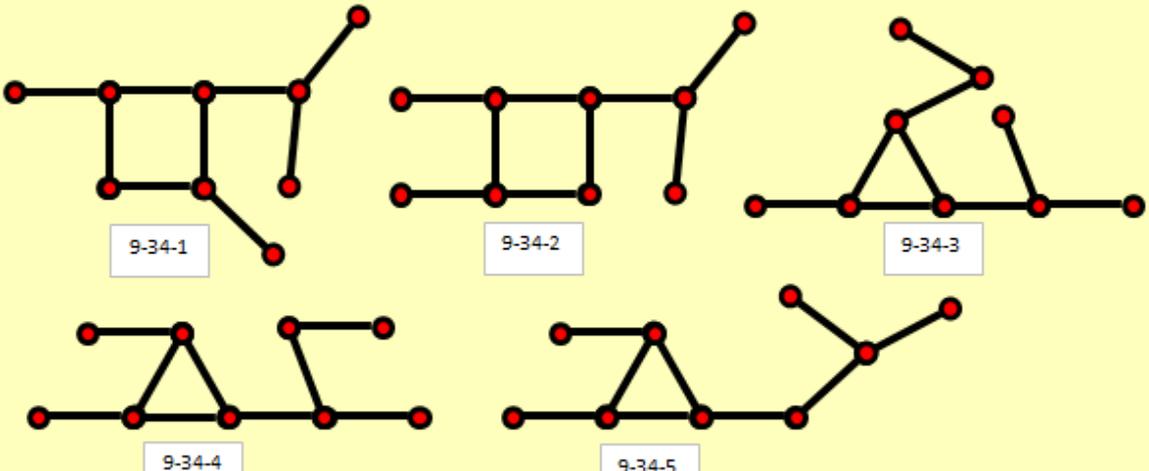 | 9-35-1 — 1 |
| 9-36 | 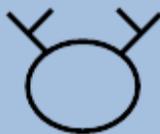 | 9-36-1, 9-36-2, 9-36-3, 9-36-4, 9-36-5 — 5 |
| 9-37 | 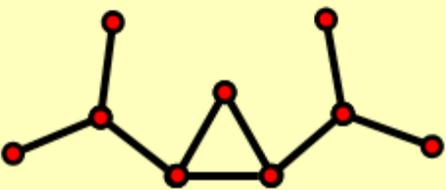 | 9-37-1 — 1 |

### 3.9.2.3. Matchstick graphs with |E|=9, *F*=2, |V|=9, Δ=4

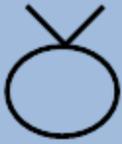

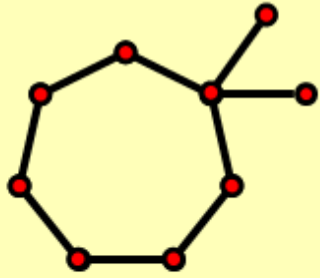
9-38-1

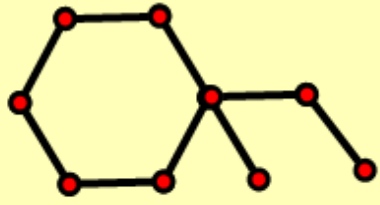
9-38-2

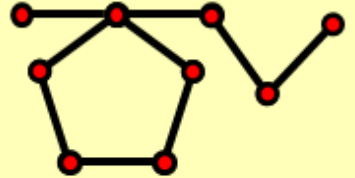
9-38-3

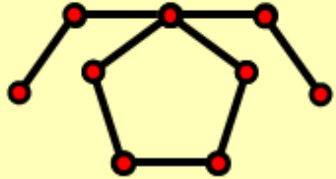
9-38-4

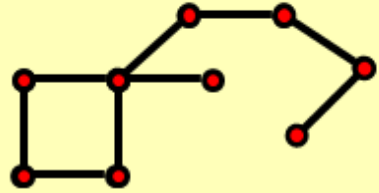
9-38-5

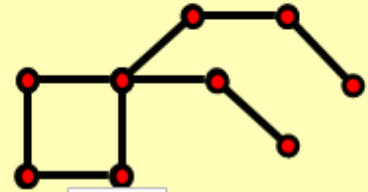
9-38-6

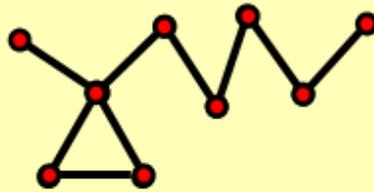
9-38-7

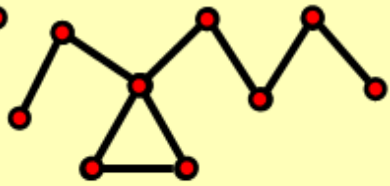
9-38-8

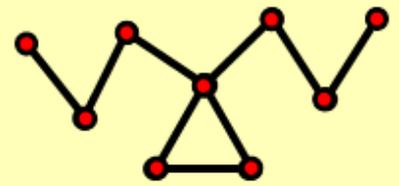
9-38-9

9-38

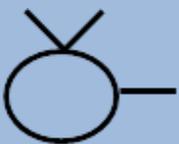

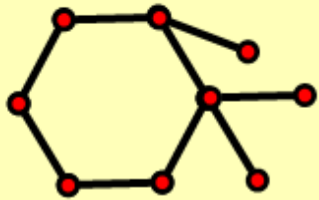
9-39-1

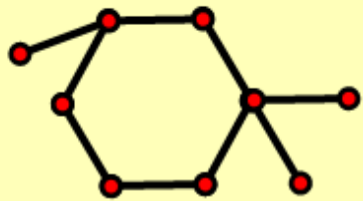
9-39-2

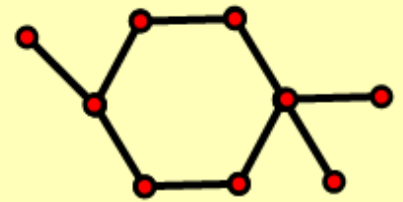
9-39-3

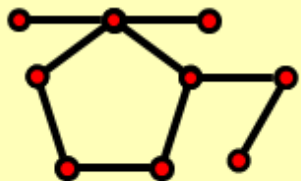
9-39-4

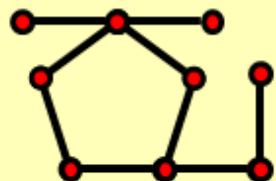
9-39-5

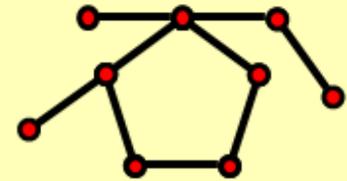
9-39-6

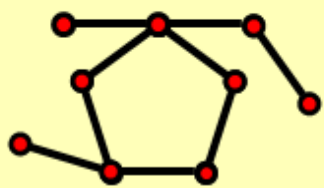
9-39-7

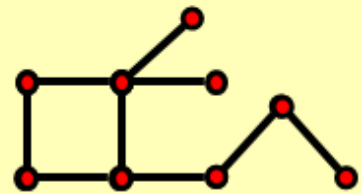
9-39-8

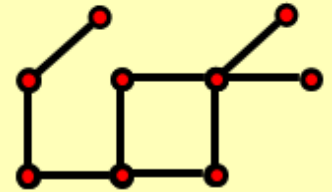
9-39-9

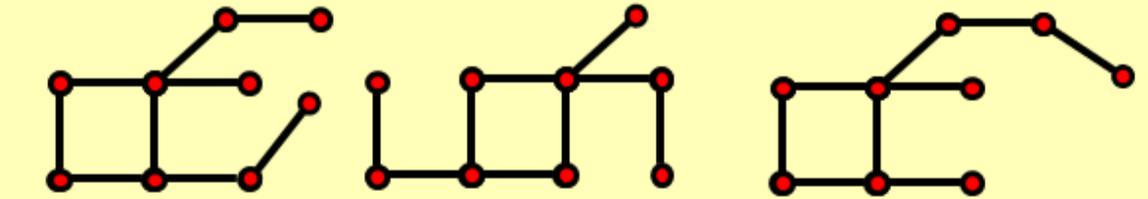

9-39-10   9-39-11   9-39-12

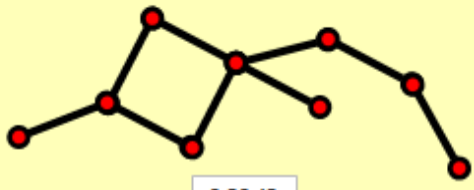
9-39-13

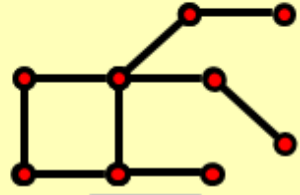
9-39-14

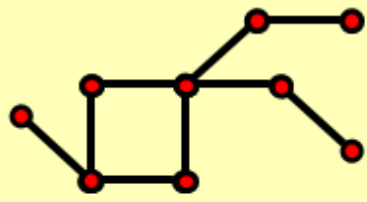
9-39-15

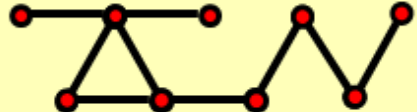
9-39-16

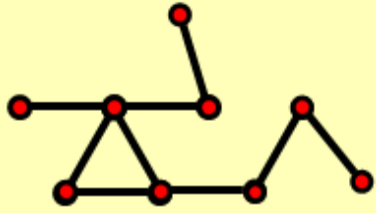
9-39-17

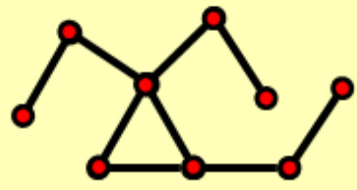
9-39-18

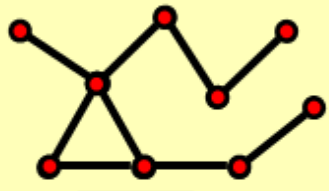
9-39-19

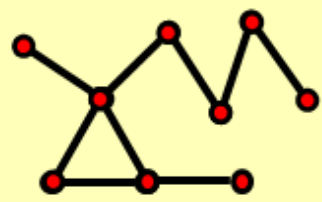
9-39-20

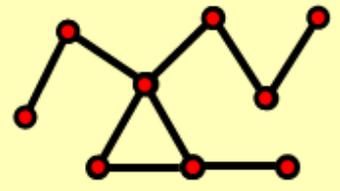
9-39-21

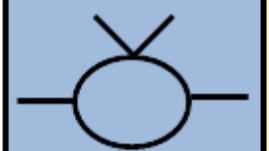

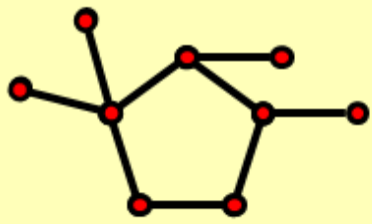
9-40-1

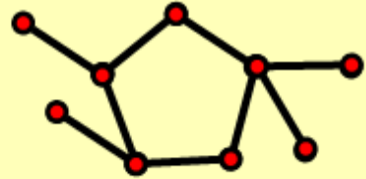
9-40-2

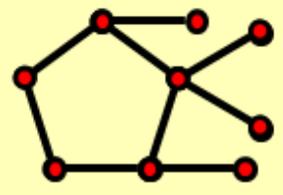
9-40-3

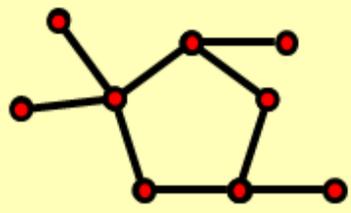
9-40-4

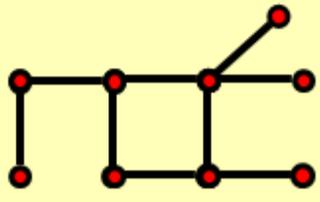
9-40-5

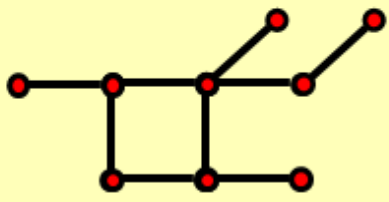
9-40-6

| 9-40 | 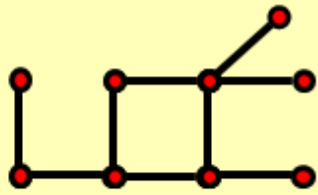 9-40-7 | 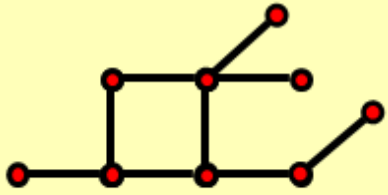 9-40-8 | 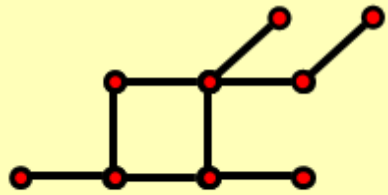 9-40-9 |
| --- | --- | --- | --- |
| | 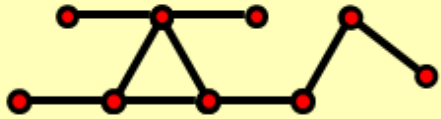 9-40-10 | 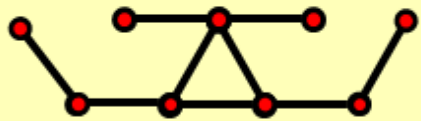 9-40-11 | |
| | 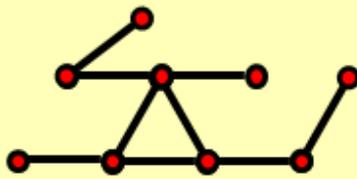 9-40-12 | 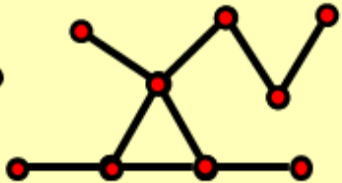 9-40-13 | 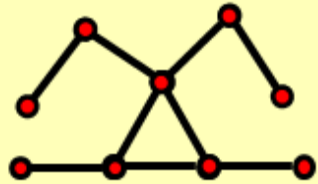 9-40-14 |

| 9-41 | 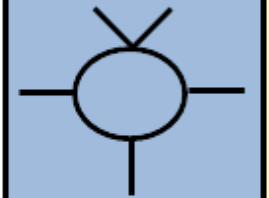 | 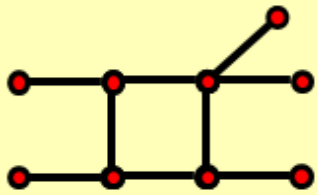 9-41-1 |
| --- | --- | --- |

| 9-42 | 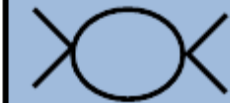 | 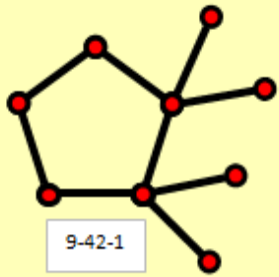 9-42-1 | 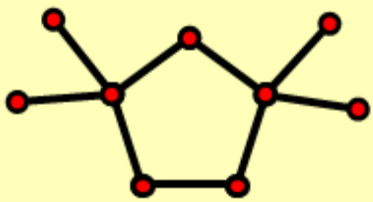 9-42-2 | 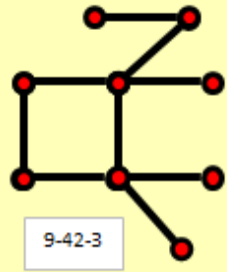 9-42-3 |
| --- | --- | --- | --- | --- |
| | | 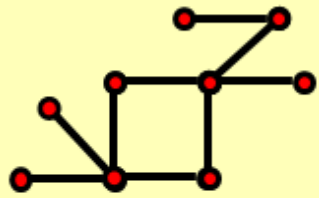 9-42-4 | 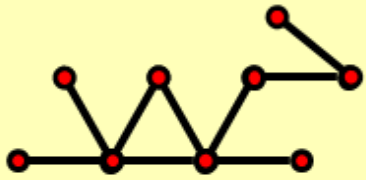 9-42-5 | 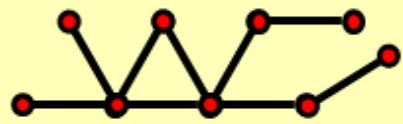 9-42-6 |
| | | 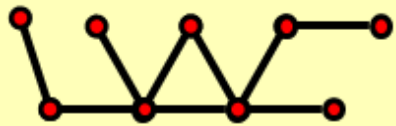 9-42-7 | | |

| | | |
|---|---|---|
| 9-43 | 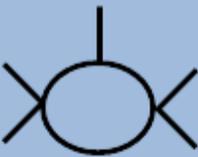 | 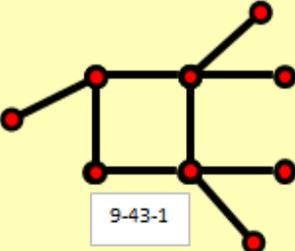 9-43-1  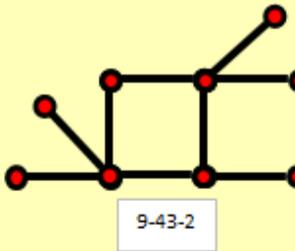 9-43-2  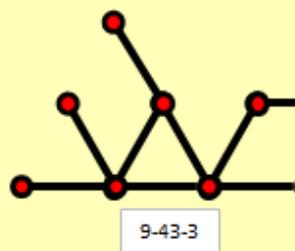 9-43-3 <br> 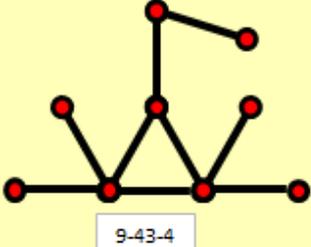 9-43-4 |
| 9-44 | 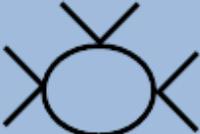 | 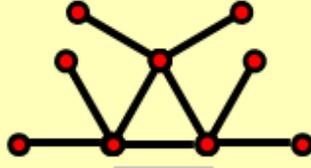 9-44-1 |
| 9-45 | 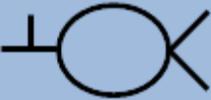 | 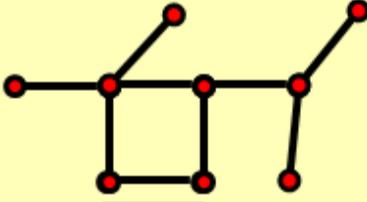 9-45-1  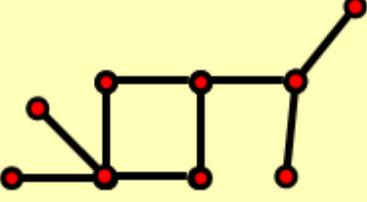 9-45-2  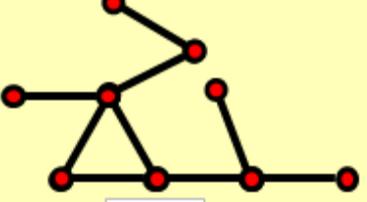 9-45-3 <br> 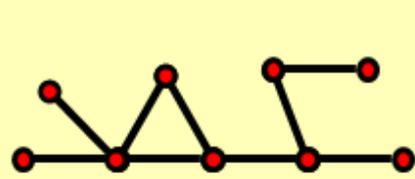 9-45-4  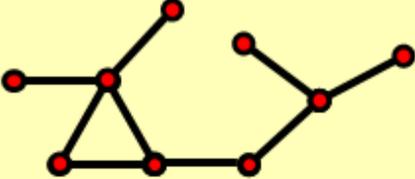 9-45-5 |
| 9-46 | 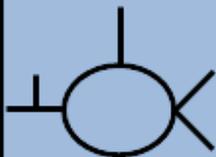 | 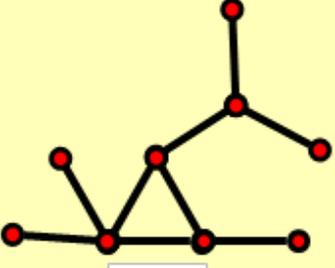 9-46-1 |
| | 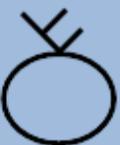 | 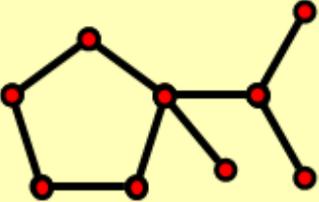 9-47-1  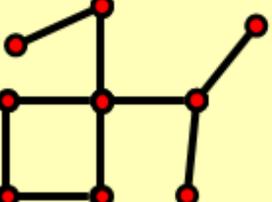 9-47-2  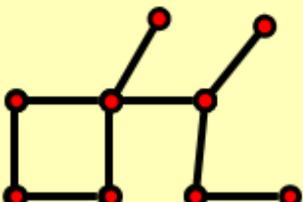 9-47-3 |
| 9-47 | 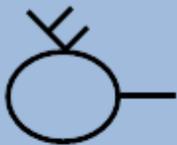 | 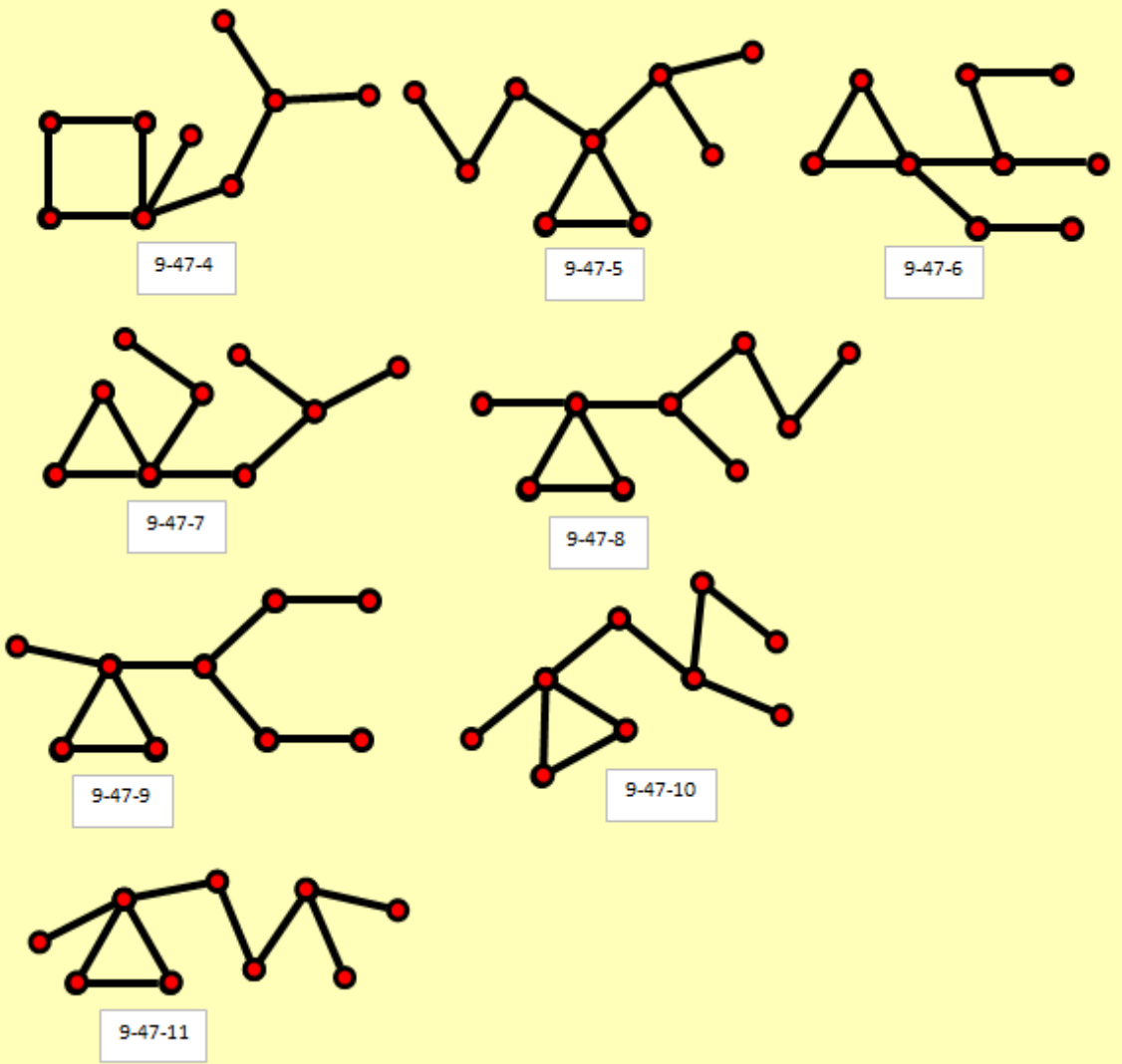 |
| 9-48 | | 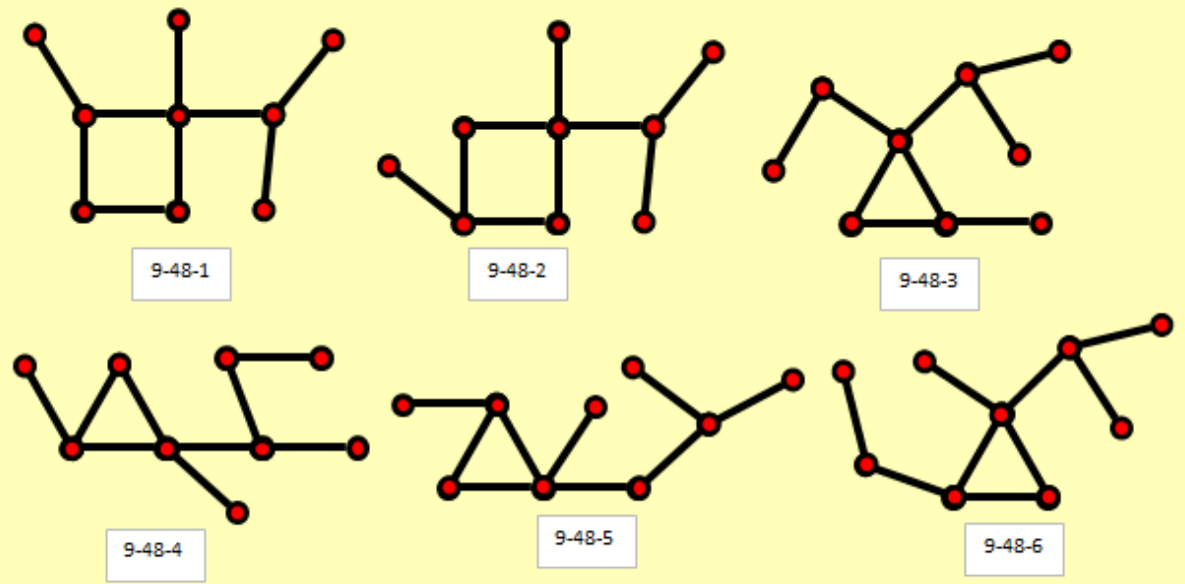 |
| 9-49 | 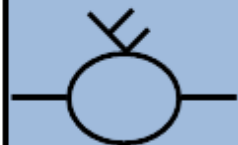 | 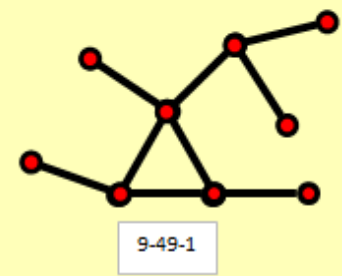 |
| 9-50 | 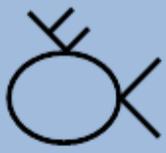 | 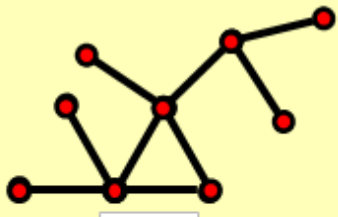 9-50-1 |
| 9-51 | 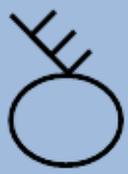 | 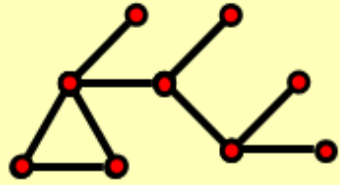 9-51-1 |
| 9-52 | 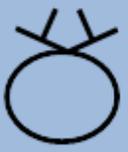 | 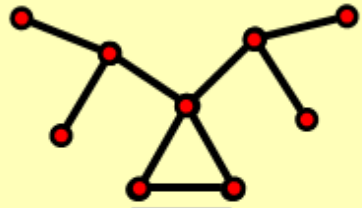 9-52-1 |
| 9-53 | 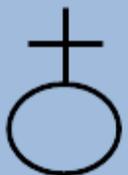 | 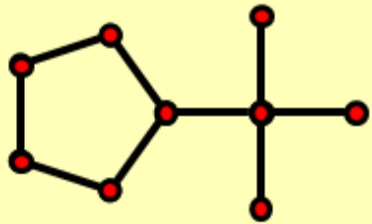 9-53-1 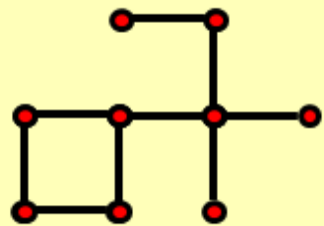 9-53-2 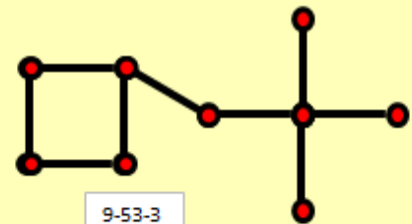 9-53-3 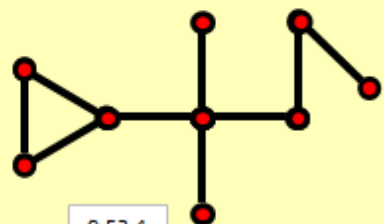 9-53-4 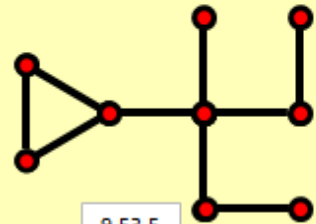 9-53-5 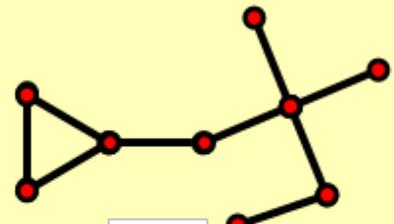 9-53-6 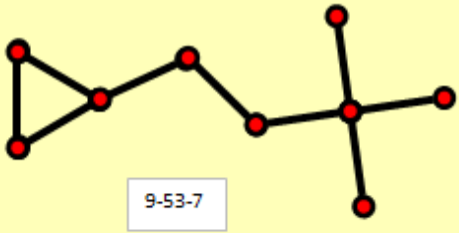 9-53-7 |
| 9-54 | | 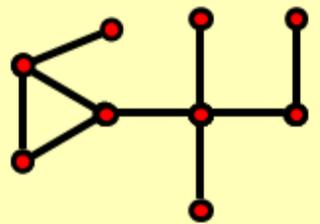 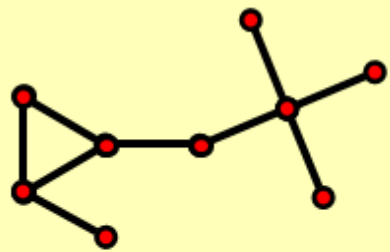 9-54-4    9-54-5   5 |
| 9-55 | 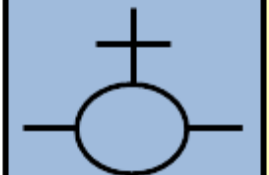 | 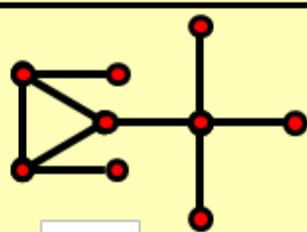 9-55-1   1 |
| 9-56 | 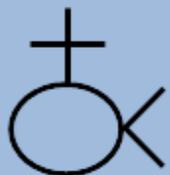 | 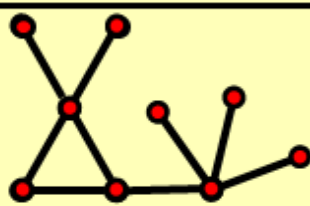 9-56-1   1 |
| 9-57 | 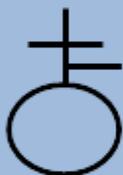 | 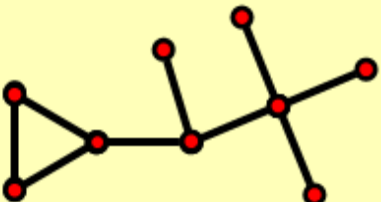 9-57-1   1 |
| 9-58 | 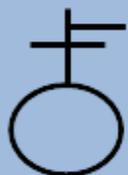 | 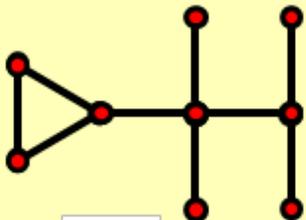 9-58-1   1 |
| 9-59 | 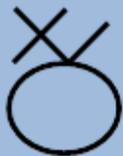 | 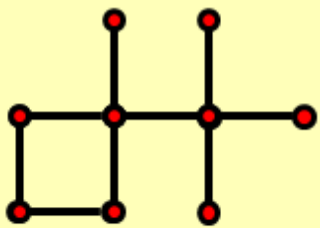 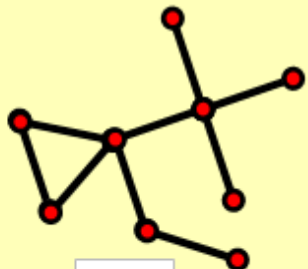 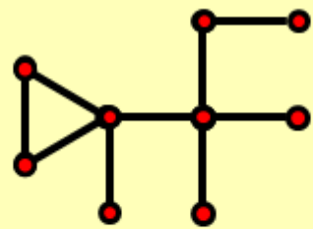 9-59-1   9-59-2   9-59-3 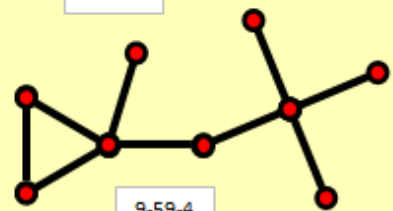 9-59-4   4 |

| 9-60 | 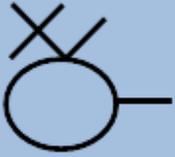 | 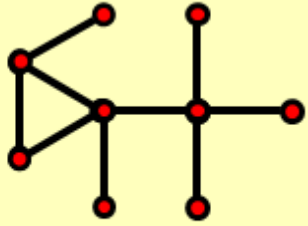  9-60-1 | 1 |

### 3.9.2.4. Matchstick graphs with |E|=9, $F$=2, |V|=9, $\Delta$=5

| 9-61 | 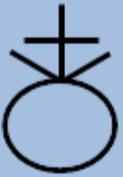 | 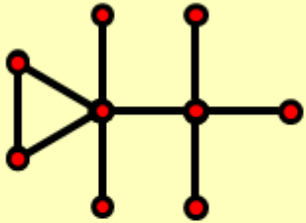  9-61-1 | 1 |

| 9-62 | 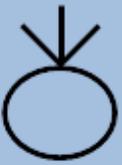 | 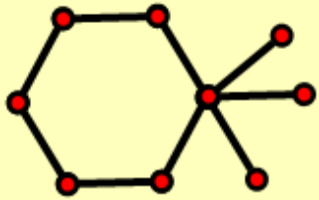 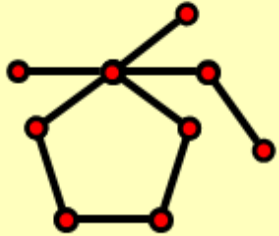 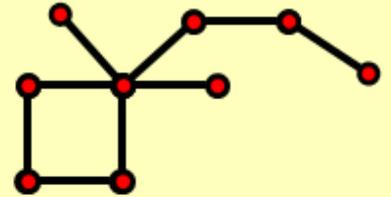  9-62-1  9-62-2  9-62-3  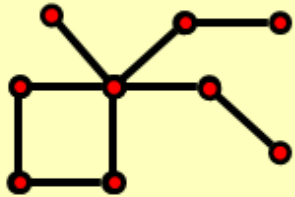 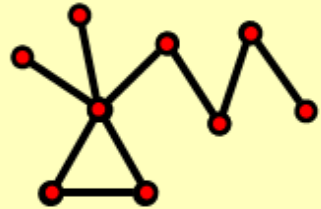 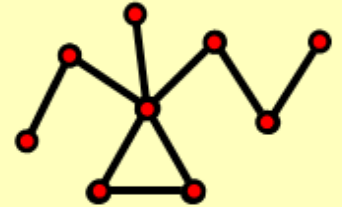  9-62-4  9-62-5  9-62-6  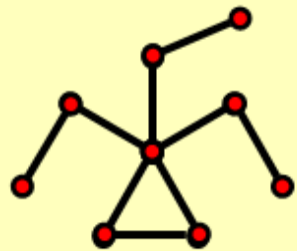  9-62-7 | 7 |

| 9-63 | 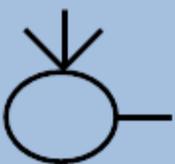 | 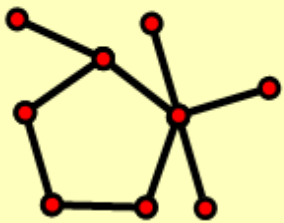 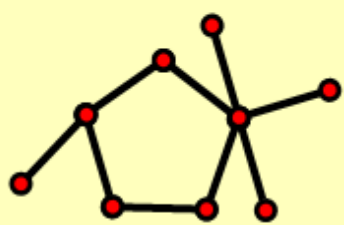 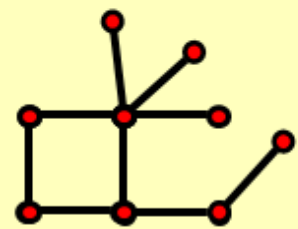  9-63-1  9-63-2  9-63-3 | |

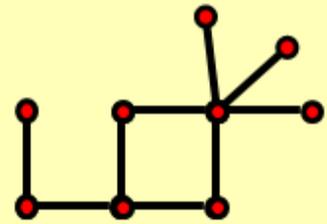
9-63-4

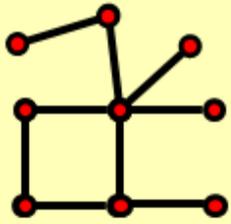
9-63-5

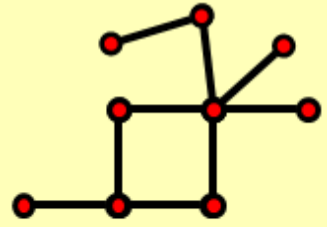
9-63-6

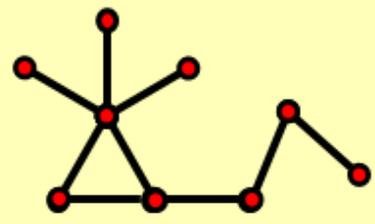
9-63-7

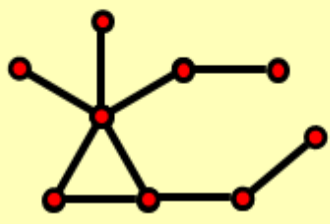
9-63-8

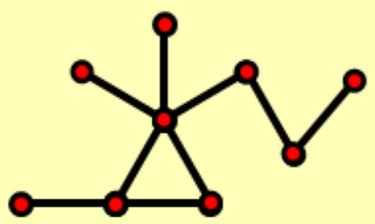
9-63-9

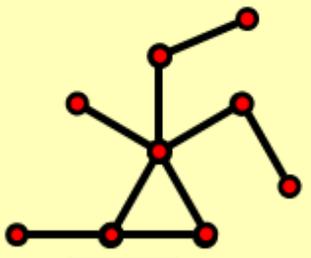
9-63-10

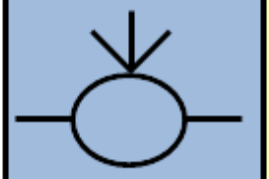

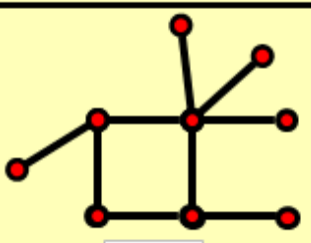
9-64-1

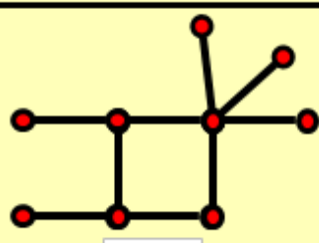
9-64-2

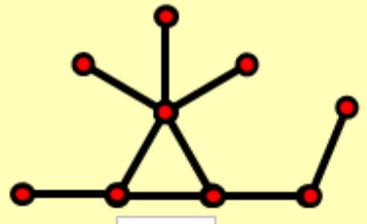
9-64-3

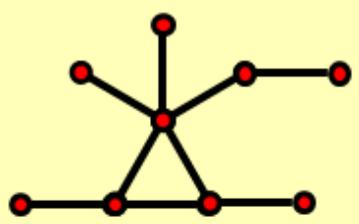
9-64-4

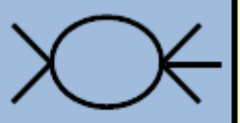

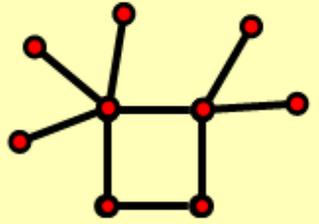
9-65-1

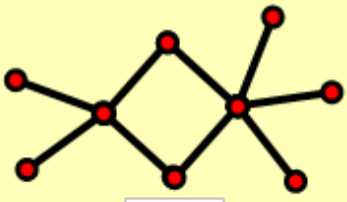
9-65-2

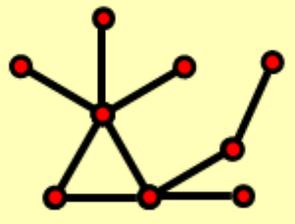
9-65-3

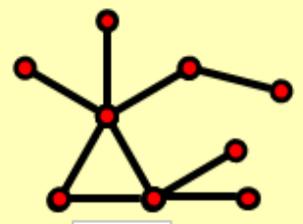
9-65-4

| | | |
|---|---|---|
| 9-66 | 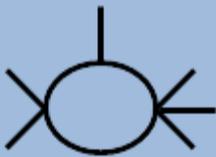 | 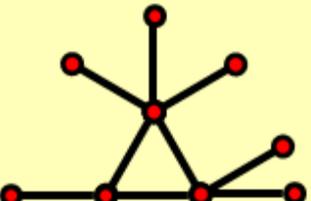<br>9-66-1 |
| 9-67 | 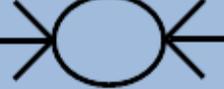 | 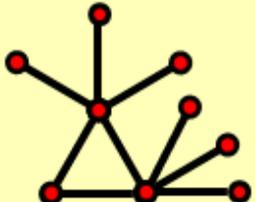<br>9-67-1 |
| 9-68 | 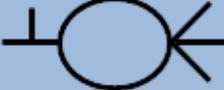 | 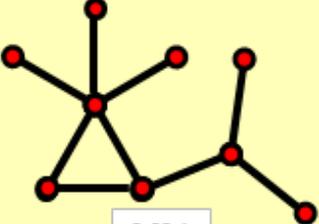<br>9-68-1 |
| 9-69 | 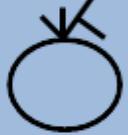 | 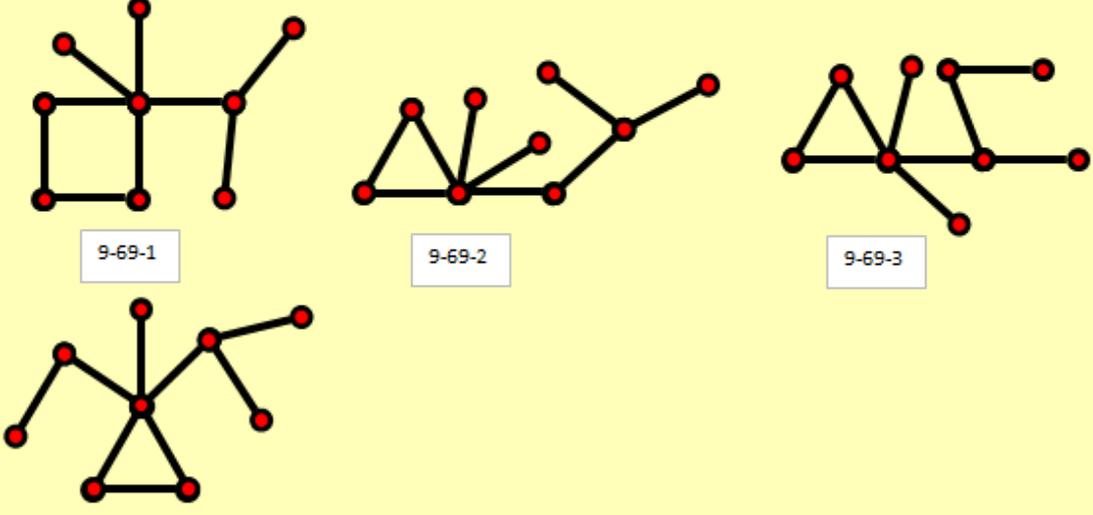<br>9-69-1, 9-69-2, 9-69-3, 9-69-4 |
| 9-70 | 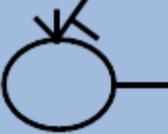 | 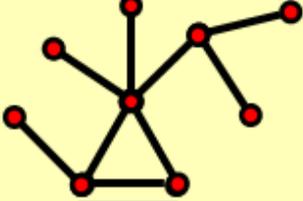<br>9-70-1 |
| 9-71 | 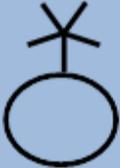 | 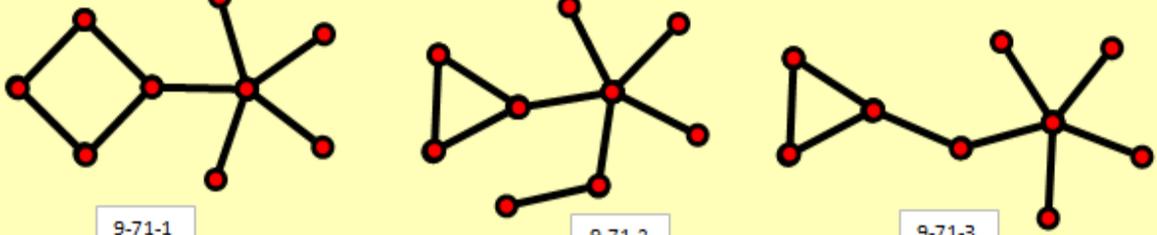<br>9-71-1, 9-71-2, 9-71-3 |

| 9-72 | 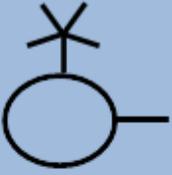 | 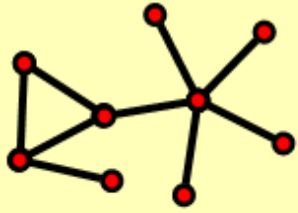

9-72-1 | 1 |

| 9-73 | 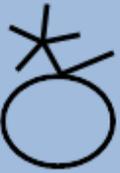 | 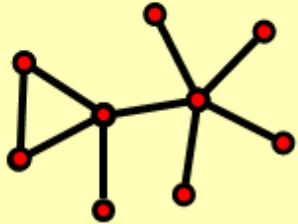

9-73-1 | 1 |

### 3.9.2.5. Matchstick graphs with |E|=9, $F$=2, |V|=9, $\Delta$=6

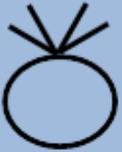

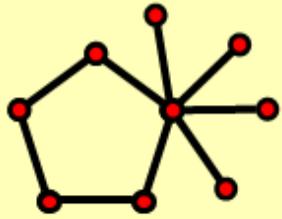

9-74-1

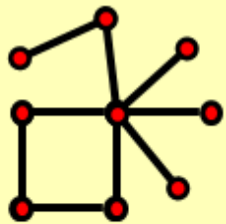

9-74-2

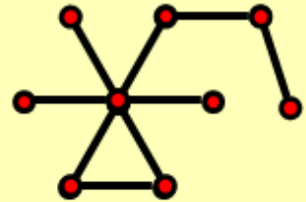

9-74-3

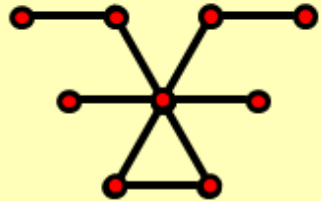

9-74-4

| 9-74 | | | 4 |

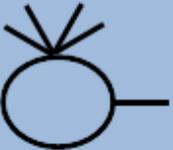

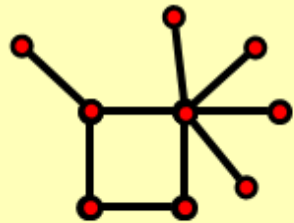

9-75-1

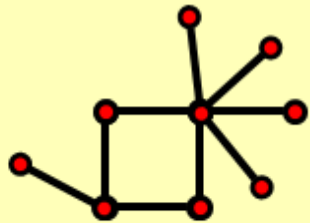

9-75-2

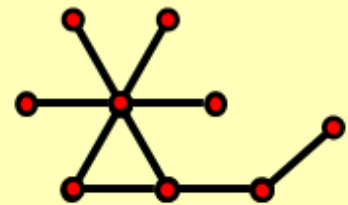

9-75-3

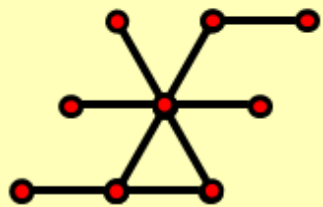

9-75-4

| 9-75 | | | 4 |

| 9-76 | 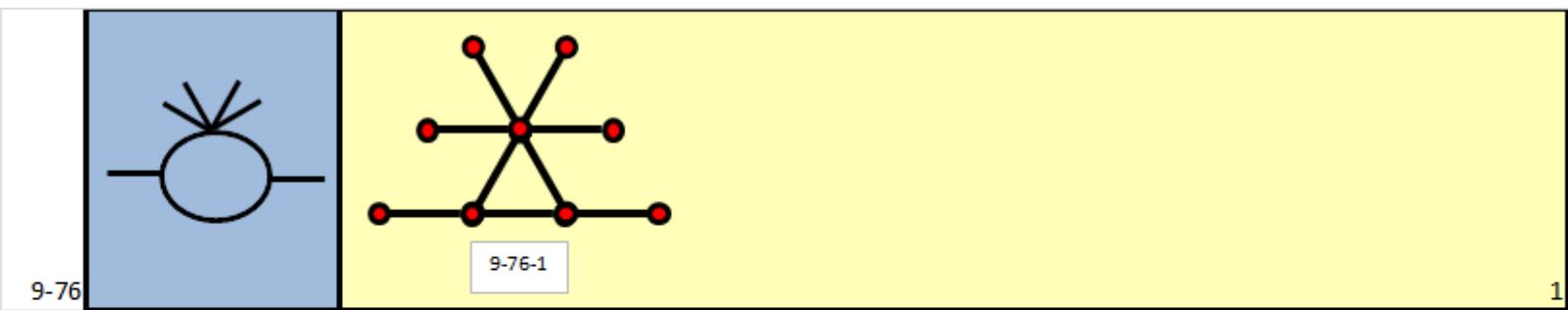 | | 1 |
|---|---|---|---|
| | 9-76-1 | | |

| 9-77 | 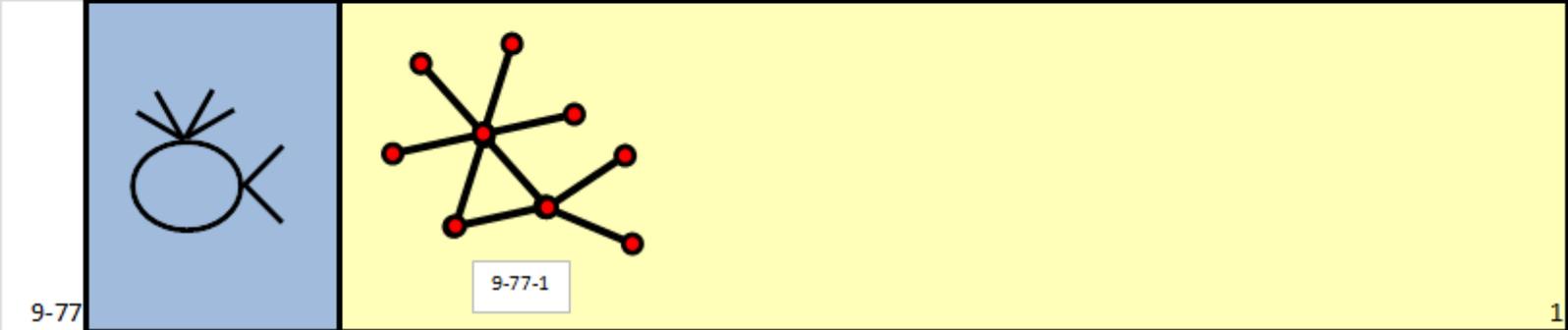 | | 1 |
|---|---|---|---|
| | 9-77-1 | | |

| 9-78 | 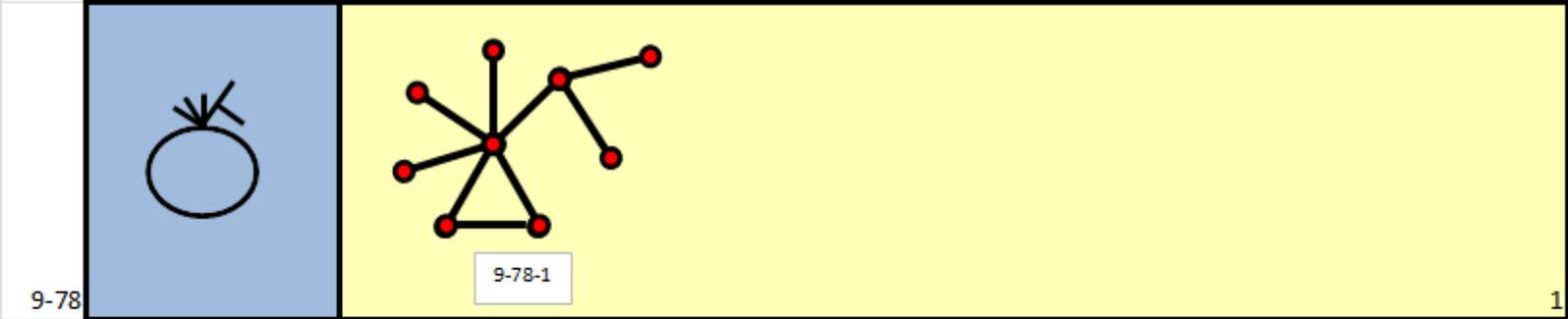 | | 1 |
|---|---|---|---|
| | 9-78-1 | | |

| 9-79 | 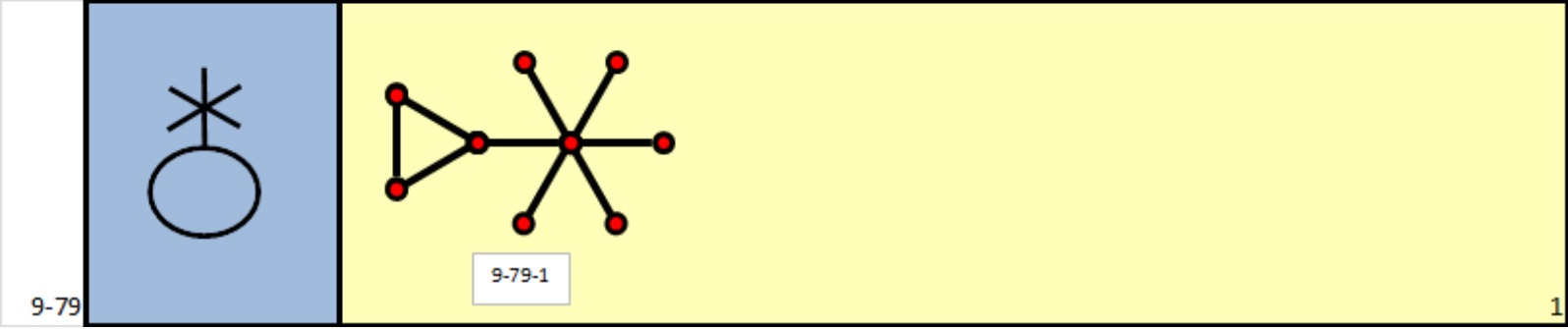 | | 1 |
|---|---|---|---|
| | 9-79-1 | | |

3.9.2.6. Matchstick graphs with $|E|=9$, $F=2$, $|V|=9$, $\Delta=7$

| 9-80 | 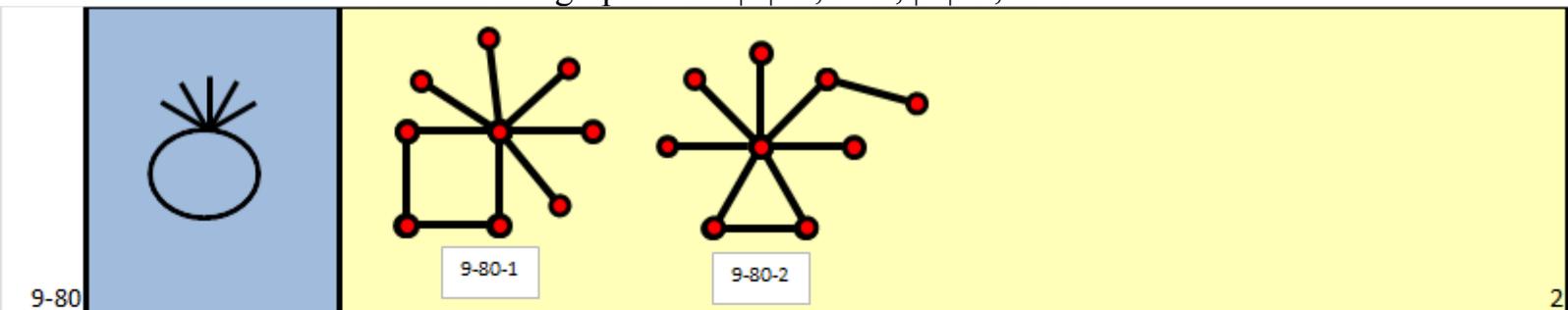 | | 2 |
|---|---|---|---|
| | 9-80-1 | 9-80-2 | |

| 9-81 | 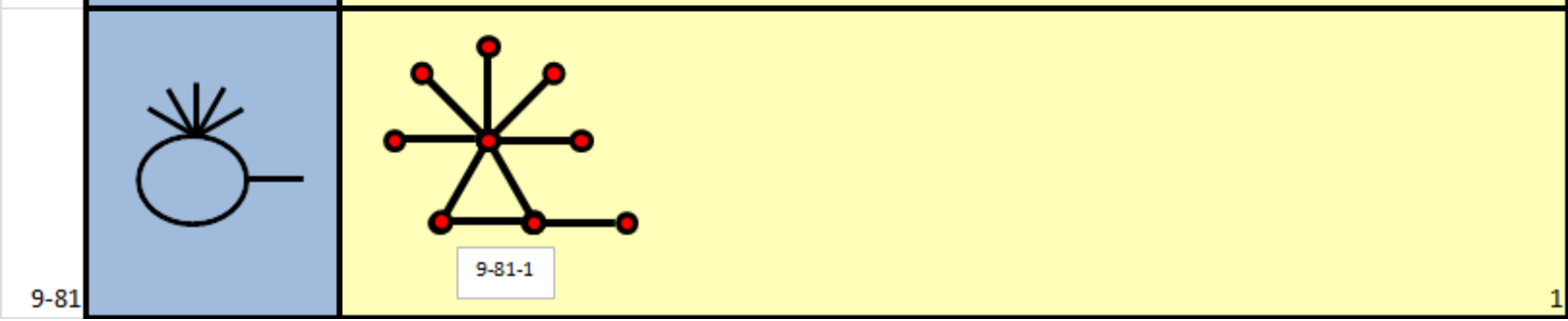 | | 1 |
|---|---|---|---|
| | 9-81-1 | | |

### 3.9.2.7. Matchstick graphs with |E|=9, F=2, |V|=9, Δ=8

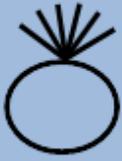
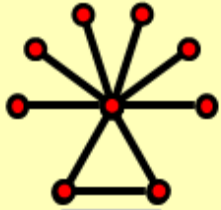

9-82-1

9-82 | | 1

## 3.9.3. Matchstick graphs with |E|=9, F=3, |V|=8
### 3.9.3.1. Matchstick graphs with |E|=9, F=3, |V|=8, Δ=3

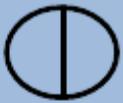
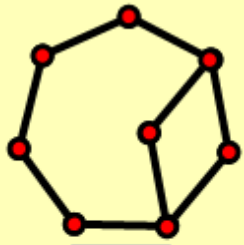
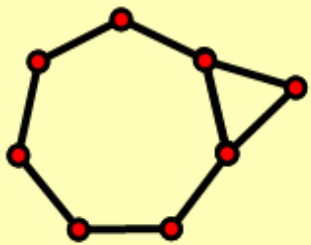
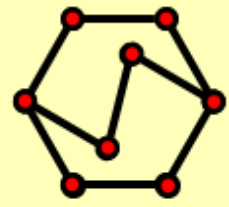

9-83-1     9-83-2     9-83-3

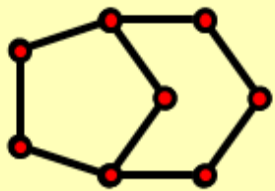
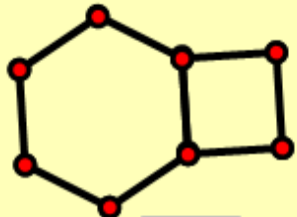
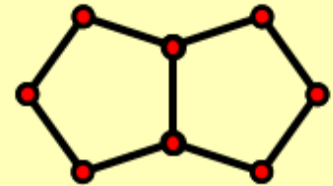

9-83-4     9-83-5     9-83-6

9-83 | | 6

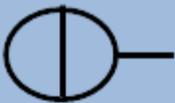
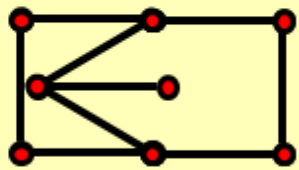
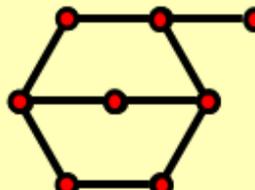
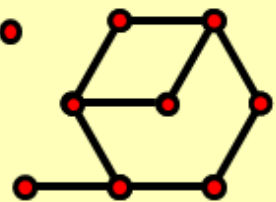
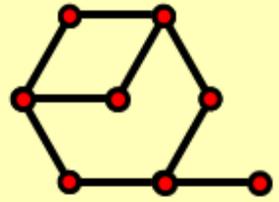

9-84-1     9-84-2     9-84-3     9-84-4

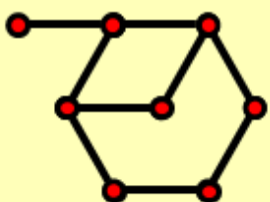
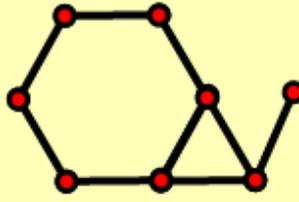
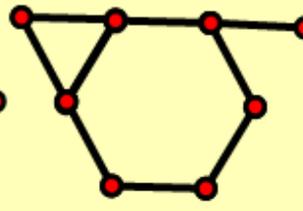
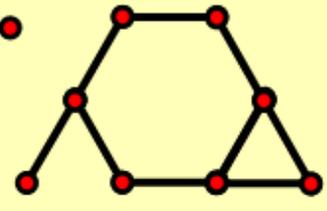

9-84-5     9-84-6     9-84-7     9-84-8

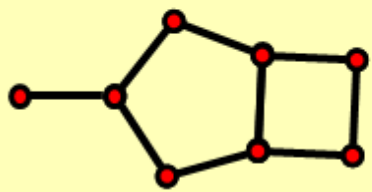
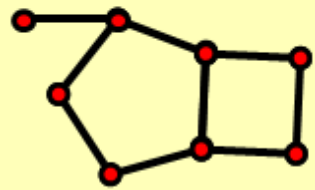
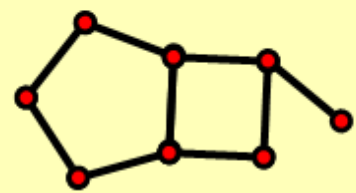

9-84-9     9-84-10     9-84-11

9-84

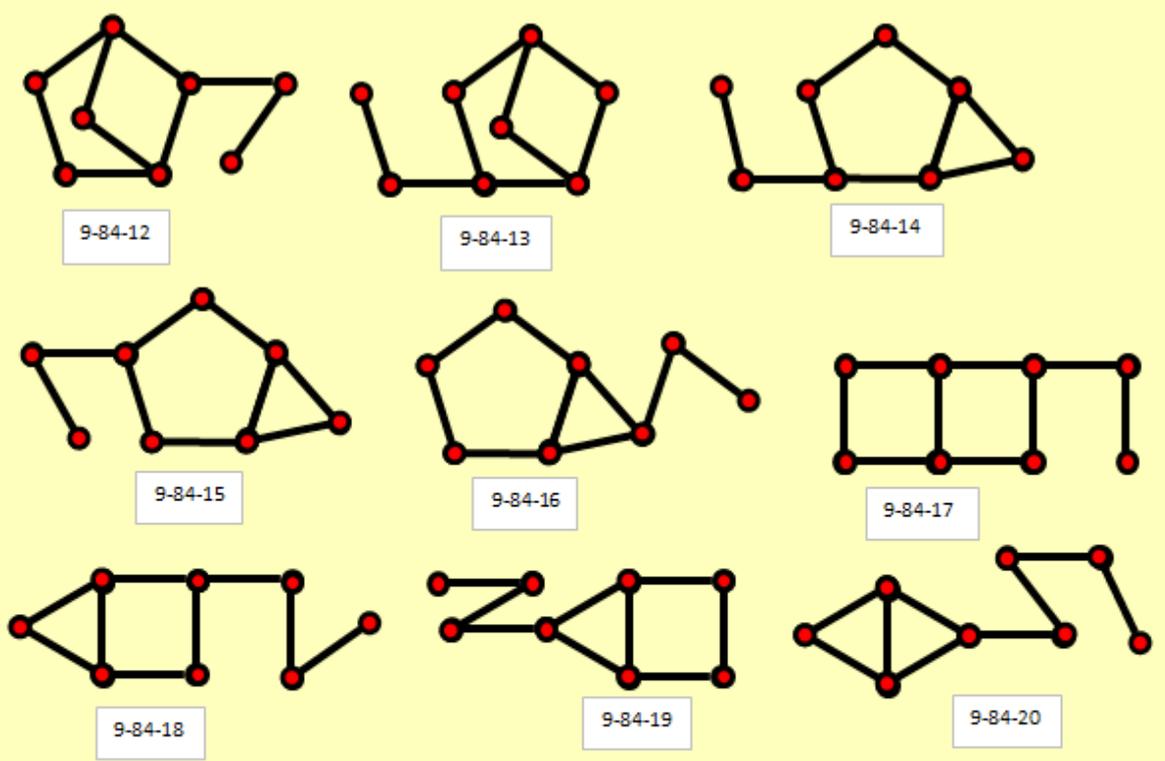

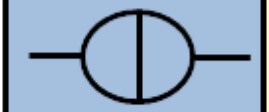

9-85

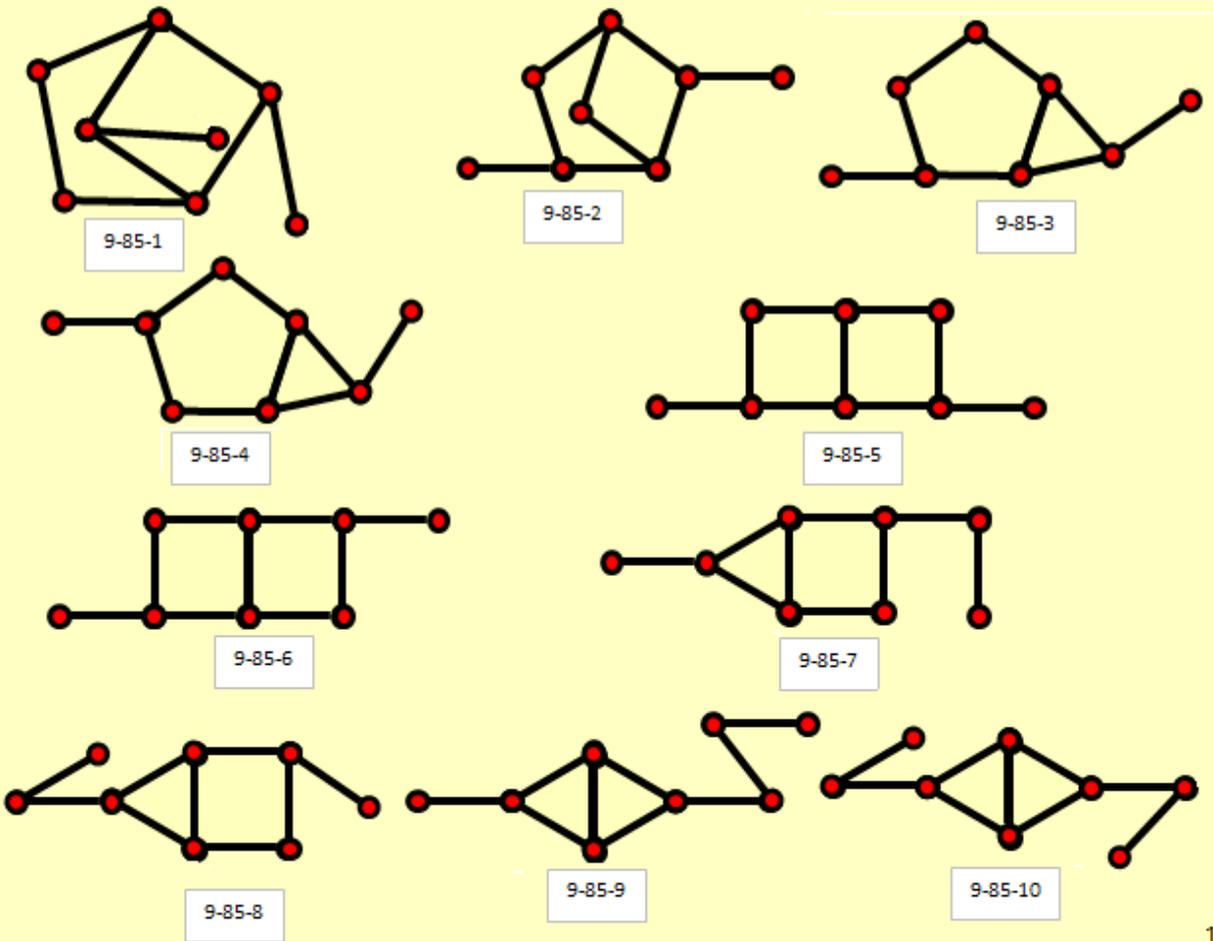

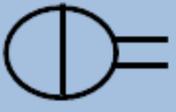

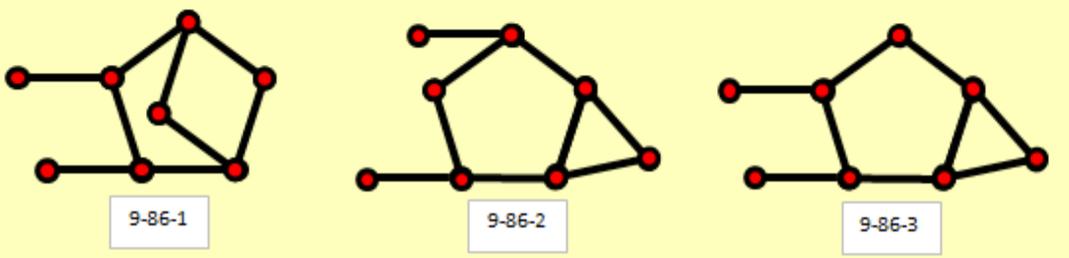

| 9-86 | 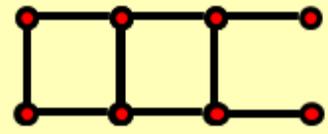 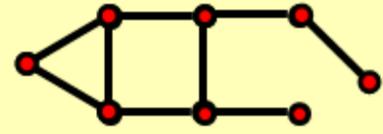 | 5 |
|---|---|---|
| 9-87 | 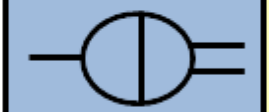 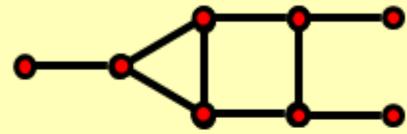 | 1 |
| 9-88 | 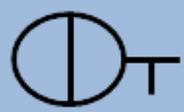 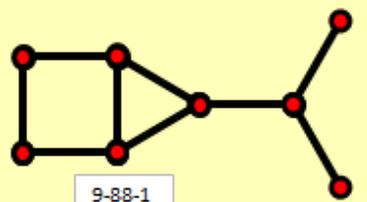 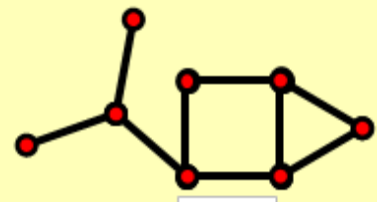 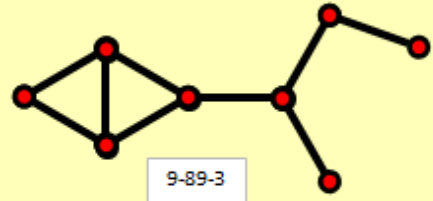 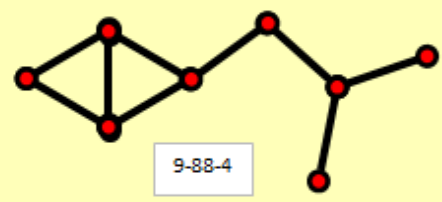 | 4 |
| 9-89 | 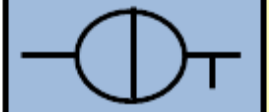 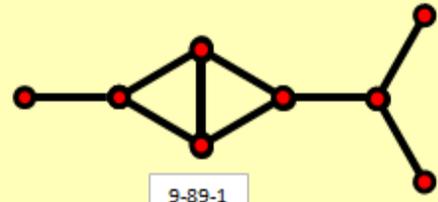 | 1 |
| 9-90 | 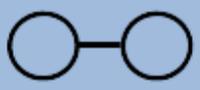 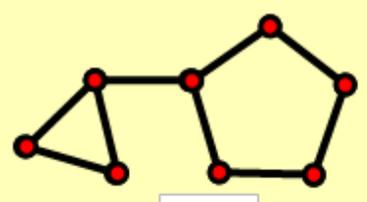 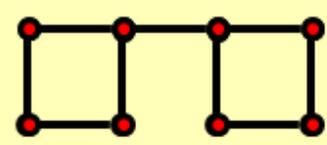 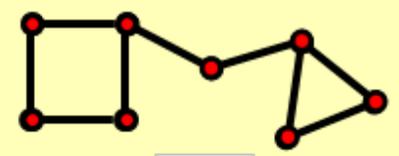 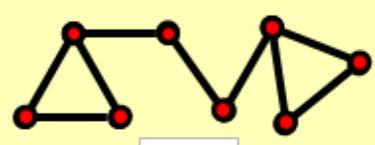 | 4 |
| 9-91 | 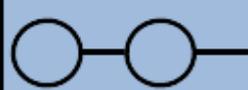 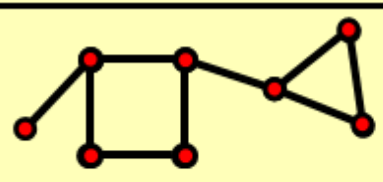 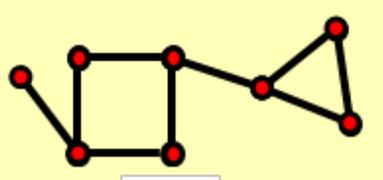 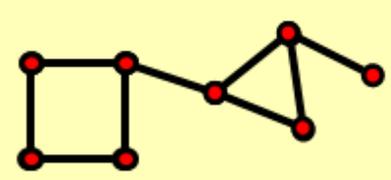 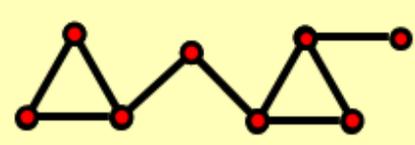 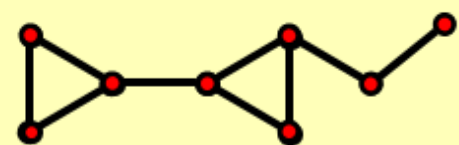 | 5 |

| 9-92 | 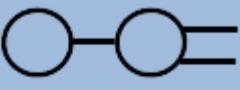 | 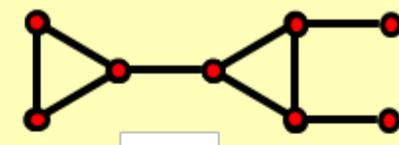 9-92-1 | 1 |
| 9-93 | 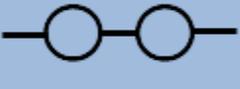 | 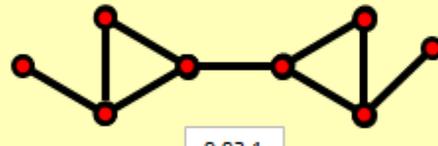 9-93-1 | 1 |
| 9-94 | 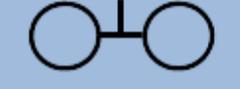 | 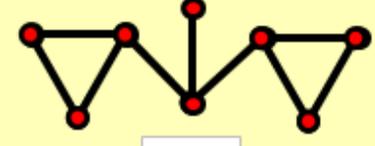 9-94-1 | 1 |

### 3.9.3.2. Matchstick graphs with |E|=9, F=3, |V|=8, Δ=4

| 9-95 | 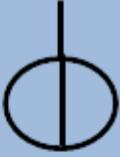 | 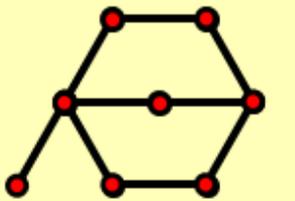 9-95-1   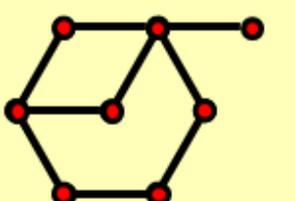 9-95-2   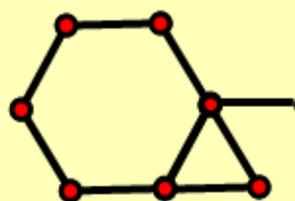 9-95-3 <br> 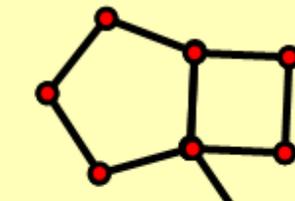 9-95-4   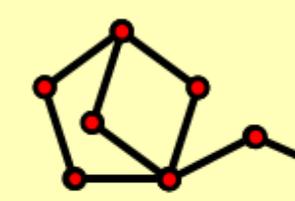 9-95-5   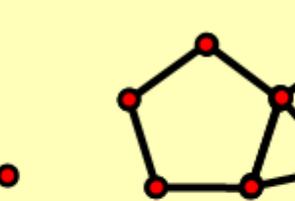 9-95-6 <br> 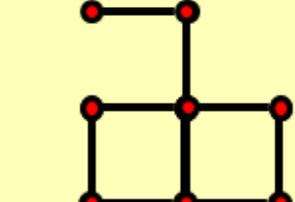 9-95-7   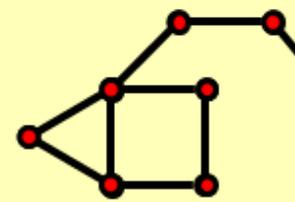 9-95-8   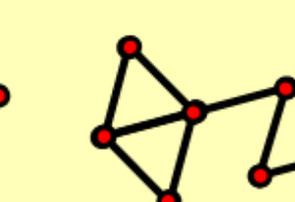 9-95-9 | 9 |
| 9-96 | 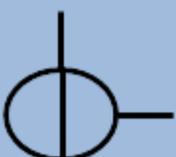 | 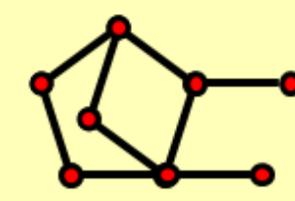 9-96-1   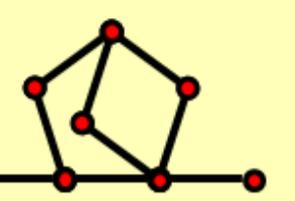 9-96-2   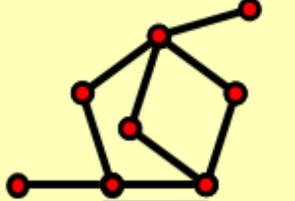 9-96-3 | |

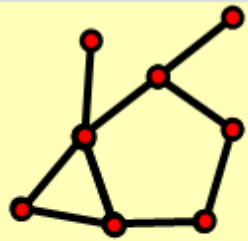
9-96-4

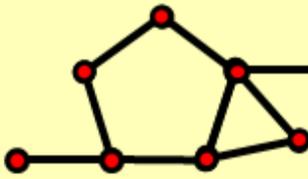
9-96-5

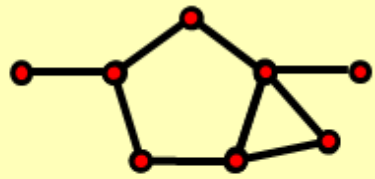
9-96-6

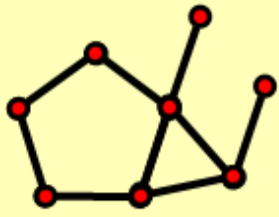
9-96-7

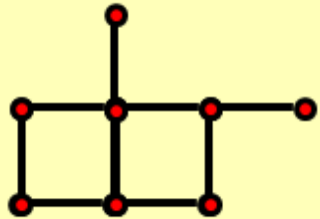
9-96-8

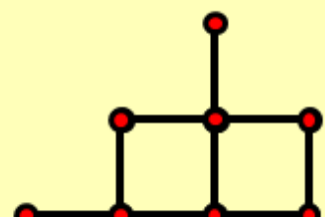
9-96-9

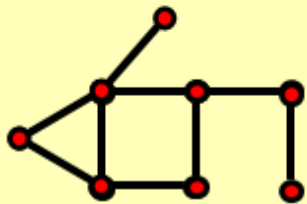
9-96-10

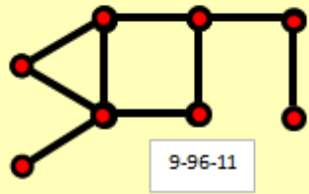
9-96-11

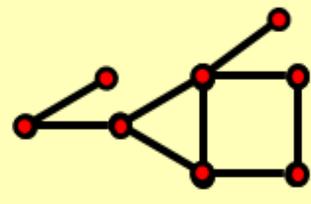
9-96-12

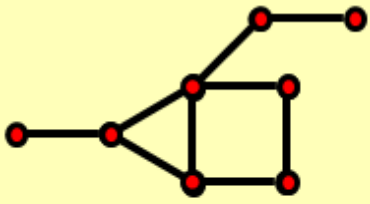
9-96-13

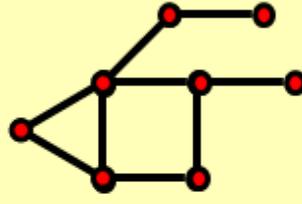
9-96-14

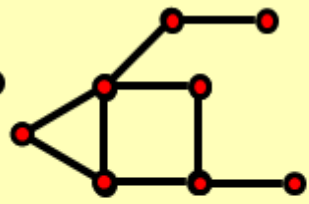
9-96-15

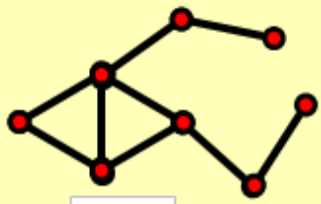
9-96-16

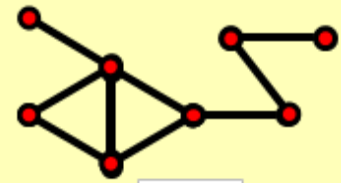
9-96-17

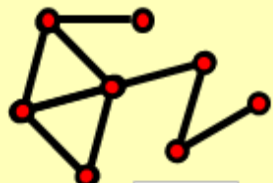
9-96-18

9-96 — 18

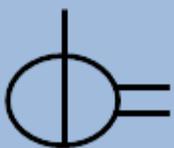
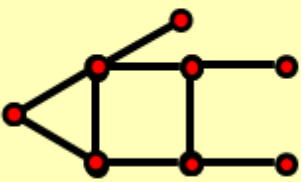
9-97-1

9-97 — 1

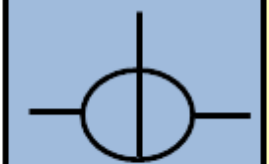
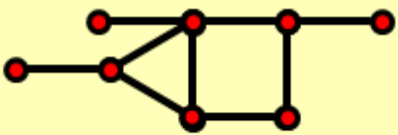
9-98-1

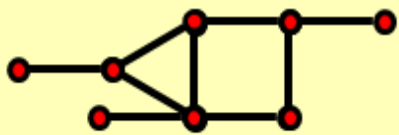
9-98-2

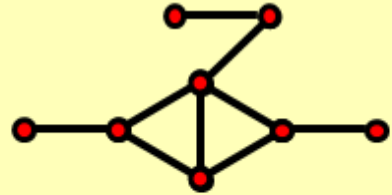
9-98-3

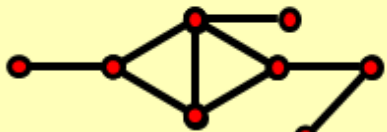
9-98-4

9-98 — 4

| 9-99 | 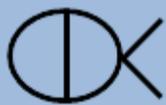 | 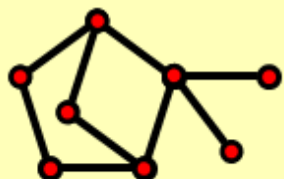 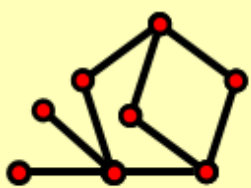 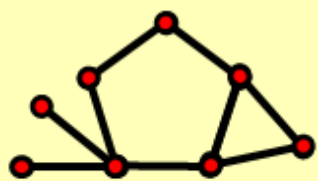 9-99-1  9-99-2  9-99-3 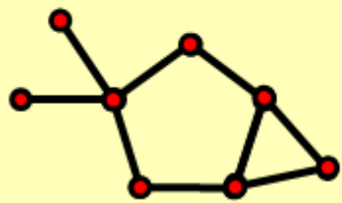 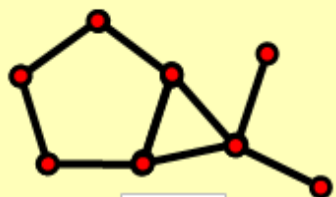 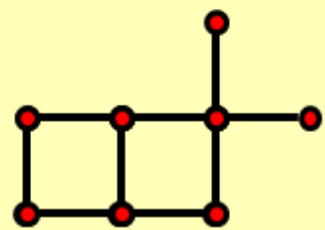 9-99-4  9-99-5  9-99-6 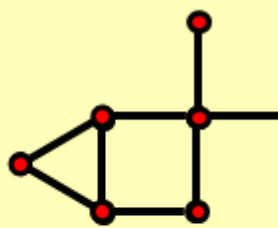 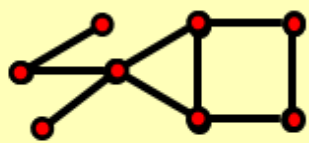 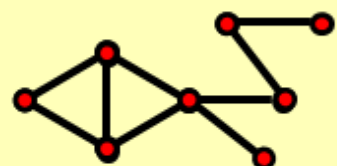 9-99-7  9-99-8  9-99-9 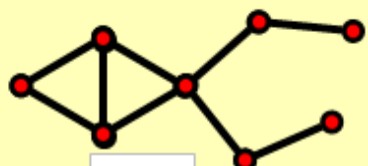 9-99-10 | 10 |
|---|---|---|---|
| 9-100 | 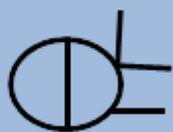 | 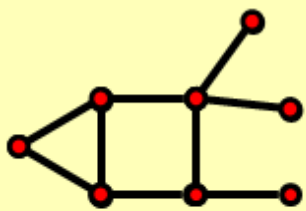 9-100-1 | 1 |
| 9-101 | 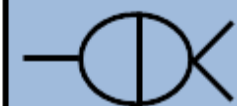 | 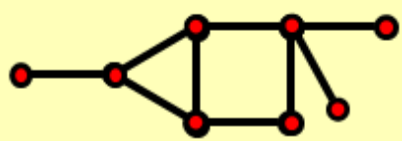 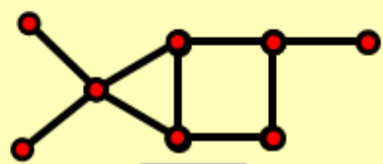 9-101-1  9-101-2 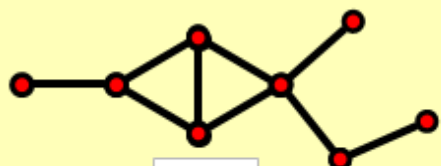 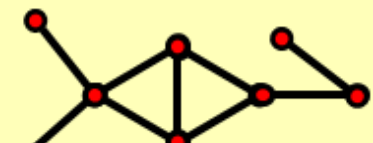 9-101-3  9-101-4 | 4 |
| 9-102 | 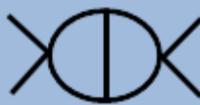 | 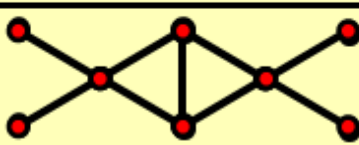 9-102-1 | 1 |

| | | |
|---|---|---|
| 9-103 | 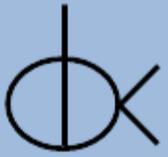 | 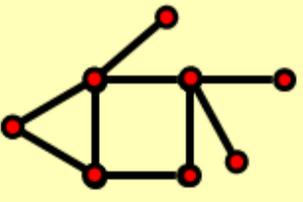 9-103-1  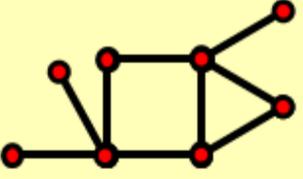 9-103-2  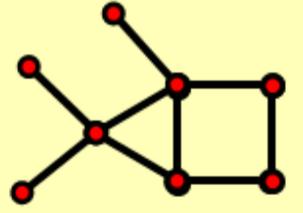 9-103-3 <br> 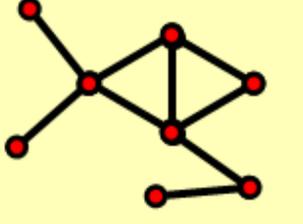 9-103-4  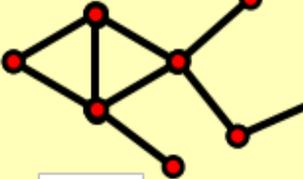 9-103-5    5 |
| 9-104 | 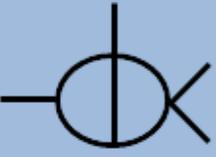 | 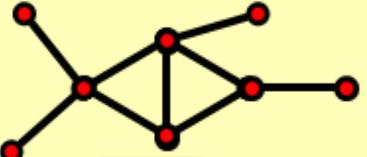 9-104-1    1 |
| 9-105 | 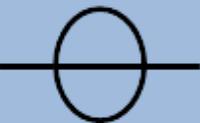 | 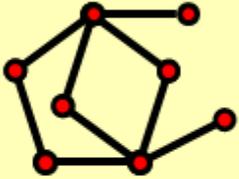 9-105-1  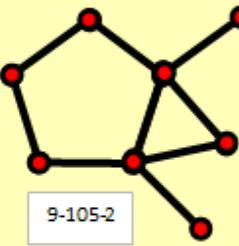 9-105-2  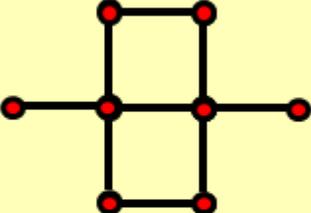 9-105-3  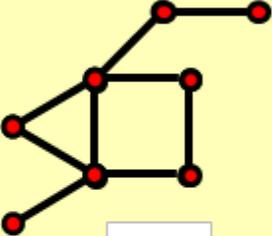 9-105-4 <br> 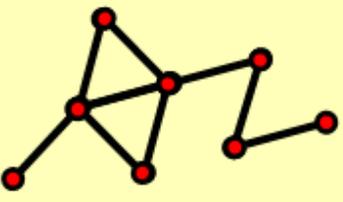 9-105-5  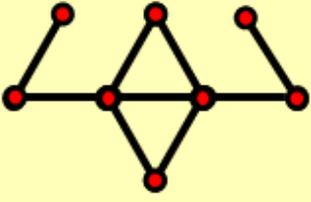 9-105-6    6 |
| 9-106 | 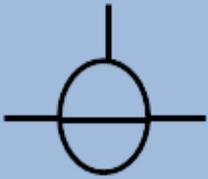 | 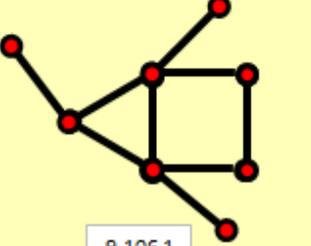 9-106-1  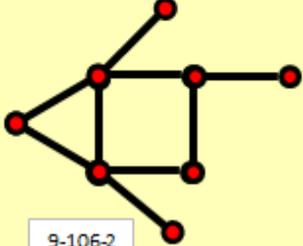 9-106-2  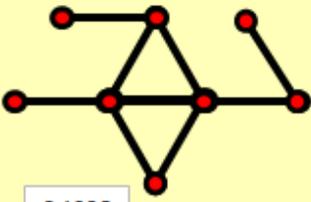 9-106-3 <br> 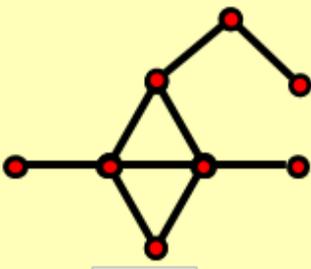 9-106-4    4 |

| | | | | |
|---|---|---|---|---|
| 9-107 | 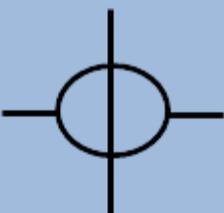 | 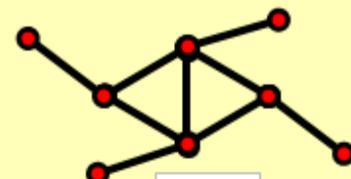 9-107-1 | | 1 |
| 9-108 | 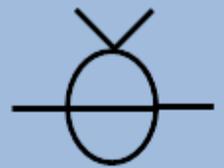 | 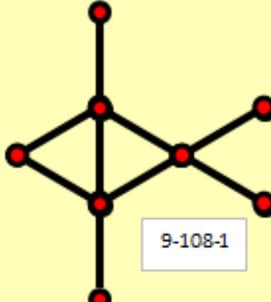 9-108-1 | | 1 |
| 9-109 | 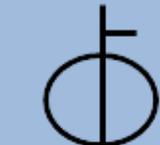 | 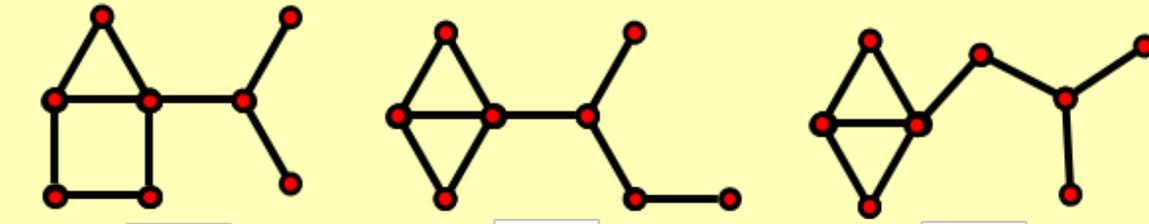 9-109-1    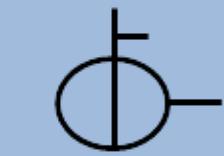 9-109-2    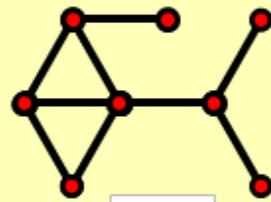 9-109-3 | | 3 |
| 9-110 | 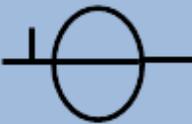 | 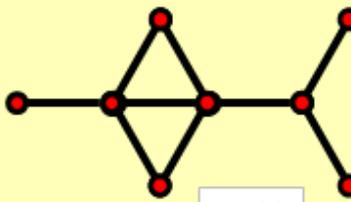 9-110-1 | | 1 |
| 9-111 | 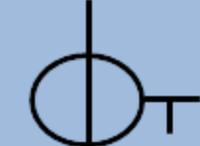 | 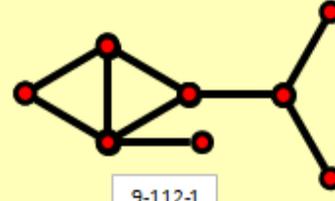 9-111-1 | | 1 |
| 9-112 | 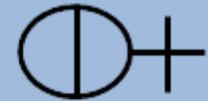 | 9-112-1 | | 1 |
| 9-113 | 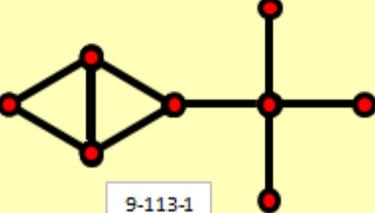 | 9-113-1 | | 1 |

| | | |
|---|---|---|
| 9-114 | 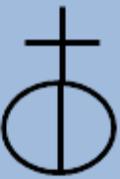 | 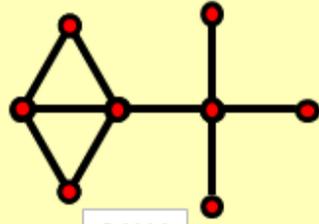<br>9-114-1 |
| 9-115 | 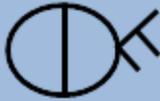 | 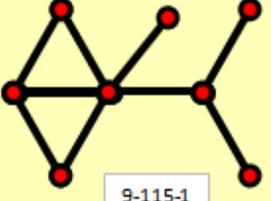<br>9-115-1 |
| 9-116 | 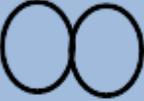 | 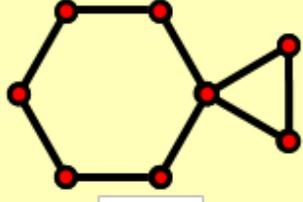 9-116-1    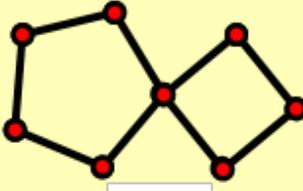 9-116-2 |
| 9-117 | 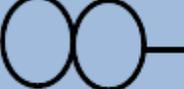 | 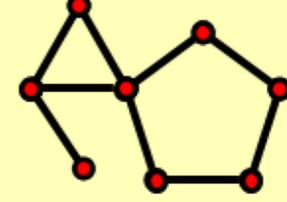 9-117-1   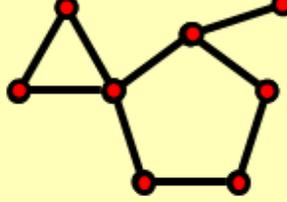 9-117-2   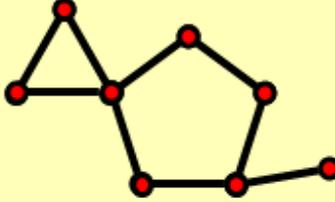 9-117-3 <br> 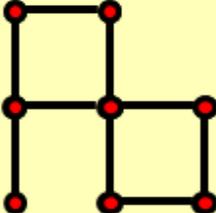 9-117-4   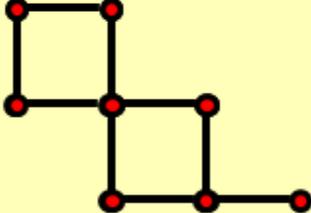 9-117-5   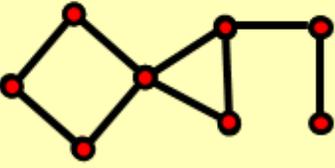 9-117-6 <br> 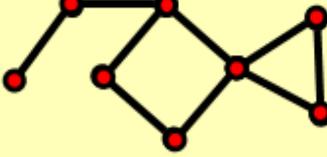 9-117-7   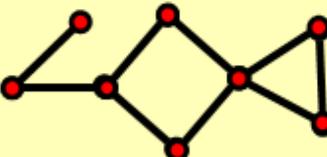 9-117-8   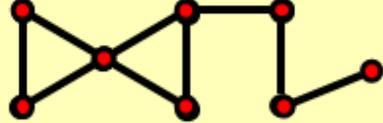 9-117-9 |
| 9-118 | 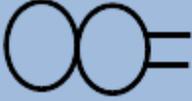 | 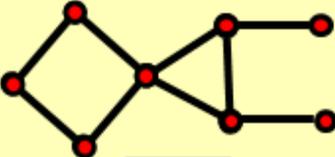 9-118-1   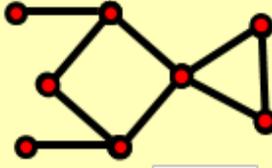 9-118-2   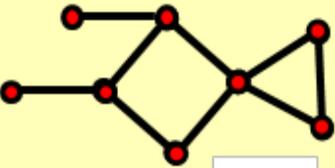 9-118-3 <br> 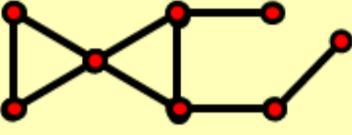 9-118-4 |
| 9-119 | 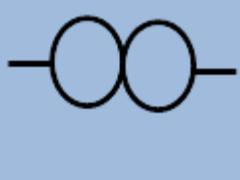 | 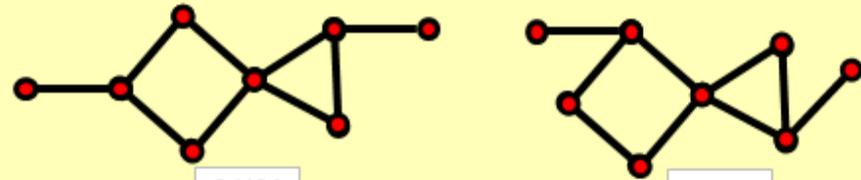 3 |
| 9-120 | 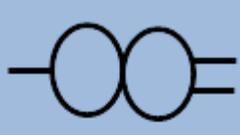 | 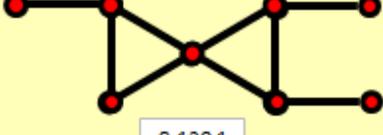 1 |
| 9-121 | 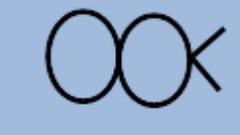 | 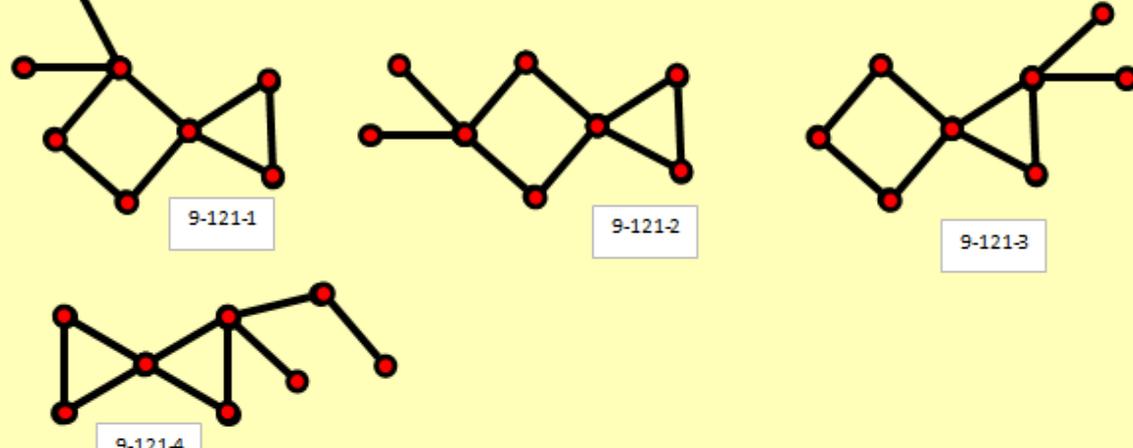 4 |
| 9-122 | 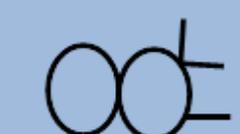 | 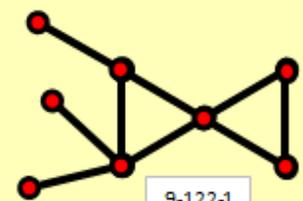 1 |
| 9-123 | 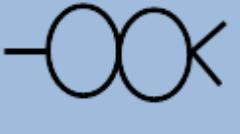 | 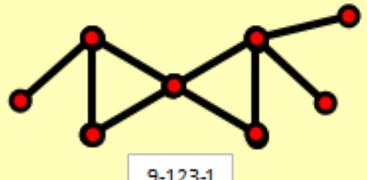 1 |
| 9-124 | 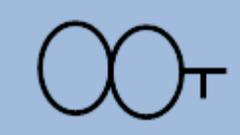 | 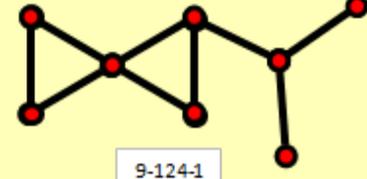 1 |
| 9-125 | 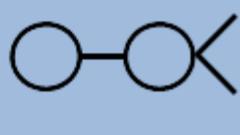 | 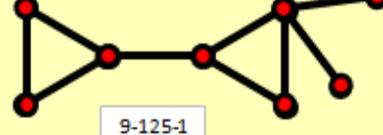 1 |

| 9-126 | 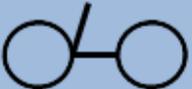 | 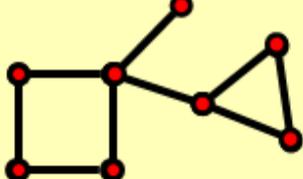 9-126-1 | 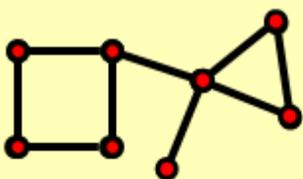 9-126-2 | 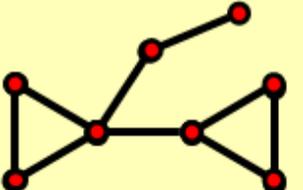 9-126-3 |
| | | 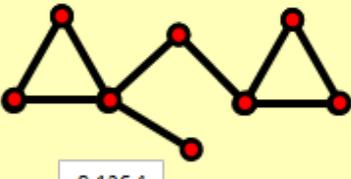 9-126-4 | | 4 |

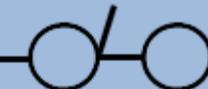

| 9-127 | 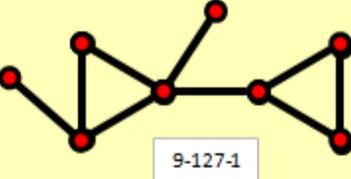 | 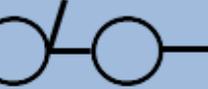 9-127-1 | | 1 |

| 9-128 | 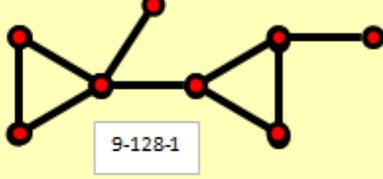 | 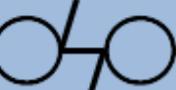 9-128-1 | | 1 |

| 9-129 | 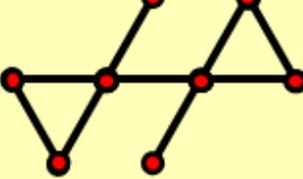 | 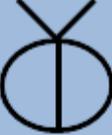 9-129-1 | | 1 |

### 3.9.3.3. Matchstick graphs with |E|=9, $F$=3, |V|=8, $\Delta$=5

| 9-130 | 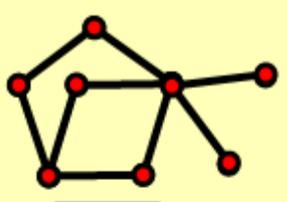 | 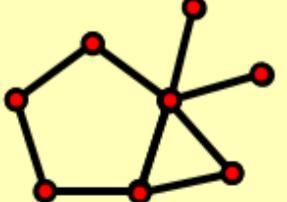 9-130-1 | 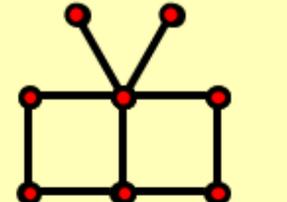 9-130-2 | 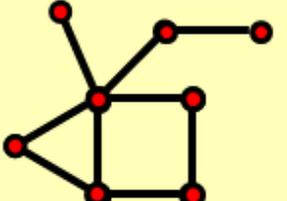 9-130-3 |
| | | 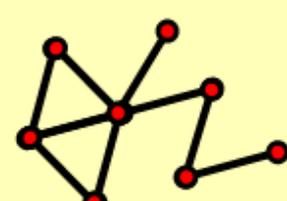 9-130-4 | 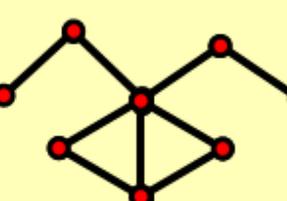 9-130-5 | 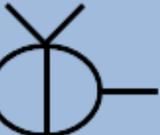 9-130-6 |
| | | | | 6 |

| 9-131 | 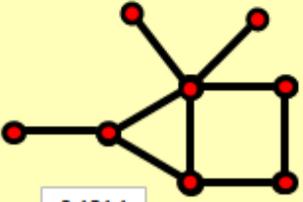 | 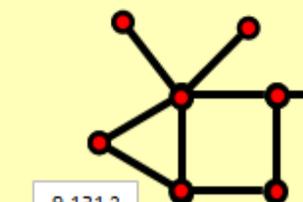 9-131-1 | 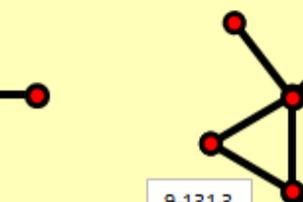 9-131-2 | 9-131-3 |

| | | |
|---|---|---|
| 9-131 | 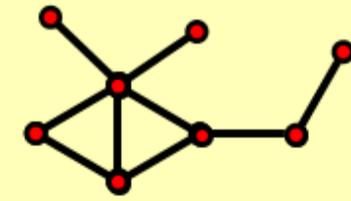 | 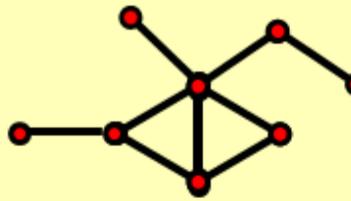 9-131-4  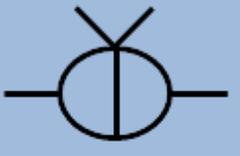 9-131-5   5 |
| 9-132 | 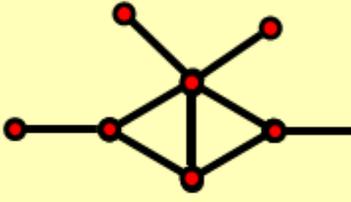 | 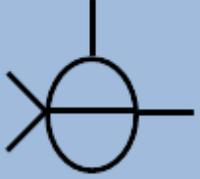 9-132-1   1 |
| 9-133 | 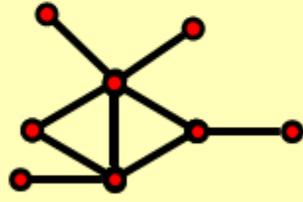 | 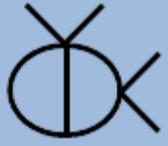 9-133-1   1 |
| 9-134 | 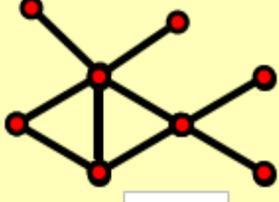 | 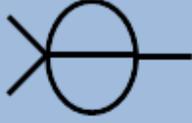 9-134-1   1 |
| 9-135 | 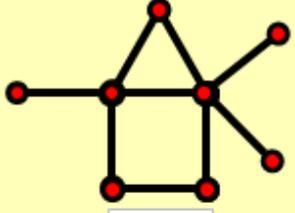 | 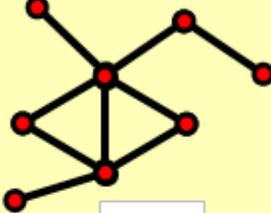 9-135-1  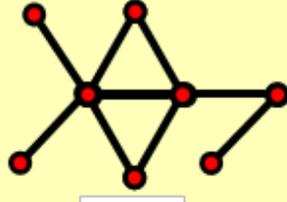 9-135-2  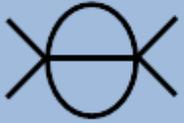 9-135-3   3 |
| 9-136 | 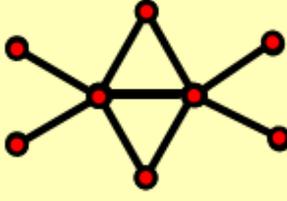 | 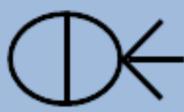 9-136-1   1 |
| 9-137 | 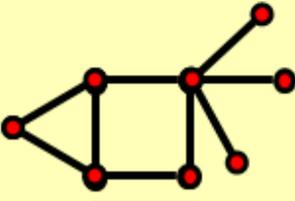 | 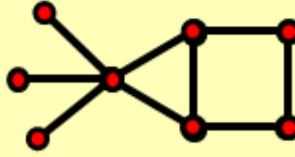 9-137-1  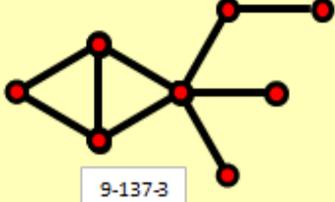 9-137-2   9-137-3   3 |
| 9-138 | 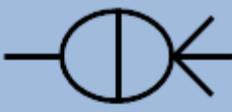 | 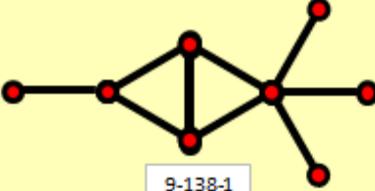<br>9-138-1   1 |
| 9-139 | 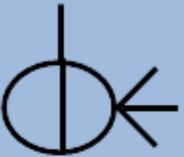 | 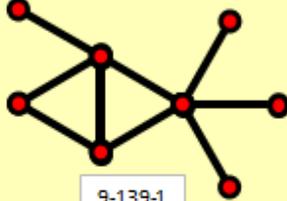<br>9-139-1   1 |
| 9-140 | 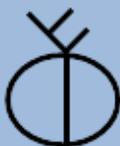 | 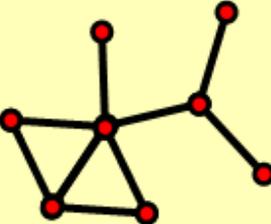<br>9-140-1   1 |
| 9-141 | 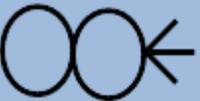 | 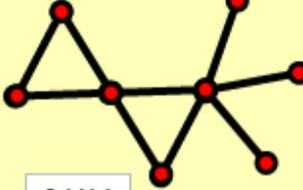<br>9-141-1   1 |
| 9-142 | 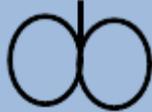 | 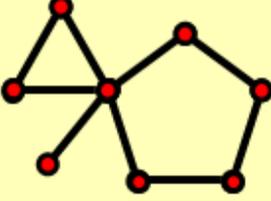 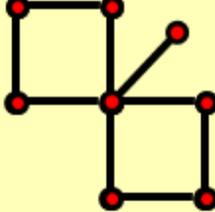 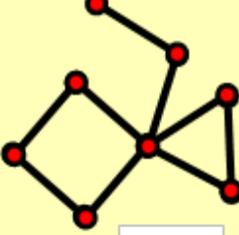<br>9-142-1   9-142-2   9-142-3<br>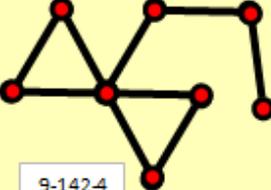<br>9-142-4   4 |
| 9-143 | 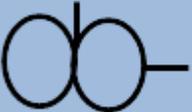 | 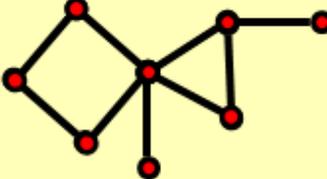 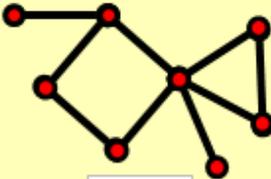 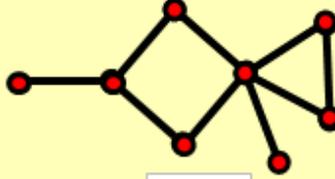<br>9-143-1   9-143-2   9-143-3<br>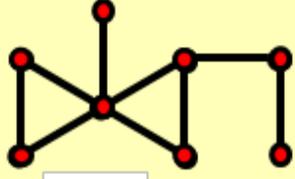 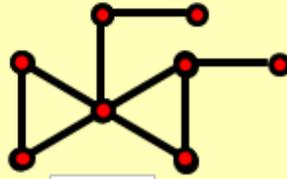<br>9-143-4   9-143-5   5 |

| | | | |
|---|---|---|---|
| 9-144 | (graph) | (matchstick graph) 9-144-1 | 1 |
| 9-145 | (graph) | (matchstick graph) 9-145-1 | 1 |
| 9-146 | (graph) | (matchstick graph) 9-146-1 | 1 |
| 9-147 | (graph) | (matchstick graph) 9-147-1 | 1 |
| 9-148 | (graph) | (matchstick graph) 9-148-1 | 1 |

3.9.3.4. Matchstick graphs with |E|=9, $F$=3, |V|=8, $\Delta$=6

| | | | |
|---|---|---|---|
| 9-149 | (graph) | (matchstick graphs) 9-149-1, 9-149-2 | 2 |
| 9-150 | (graph) | (matchstick graph) 9-150-1 | 1 |

| 9-151 | 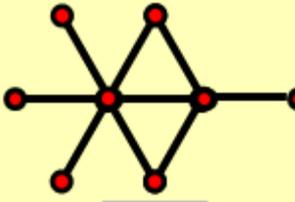 | 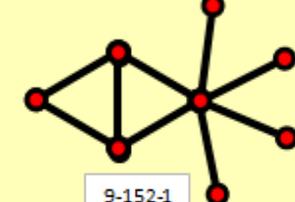  9-151-1 | 1 |
| 9-152 | 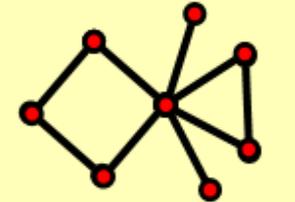 | 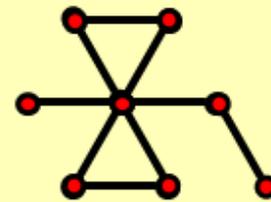  9-152-1 | 1 |
| 9-153 | 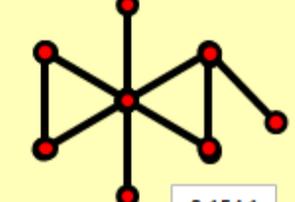 | 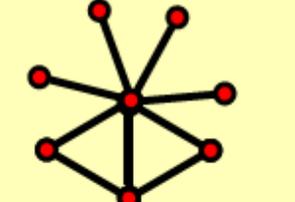  9-153-1  9-153-2 | 2 |
| 9-154 | 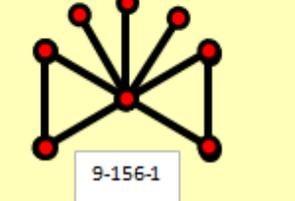 | 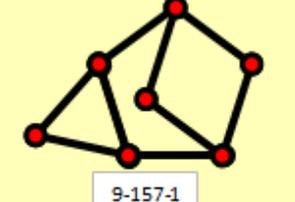  9-154-1 | 1 |

### 3.9.3.5. Matchstick graphs with |E|=9, F=3, |V|=8, Δ=7

| 9-155 | 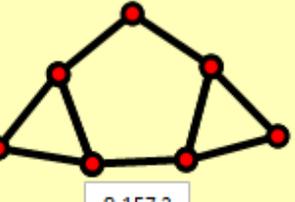 | 9-155-1 | 1 |
| 9-156 | | 9-156-1 | 1 |

### 3.9.4. Matchstick graphs with |E|=9, F=4, |V|=7
### 3.9.4.1. Matchstick graphs with |E|=9, F=4, |V|=7, Δ=3

| 9-157 | 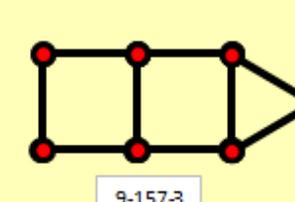 | 9-157-1  9-157-2  9-157-3 | 3 |

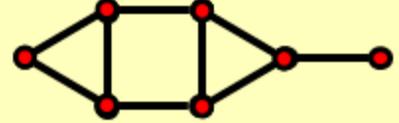

3.9.4.2. Matchstick graphs with |E|=9, $F$=4, |V|=7, $\Delta$=4

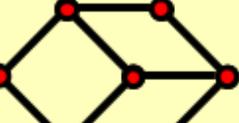

| | | |
|---|---|---|
| 9-164 | 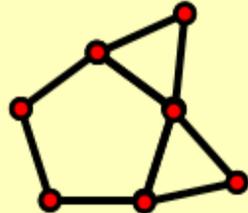 | 9-164-4, 9-164-5, 9-164-6, 9-164-7 — 7 |
| 9-165 | 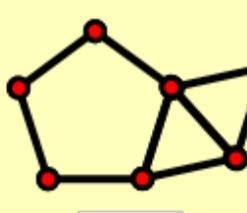 | 9-165-1, 9-165-2, 9-165-3, 9-165-4, 9-165-5 — 5 |
| 9-166 | 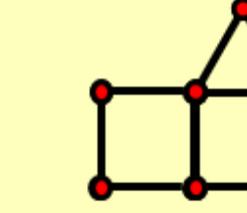 | 9-166-1 — 1 |
| 9-167 | 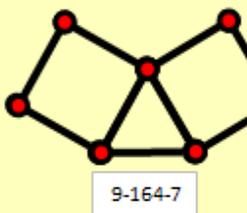 | 9-167-1 — 1 |
| 9-168 | 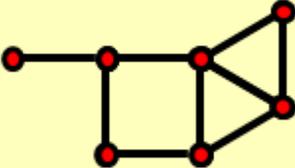 | 9-168-1, 9-168-2, 9-168-3, 9-168-4 — 4 |

| | | | | |
|---|---|---|---|---|
| 9-169 | 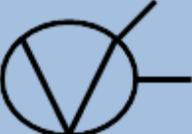 | 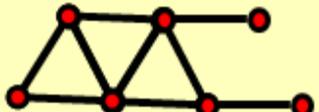 9-169-1 | | 1 |
| 9-170 | 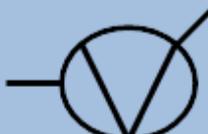 | 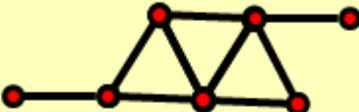 9-170-1 | | 1 |
| 9-171 | 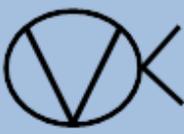 | 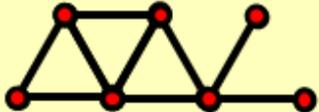 9-171-1 | | 1 |
| 9-172 | 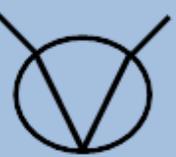 | 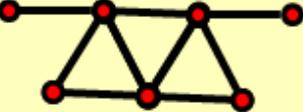 9-172-1 | | 1 |
| 9-173 | 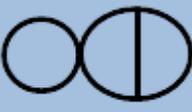 | 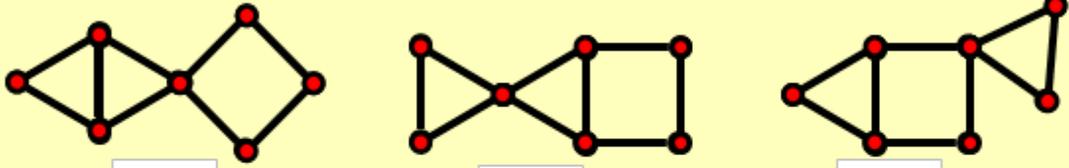 9-173-1   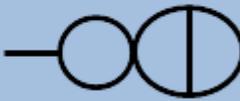 9-173-2   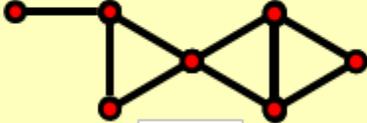 9-173-3 | | 3 |
| 9-174 | 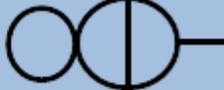 | 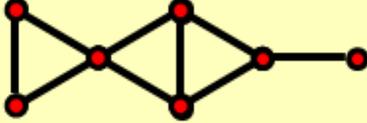 9-174-1 | | 1 |
| 9-175 | 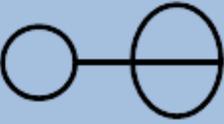 | 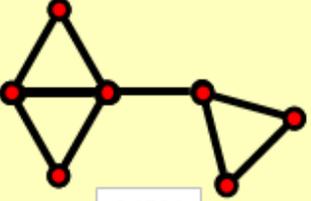 9-175-1 | | 1 |
| 9-176 | | | | 1 |

| 9-177 | 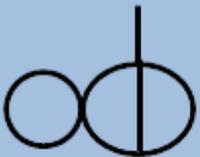 | 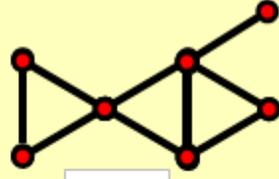9-177-1 | 1 |
| 9-178 | 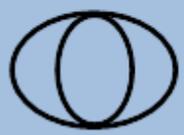 | 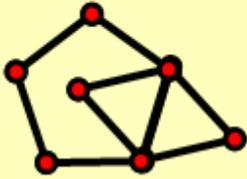9-178-1　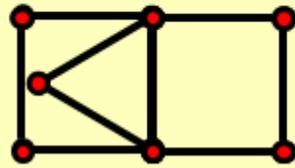9-178-2 | 2 |
| 9-179 | 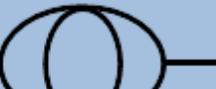 | 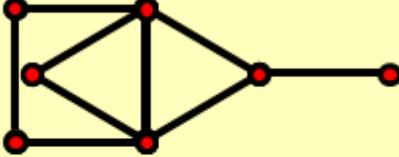9-179-1　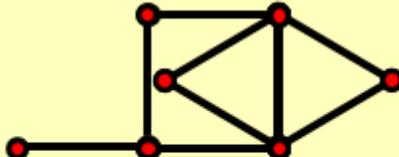9-179-2 | 2 |
| 9-180 | 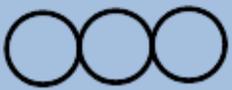 | 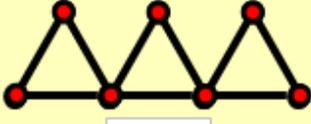9-180-1 | 1 |

3.9.4.3. Matchstick graphs with $|E|=9$, $F=4$, $|V|=7$, $\Delta=5$

| 9-181 | 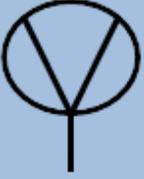 | 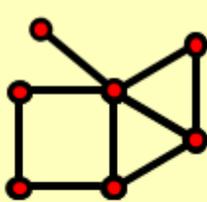9-181-1　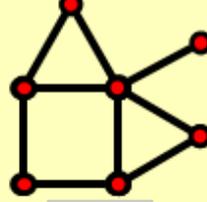9-181-2　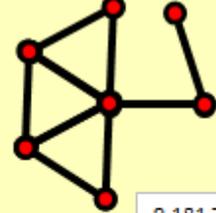9-181-3 | 3 |
| 9-182 | 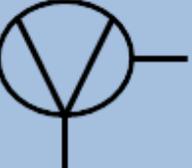 | 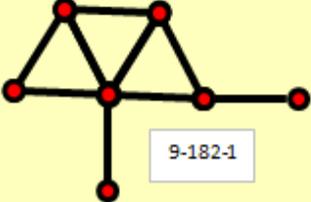9-182-1 | 1 |
| 9-183 | 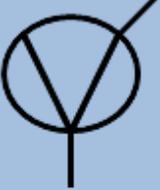 | 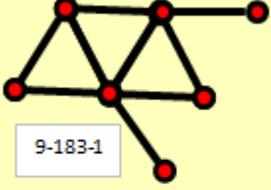9-183-1 | 1 |
| 9-184 | 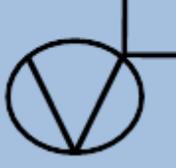 | 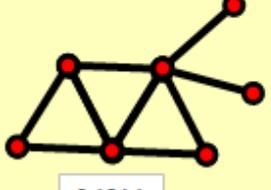9-184-1 | 1 |

| | | |
|---|---|---|
| 9-185 | 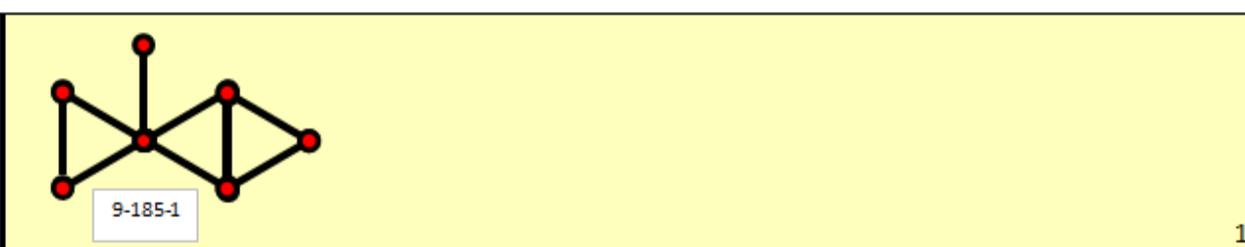 | 1 |
| 9-186 | 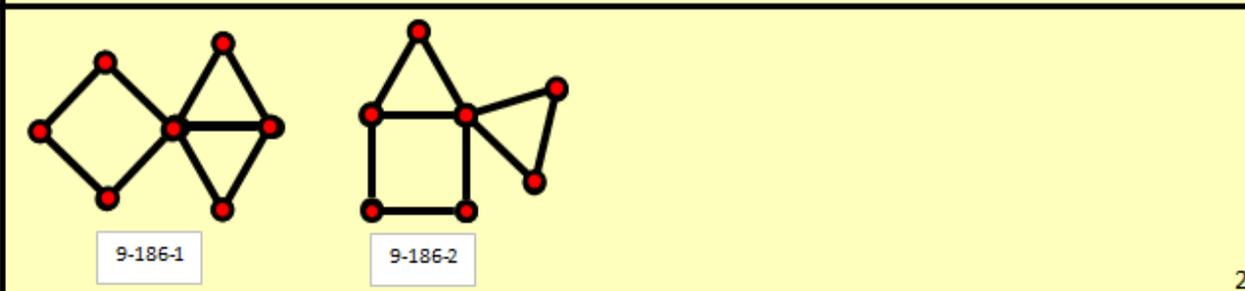 | 2 |
| 9-187 | 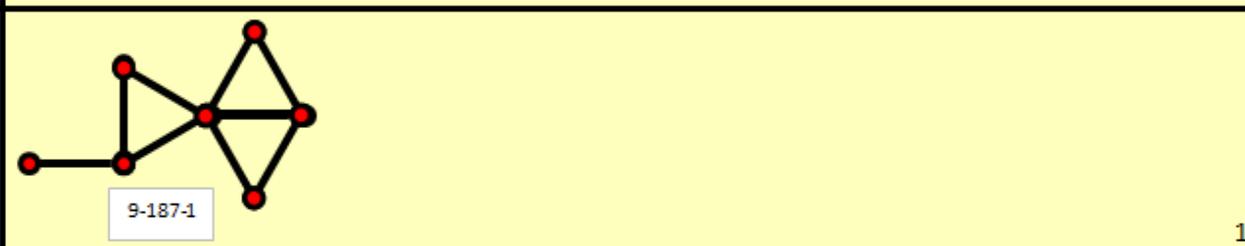 | 1 |
| 9-188 | 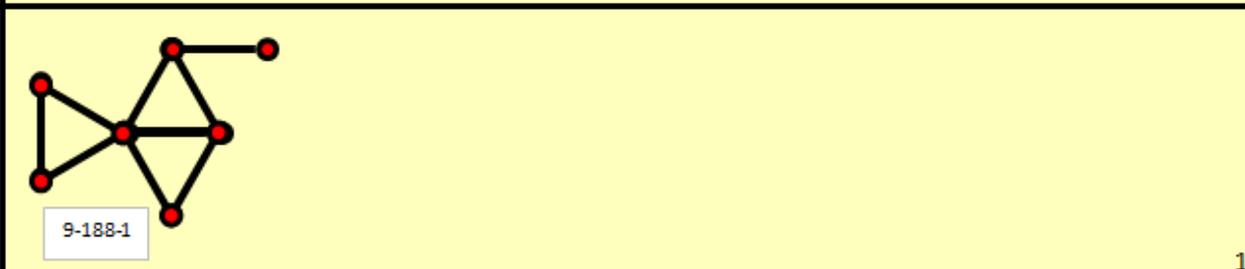 | 1 |
| 9-189 | 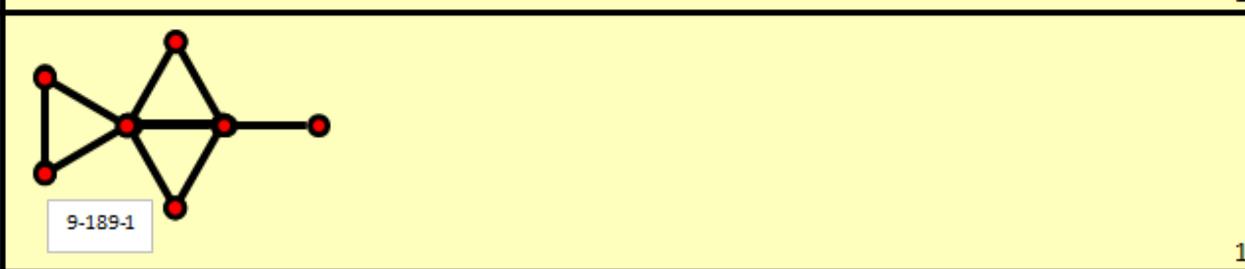 | 1 |
| 9-190 | 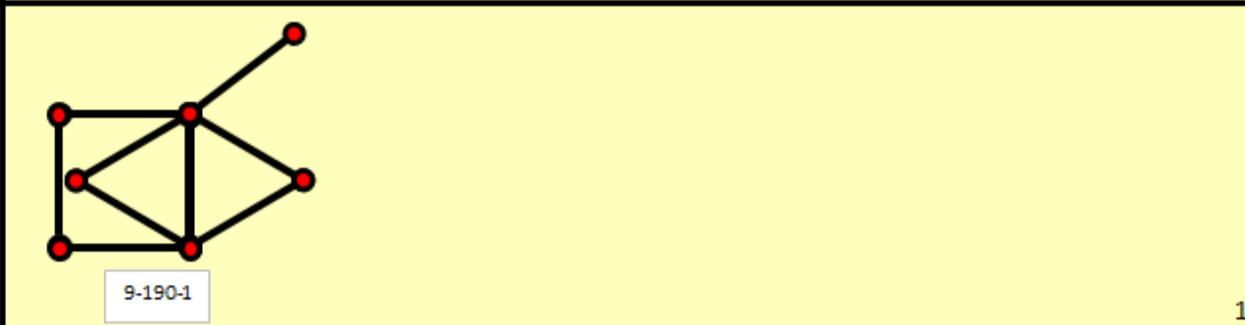 | 1 |

3.9.4.4. Matchstick graphs with |E|=9, $F$=4, |V|=7, $\Delta$=6

| | | |
|---|---|---|
| 9-191 | 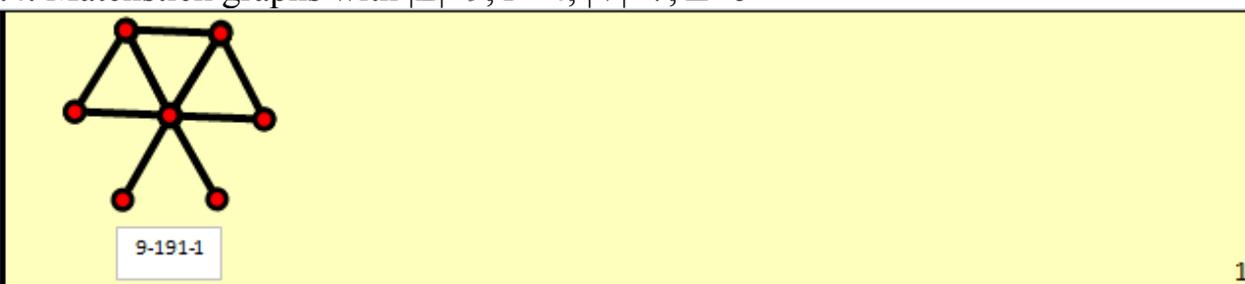 | 1 |
| 9-192 | 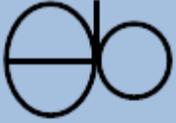 | 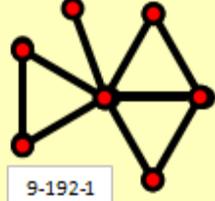<br>9-192-1 |
| 9-193 | 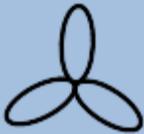 | 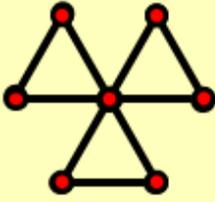<br>9-193-1 |

### 3.9.5. Matchstick graphs with |E|=9, $F$=5, |V|=6

#### 3.9.5.1. Matchstick graphs with |E|=9, $F$=5, |V|=6, $\Delta$=4

| | | |
|---|---|---|
| 9-194 | 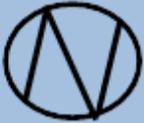 | 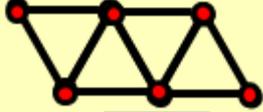<br>9-194-1 |
| 9-195 | 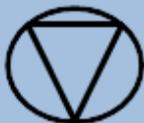 | 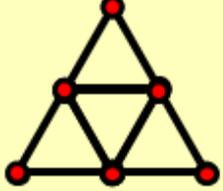<br>9-195-1 |

#### 3.9.5.2. Matchstick graphs with |E|=9, $F$=5, |V|=6, $\Delta$=5

| | | |
|---|---|---|
| 9-196 | 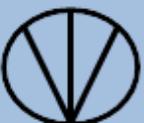 | 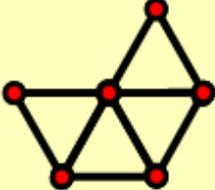<br>9-196-1 |

# 4. Plots

Plots 1 and 2 displays the value of ϙ (n), respectively in linear and logarithmic scale.

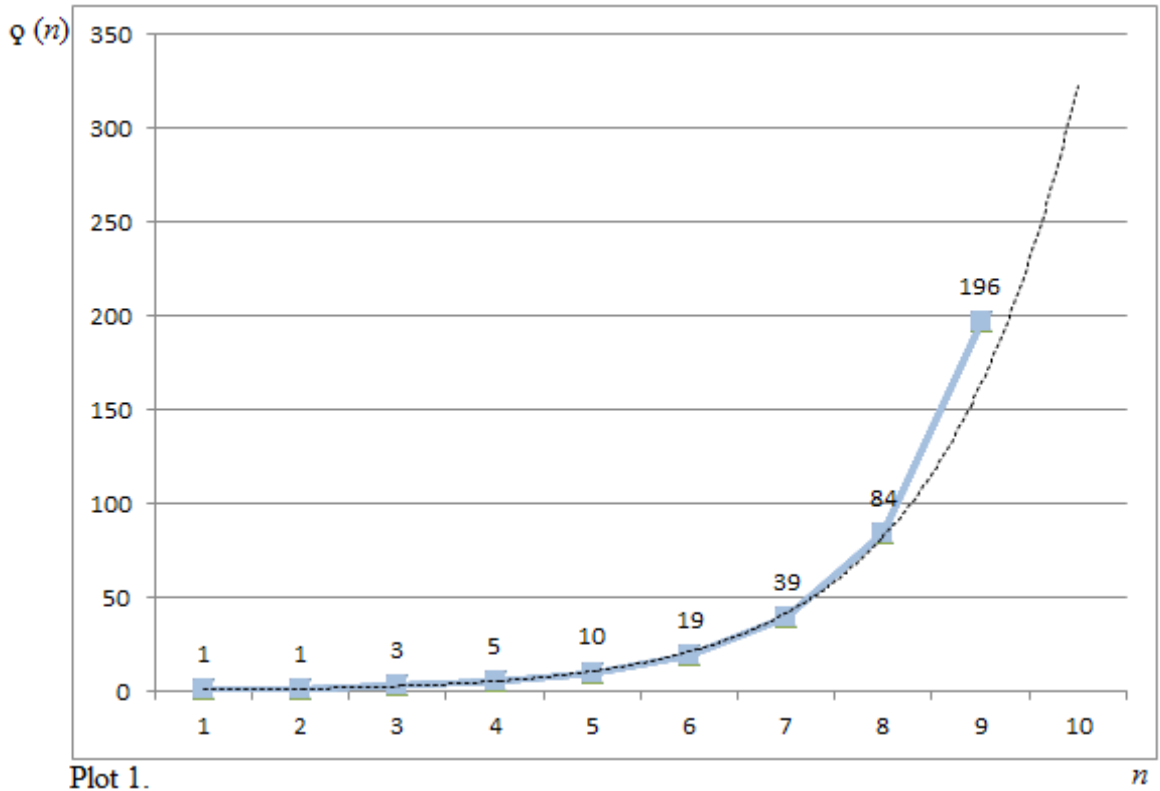

Plot 1.

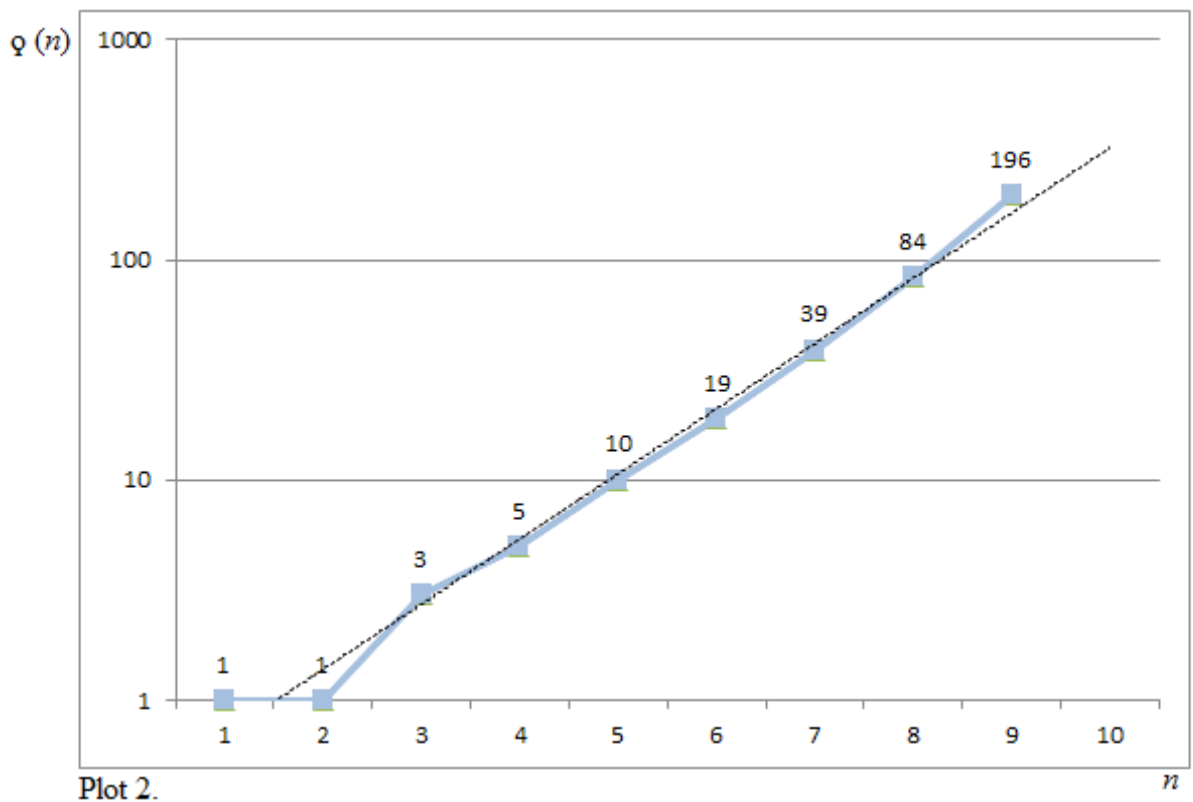

Plot 2.

The same is done for þ (n) in plots 3 and 4.

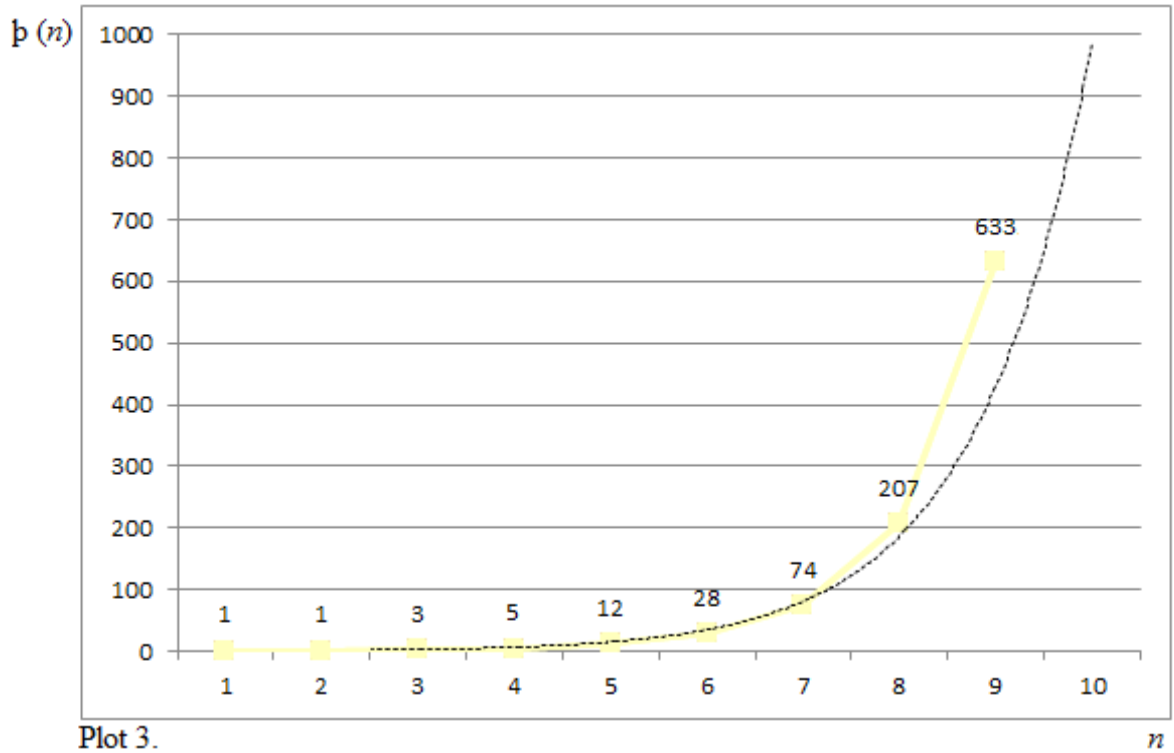

Plot 3.

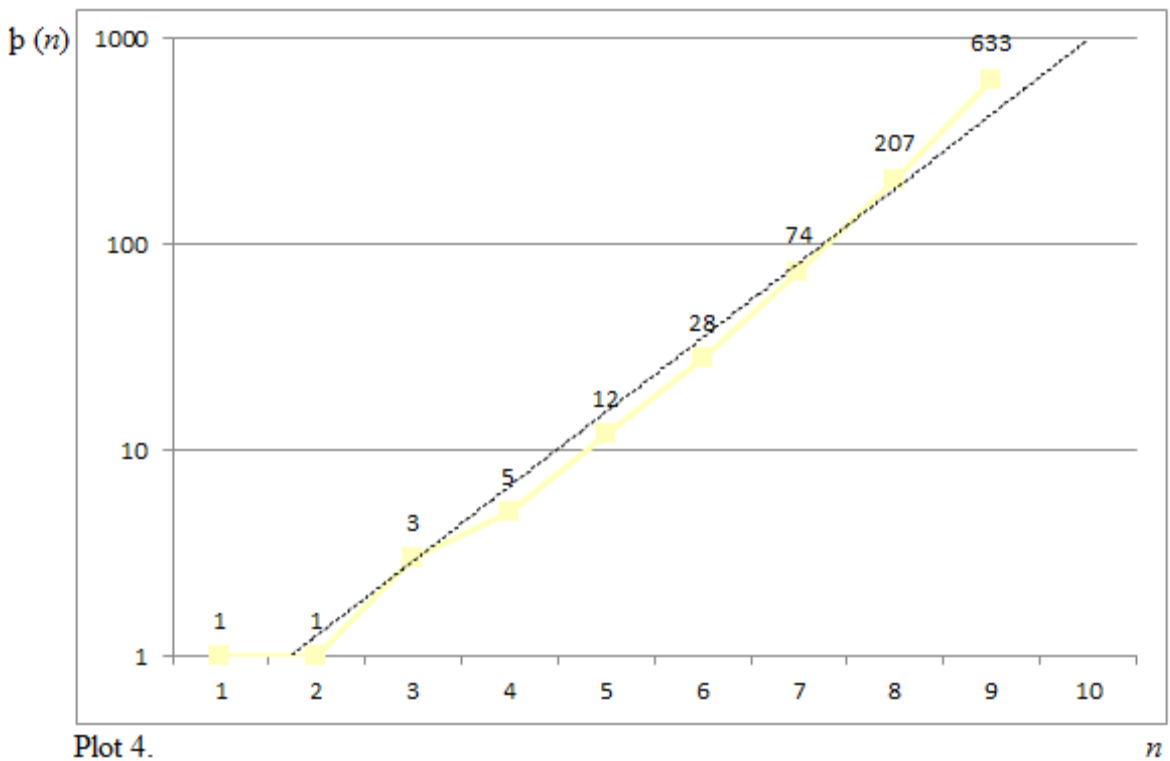

Plot 4.

The incremental ratios between consecutive terms of ϙ(*n*) and þ(*n*) are shown in plots 5 and 6. After an initial oscillation, with *n*>4 both quotients grow monotonically as *n* increases.

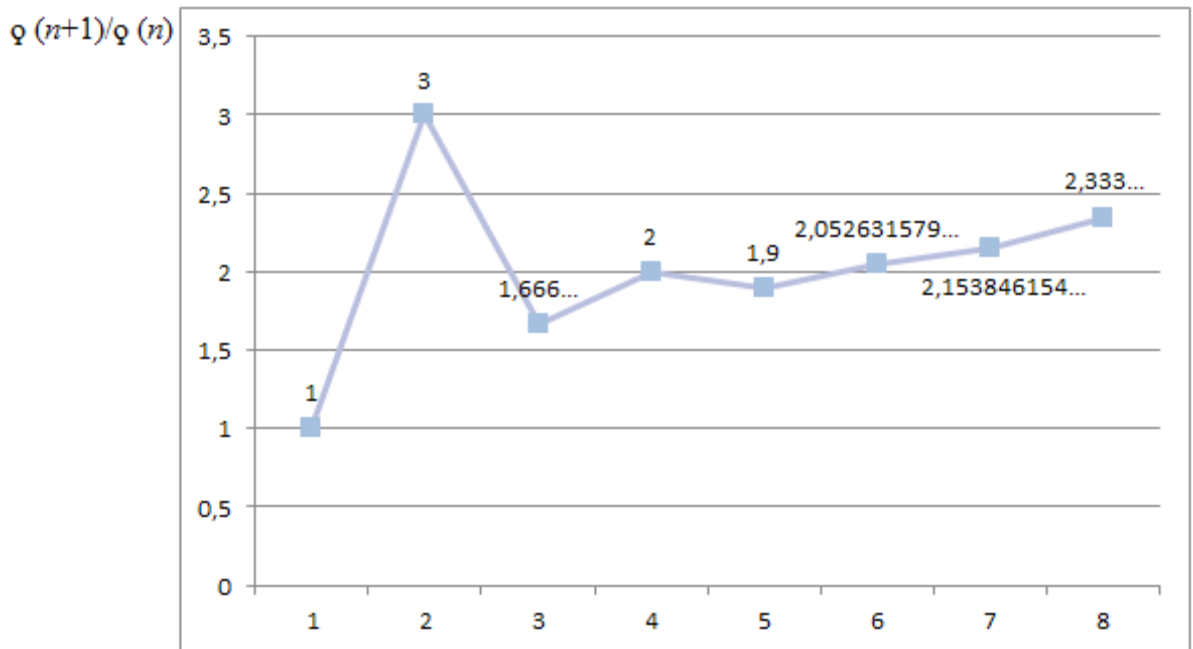
Plot 5.

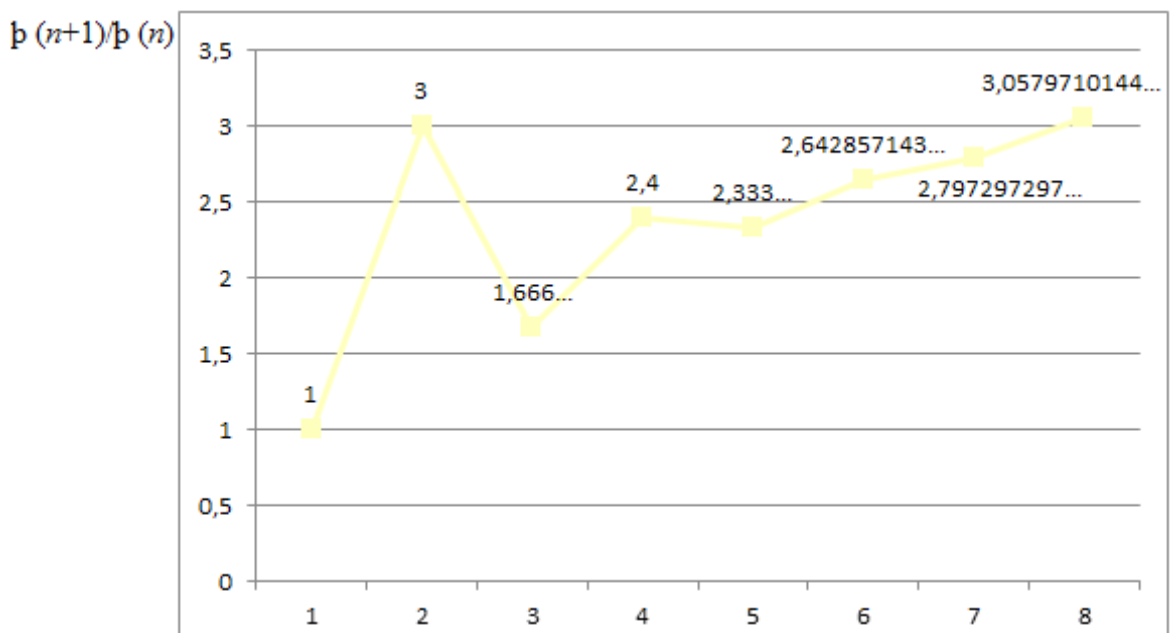
Plot 6.

If we suppose that the trend continues also for *n*=10, this would place a lower limit of 458 for ǫ(10) and of 1936 for þ(10). It seems reasonable to expect ǫ(10) to be about from 500 to 600.

Plot 7 reports the average number of non-isomorphic matchstick graphs for homeomorphism class. Also this sequence shows a monotonically growing behavior with *n*>4.

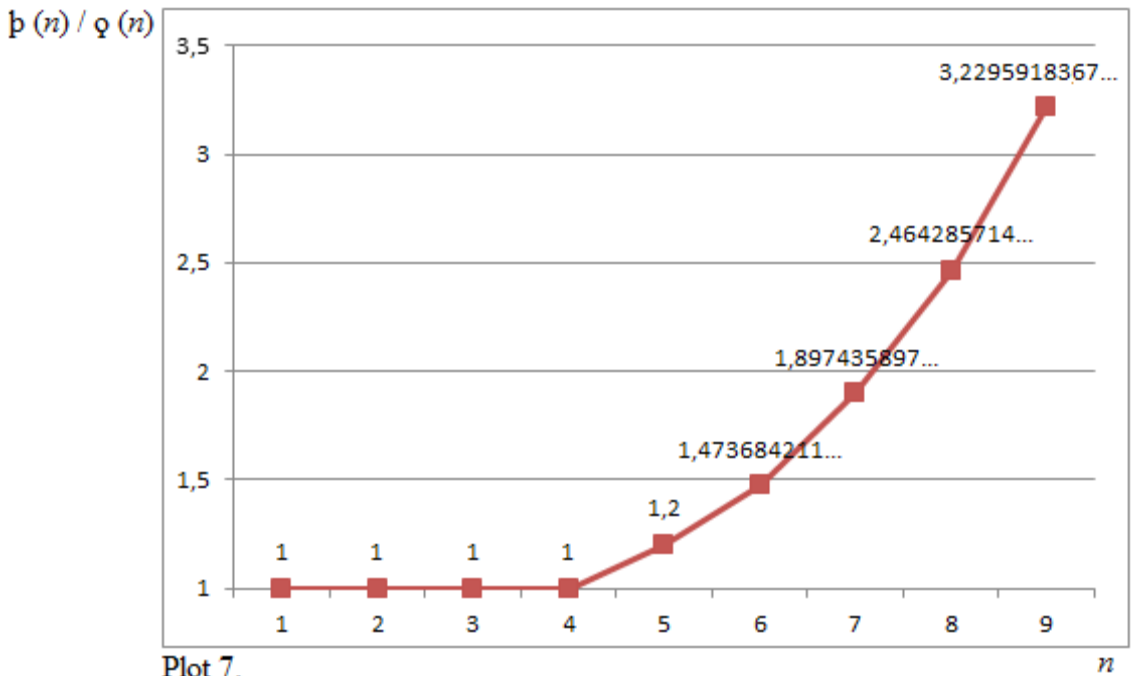
Plot 7.

Plot 8 shows what fraction of all the connected planar graphs consists in matchstick graphs. Up to 5 edges, every planar graph is a matchstick graph as well. Then the quotient starts decreasing quite irregularly. Since we know there are 2318 connected planar graphs with 10 edges (Steinbach 1990 [7]), to made a coarse intepolation, would be plausible to expect a value of about 2000-2100 for þ (*n*).

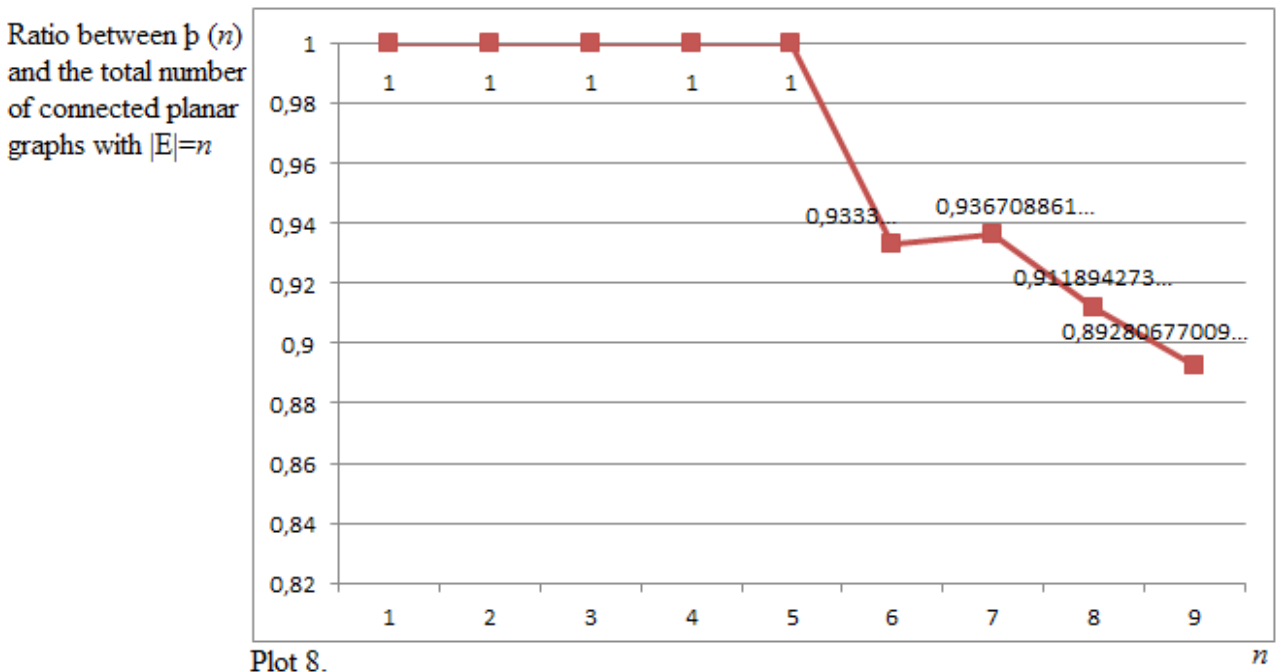
Plot 8.

# 5. Acknowledgement

I gratefully thank Alexis Vaisse (Ubisoft Entertainment) for the correction of some errors in a previous version of this catalog.

RAFFAELE SALVIA

*E-mail address:* raffaelesalvia@alice.it